\documentclass[3p,times,nopreprintline]{elsarticle}

\usepackage{etoolbox}

\usepackage{amssymb}
\usepackage{amsmath}
\usepackage{amsthm}
\usepackage{lipsum}
\usepackage{subcaption}

\usepackage{graphicx}
\usepackage{natbib}
\biboptions{sort&compress}
\usepackage{hyperref}
\usepackage{orcidlink}
\hypersetup{
	colorlinks,
	citecolor=blue,
	filecolor=blue,
	linkcolor=blue,
	urlcolor=blue
}

\usepackage{amsfonts}
\usepackage{float}
\usepackage{xcolor}
\usepackage{booktabs}
\usepackage{bm}

\theoremstyle{definition}
  
\usepackage{mathtools}

\usepackage{setspace} % For adjusting line spacing
\usepackage{parskip} % 
\usepackage{fancyhdr} %
\journal{}%

\begin{document}

\begin{frontmatter}

\title{Projection-based Reduced Order Modelling for Unsteady Parametrized Optimal Control Problems in 3D Cardiovascular Flows}

%% Author's Name
\author[1]{Surabhi Rathore}%\ead{surabhi.rathore@sissa.it}
\author{Pasquale C. Africa\corref{cor1}\fnref{1}}
\cortext[cor1]{Corresponding author.  \textit{E-mail address:} \url{pafrica@sissa.it}}
\author[2]{Francesco Ballarin}%\ead{francesco.ballarin@unicatt.it}
\author[1]{Federico Pichi}%\ead{fpichi@sissa.it}
\author[1]{Michele Girfoglio}%\ead{mgirfogl@sissa.it}
\author[1]{Gianluigi Rozza}%\ead{grozza@sissa.it}

%% Author's affiliation
\affiliation[1]{organization={mathLab, Mathematics Area},%Department and Organization
            addressline={SISSA Scuola Internazionale Superiore di Studi Avanzati, Via Bonomea 265}, 
            city={Trieste},
            postcode={34136}, 
            country={Italy}}

\affiliation[2]{organization={Department of Mathematics and Physics},%Department and Organization
            addressline={Università Cattolica del Sacro Cuore, Via Garzetta 48}, 
            city={Brescia},
            postcode={25133}, 
            country={Italy}}

%% Abstract
\begin{abstract}

\noindent \textit{Background and Objective:} Accurately defining outflow boundary conditions in patient-specific models poses significant challenges due to complex vascular morphologies, physiological conditions, and high computational demands. These challenges hinder the computation of realistic and reliable cardiovascular (CV) haemodynamics by incorporating clinical data such as 4D magnetic resonance imaging. The objective is to control the outflow boundary conditions to optimize CV haemodynamics and minimize the discrepancy between target and computed flow velocity profiles. This paper presents a projection-based reduced order modelling (ROM) framework for unsteady parametrized optimal control problems (OCP$_{(\bm{\mu})}$s)  arising from CV applications.

\noindent \textit{Methods:} Numerical solutions of OCP$_{(\mathbf{\mu})}$s require substantial computational resources, highlighting the need for robust and efficient ROMs to perform real-time and many-query simulations. We investigate the performance of a projection-based reduction technique that relies on the offline-online paradigm, enabling significant computational cost savings.  In this study, the fluid flow is governed by unsteady Navier--Stokes equations with physical parametric dependence, \textit{i.e.} the Reynolds number. The Galerkin finite element method is used to compute the high-fidelity solutions in the offline phase. We implemented a nested-proper orthogonal decomposition (\textit{nested-POD}) for fast simulation of OCP$_{(\bm{\mu})}$s that encompasses two stages: temporal compression for reducing dimensionality in time, followed by parametric-space compression on the precomputed POD modes.

\noindent \textit{Results:} %We aim to investigate the performance of a projection-based reduction technique that enables significant computational cost savings. 
We tested the efficacy of the proposed methodology on vascular models, namely an idealized bifurcation geometry and a patient-specific coronary artery bypass graft, incorporating stress control at the outflow boundary and observing consistent speed-up with respect to high-fidelity strategies.  We observed the inter-dependency between the state, adjoint, and control solutions and presented detailed flow field characteristics, providing valuable insights into factors such as atherosclerosis risk.

\noindent \textit{Conclusion:} The projection-based ROM framework provides an efficient and accurate approach for simulating parametrized CV flows. By enabling real-time, patient-specific modelling, this advancement supports personalized medical interventions and improves the predictions of disease progression in vascular regions.
\end{abstract}

%% Keywords
\begin{keyword}
 Cardiovascular Flows \sep Parametrized Partial Differential Equations \sep Optimal Control \sep Galerkin Finite Element Method \sep Lagrange Multiplier \sep Nested-Proper Orthogonal Decomposition
\vspace{3mm}
\MSC [2020] 49M41 \sep 49K20 \sep 65M60 \sep 76-10 \sep 92C50
\end{keyword}
\end{frontmatter}
%% main text

%%----------------------------------------------------%%
\section{Introduction}
\label{sec1:Intro}

In recent years, the integration of well-established computational methodologies with optimization techniques has gained popularity in the realm of computational flow control and optimization \cite{gunzburger2000,gunzburger2002}. Many advances have been made in the analysis of optimal control problems (OCPs), especially for viscous, incompressible flows governed by partial differential equations (PDEs) with applications in engineering, environmental, and biomedical fields \cite{bochev2004,hinze2008,bochev2009,quarteroni2006,gerdts2008}. These advances rely on developing models and computational algorithms to solve them. OCPs strive to identify control inputs that effectively uphold the extrema of an objective function while ensuring all constraints are met, and have been widely used in engineering to optimize the performance of complex systems. The analysis of OCPs governed by PDEs is based on the theory developed by J.~L.\ Lions \cite{lions1971,lions1972}. Parametrized optimal control problems (OCP$_{(\bm{\mu})}$s) represent a broader class of optimization problems characterized by a set of parameters \cite{lin2014,elnagar1997,teo2021,manzoni2015,nobile2024combination}, wherein the parameters, denoted by $\bm{\mu}$ $\in \mathcal{D}$, characterize some physical and/or geometrical properties of the system. The control may act as a forcing term, a boundary condition, an initial condition, or even as a coefficient in the equation. These research studies have developed advanced control strategies to provide flexibility, adaptability, and versatility, which have significant implications for various engineering applications. In particular, the works \cite{nagaiah2011,williams2019} offer an advanced perspective for comprehending system behaviour and control strategies in the realm of biomedical applications. 

Over the decades, mathematical modelling and scientific computing have proven to be vital tools in advancing our understanding of the complex physiology of cardiovascular (CV) systems and their underlying causes. The relationships between haemodynamics and CV diseases are complex and multifaceted, as discussed in \cite{malek1999,loth2008blood,kwak2014,bluestein2017}. The challenges associated with CV modelling are mainly due to: $(i)$ high computational costs; $(ii)$ complex vascular morphologies and their meshing; and $(iii)$ the setting of accurate boundary conditions, which is vital for blood flow analysis \cite{quarteroni2000,vignon2010outflow,esmaily2011comparison,rathore2021}. The studies \cite{vignon2010outflow,rathore2021,africa2024lifex} underscore the significance of outflow boundary conditions in the computational modelling of patient-specific models, including the consideration of numerical backflow. Hence, it is imperative to account for many different factors when developing numerical techniques to provide quantitative insights into blood flow dynamics within vascular regions. In this framework, \textit{coronary artery bypass grafting} (CABG) is a widely adopted surgical technique used globally to treat patients with coronary artery diseases \cite{go2014heart,sankaran2012patient,agoshkov2006shape}. This technique involves creating new pathways around blocked or narrowed coronary arteries to restore proper blood flow to the heart muscle. Numerous computational studies have been performed for CABG models  \cite{scott2000isolated,ballarin2017numerical,rosenblum2021,siena2023data}, where the outflow boundary condition is typically set as a ``do-nothing" (\textit{i.e.}, homogeneous Neumann) condition. However, incorporating the outflow boundary condition as a controlled parameter in CV modelling entails applying principles from OCPs to optimize blood flow dynamics. In \cite{rozza2005optimization,fevola2021}, the outflow boundary is estimated as a controlled parameter for a two-dimensional bifurcation vascular model and for an aortic model, respectively, using the OCP framework to optimize physiological simulations for steady Stokes flow. The intricate connection between vessel geometries and haemodynamics is vital for the analysis and can be better understood through advanced computational methods, but this requires high computational costs. Such a task becomes even more challenging if \textit{real-time} solutions of the OCP$_{(\bm{\mu})}$s are needed, spanning across the whole parameter domain $\mathcal{D} \subset \mathbb{R}^P$ with $P \geq 1$, known as the \textit{many-query} context.

To achieve this objective, we employ the reduction approaches for OCP$_{(\bm{\mu})}$s governed by PDE$_{(\bm{\mu})}$s in which the \textit{high-fidelity} numerical system of dimension, say $\mathcal{N}$, is replaced by constructing \textit{low-dimensional} problem-specific approximation spaces $(N^{rd} \ll \mathcal{N})$, discussed in \cite{rozza2012reduction,quarteroni2014reduced,lassila2014model,strazzullo2022,pichi2022}. This framework adheres to the \textit{offline-online} paradigm, where the high-fidelity solution manifold is computed using the Galerkin Finite Element (FE) method (or other classical numerical approaches) during the offline phase, requiring considerable computational resources to properly capture physical and/or geometric variability. On the other hand, the \textit{online} phase leverages the precomputed quantities to generate a low-dimensional approximation manifold, which is more efficient to query. Enhancing the computational efficiency, these reduction methodologies have gained considerable attention and researchers have delved into various model reduction strategies, including certified reduced basis \cite{hesthaven2016certified,hesthaven2022reduced}, proper orthogonal decomposition (POD)  \cite{quarteroni2014reduced,baiges2014reduced}, and non-intrusive approaches\footnote{Here, we mean that the physics is not known or the solver is a black box model, from which we cannot access the equation's operators and/or their discretization.}\cite{rozza2022advanced,hesthaven2018non,PichiArtificialNeuralNetwork2023}. These approaches are often applied to engineering and CV applications involving PDE$_{(\bm{\mu})}$s \cite{ito2001,rozza2008,manzoni2012,dede2012,strazzullo2022,BallarinChapter9Spacetime2022,zainib2021,choi2020gradient,balzotti2022data,prusak2023}, to provide very efficient solutions. However, projection-based intrusive ROMs, exploiting the knowledge of the physical model, offer significantly improved reliability and accuracy in computing the reduced-order solutions. These methods project the high-dimensional system onto a lower-dimensional subspace, typically spanned by dominant modes computed from the high-fidelity solutions; they effectively capture the essential properties of the system's dynamics while significantly reducing the computational complexity of the online phase. In \cite{ballarin2015supremizer}, such intrusive strategy has been used to investigate potential sources of pressure instabilities, exploiting the supremizer stabilization to prevent spurious pressure modes in the POD approximation of parametrized flows. By leveraging these techniques, the authors have significantly improved the efficiency and accuracy of their computational analyses of vascular models, as detailed in \cite{zainib2021}. In \cite{strazzullo2020}, OCP$_{(\bm{\mu})}$s using classical POD techniques have been used to understand the flow characteristics of a two-dimensional bifurcation model governed by the unsteady Stokes equations. The authors of \cite{ballarin2016fast} introduced an intrusive reduction technique for computing haemodynamics in CABGs using the Navier--Stokes (N--S) equations, which employs an advanced POD technique to extract dominant flow modes from high-fidelity solutions.

This research work contributes to the realm of real-time CV modelling through a novel approach for simulating intricate haemodynamics,  including the following notable advancements:
\begin{enumerate}
    \item \textbf{OCP$_{(\bm{\mu})}$ for Haemodynamic Optimization:} Developing a comprehensive framework to optimize blood flow dynamics by estimating outflow boundary conditions, treated as control variables that minimize discrepancies of the flow fields from a target velocity profile. This entails formulating OCP$_{(\bm{\mu})}$ that integrates the unsteady N--S equations along with boundary conditions, as well as assessing the influence of physical parameters, such as Reynolds number, and computing the high-fidelity solutions with varying inlet flow profiles. 

    \item \textbf{Efficient Reduced Order Modelling:}  
   Given the substantial computational cost and iterative procedures involved in solving OCP$_{(\bm{\mu})}$, we have implemented a projection-based reduction strategy to address these challenges. This methodology effectively addresses the nonlinearities intrinsic to the unsteady N--S equations due to advection velocity. This framework leverages dimensional reduction techniques, particularly the nested Proper Orthogonal Decomposition \textit{(nested-POD)}, to enhance computational efficiency while maintaining high accuracy. This advancement facilitates real-time and many-query simulations of complex, time-dependent blood flow within realistic vascular geometries.  
\end{enumerate}
%This paper presents a novel approach to reducing the computational complexity of simulating spatio-temporal flow in complex three-dimensional vascular models formulated as parametrized OCPs. We implemented a projection-based reduction strategy for OCP$_{(\bm{\mu})}$, effectively addressing the nonlinearities inherent in the unsteady N--S equations. The primary goal of the proposed methodology is to optimize blood flow dynamics by estimating outflow boundary conditions, treated as control variables, that minimize variations of the flow fields from a target velocity profile. This entails formulating parametrized OCPs that integrate the unsteady N--S equations along with boundary conditions, as well as assessing the impact of various physical parameters governing the system. Given the substantial computational costs and iterative procedures involved in solving the optimal control problems, to address the computational challenges associated with solving these complex models, we strategically leverage dimensional reduction techniques, particularly the \textit{nested-POD}, to enhance the computational efficiency while keeping a high level of accuracy.

The structure of this paper is as follows: Section~\ref{sec2:Modelling} introduces the framework of OCP$_{(\bm{\mu})}$s and discusses the Lagrangian formulation. Section~\ref{sec3:Numerical} discusses the numerical methods employed, beginning with Galerkin FE formulation in Section~\ref{sec3.1:FEM}, which is pivotal for computing high-fidelity solutions, as well as presenting a projection-based ROM  for OCP$_{(\bm{\mu})}$ in Section~\ref{sec3.2:ROM}. Moreover, Section~\ref{sec4:Results} presents a detailed analysis of the numerical simulations, focusing on the comparison between high-fidelity and reduced-order approximations and examining the flow field characteristics and unsteady flow dynamics within vascular models. Finally, a detailed discussion has been presented in Section~\ref{sec5:Con}.
%Finally, the conclusions are outlined in Section~\ref{sec5:Con}.

%%================================================================================================%%
\section{Mathematical Modelling}
\label{sec2:Modelling}

This section presents a mathematical modelling framework for OCP$_{(\bm{\mu})}$s, focusing on determining the optimal control strategy for dynamical systems. These problems comprise multiple parts: state equations governing the fluid flow, a cost functional to be minimized, dependency on the physical properties, and a control variable that influences fluid flow. The Lagrangian formulation exploits the so-called Lagrange multipliers to incorporate constraints into the optimization problem, yielding the necessary conditions for optimality.

%%-------------------------------------------------------%%
\subsection{Parametrized OCPs for Unsteady Navier-Stokes Equations}
\label{sec2.1:OCPs}

We aim at solving OCPs with parametric dependence, in which optimal solutions depend on parameters $\bm{\mu} \in \mathcal{D}$. The abstract formulation of OCP$_{(\bm{\mu})}$s is expressed as~\cite{gunzburger2002,negri2015}:
\begin{equation}\label{eq:2.1}
\begin{gathered}
\textit{given}\, \bm{\mu} \in \mathcal{D},\, \textit{find the optimal control}\,\, \bm{u}\left(\bm{x}, t; \bm{\mu} \right)\, \textit{and the state variable}\,\, \bm{y} \left(\bm{x}, t; \bm{\mu} \right) \\\textit{s.t.\ cost functional}\,\, \mathcal{J} \left(\bm{y} \left(\bm{x}, t; \bm{\mu} \right), \bm{u} \left(\bm{x}, t; \bm{\mu} \right); t; \bm{\mu} \right)\,\textit{is minimized subject to}\,\, \mathcal{E} \left(\bm{y} \left(\bm{x}, t; \bm{\mu} \right), \bm{u}\left(\bm{x}, t; \bm{\mu} \right); t; \bm{\mu} \right) = 0.
\end{gathered}
\end{equation}
In problem~\eqref{eq:2.1}, $\bm{\mu} \in \mathcal{D} \subset \mathbb{R}^P$ represents the physical parameter vector of the system, with dimension $P \geq 1.$ The state equation $ \displaystyle \mathcal{E} \left(\bm{y} \left(\bm{x}, t; \bm{\mu} \right),  \bm{u} \left(\bm{x}, t; \bm{\mu} \right); t; \bm{\mu} \right)$ describes how the state solution $\bm{y}\left(\bm{x}, t; \bm{\mu} \right) \in \mathbf{Y}$ evolves with space and time under the influence of the control $\bm{u} \left(\bm{x}, t; \bm{\mu} \right) \in \mathcal{U}.$ The control may correspond to a forcing term, boundary condition, or initial condition, with the spaces $\mathcal{U}$ and $ \mathbf{Y}$ representing the \textit{control space} and \textit{state space}, respectively. In this study, the dynamical system is a CV system, which is governed by the unsteady N--S equations with some specified initial and/or boundary conditions.
In these optimization problems, the observation of the state variable is often represented by a linear mapping, $z \left( \bm{y} \left(\bm{x}, t; \bm{\mu} \right) \right) = \bm{C} \bm{y} \left(\bm{x}, t; \bm{\mu} \right),$ where $\bm{C}$ is a suitable operator.  The cost functional $ \displaystyle \mathcal{J} \left(\bm{y} \left(\bm{x}, t; \bm{\mu} \right), \bm{u}\left(\bm{x}, t; \bm{\mu} \right); t; \bm{\mu} \right)$ is defined based on the observation, and potentially a desired profile, within the observation space $\left(\mathcal{Z} \supseteq \mathbf{Y} \right)$  associated with the dynamical system.

Let us consider a spatio-temporal domain $\Omega \times \left(0, T\right)$, where $\Omega \subset \mathbb{R}^{n_\mathrm{sd}}$, with $n_\mathrm{sd}=2, 3$ representing the spatial dimension, and $T \geq 0$ is the final time. The boundary of the spatial domain $\Omega$ can be expressed as a disjoint union of the boundaries, i.e.\ $\Gamma = \Gamma_\mathrm{w} \cup \Gamma_\mathrm{in} \cup \Gamma_\mathrm{out}$, where $\Gamma_\mathrm{w}$, $\Gamma_\mathrm{in}$, and $\Gamma_\mathrm{out}$, respectively, represent the wall, inlet, and outflow boundaries of the computational models. A rigid wall assumption is considered for the computation, and the blood is modelled as an incompressible Newtonian fluid governed by the unsteady N--S equations. 

For any given $\bm \mu \in \mathcal{D}$, the OCP$_{\left( \bm{\mu} \right)}$s in Problem~\eqref{eq:2.1} can be stated as\footnote{For simplicity, the parametrized space-time variables are denoted as $\bm{v} \left(\bm{x}, t; \bm{\mu} \right) = \bm{v} \left(t; \bm{\mu} \right)$ and $p \left(\bm{x}, t; \bm{\mu} \right) = p \left(t; \bm{\mu} \right)$, with similar notation applied to other variables throughout the paper.}:
\begin{align}\label{eq:2.2}
\begin{gathered}
\displaystyle
  \min_{\bm{v} \in \mathcal{V}, \bm{u} \in \mathcal{U}} \mathcal{J} \left(\bm{v} \left(t; \bm{\mu} \right), \bm{u} \left(t; \bm{\mu} \right); t; \bm{\mu} \right), \\[0.2cm]
  \text{where} \quad
 \mathcal{J}\left(\bm{v} \left(t; \bm{\mu} \right), \bm{u} \left(t; \bm{\mu} \right); t; \bm{\mu} \right) \doteq \frac{1}{2} \int_{0}^{T} \int_{\Omega} m_d \left(\bm{v} \left(t; \bm{\mu} \right), \bm{v}_d; t; \bm{\mu} \right) \, d \Omega \, dt + 
    \frac{\alpha}{2} \int_{0}^{T} \int_{\Gamma_{out }} \mathcal{R} \left(\bm{u} \left(t; \bm{\mu} \right); t; \bm{\mu} \right) \, d \Omega \, dt.
\end{gathered}
\end{align}
subject to the state problem 
\begin{align}\label{eq:2.3}
\begin{cases}
\begin{aligned}
\displaystyle  
&\frac{\partial \bm{v}\left(t; \bm{\mu} \right)}{\partial t} +\left(\bm{v} \left(t; \bm{\mu} \right) \cdot \nabla \right)\, \bm{v} \left(t; \bm{\mu} \right)  - \nu\, \Delta \bm{v} \left(t; \bm{\mu} \right) + \nabla p \left(t; \bm{\mu} \right) = \bm{f} \left(t; \bm{\mu} \right),  && \quad   \mbox{in }  \Omega\times \left(0, T \right), \\
&\nabla \cdot \bm{v}{\left(t; \bm{\mu} \right)} = 0,  && \quad  \mbox{in } \Omega\times \left(0, T\right), \\
&\bm{v}{\left( t; \bm{\mu} \right)} = \bm{v}_\mathrm{in}\left(\bm{\mu} \right),  && \quad \mbox{on } \Gamma_{in} \times \left(0, T\right), \\
& \bm{v}{\left( t; \bm{\mu} \right)} = 0,  && \quad \mbox{on } \Gamma_{w} \times \left(0, T\right),\\
&\bm{v} \left(\bm{\mu} \right)\left(\bm{x}, 0\right) = \bm{v}_0, && \quad  \mbox{on } \Omega\times \{0\}, \\
&\left(\nu \, \nabla \bm{v} \left(t; \bm{\mu} \right) - p \left(t; \bm{\mu} \right) \right) \cdot \bm{n} = \bm{u} \left(t; \bm{\mu} \right),  && \quad   \mbox{on } \Gamma_{out }\times \left(0, T\right).
\end{aligned}
\end{cases}
\end{align}
In Eq.~(\ref{eq:2.2}), the chosen bilinear form $m_d \left(\cdot, \cdot; t; \bm{\mu} \right): \mathcal{Z} \times \mathcal{Z} \rightarrow \mathbb{R}$ is symmetric, continuous, coercive over the observation space, expressed as $\displaystyle  m_d \left(\cdot, \cdot; t; \bm{\mu} \right) :=  \| \left( \bm{v}\left(t; \bm{\mu} \right) - \bm{v}_d; t; \bm{\mu} \right)\|^2,$ and quantifies the discrepancy between the state solution $\bm{v}$ and the desired solution $\bm{v}_d$ of the dynamical system. The term $\mathcal{R}\left(\cdot; t; \bm{\mu} \right): \mathcal{U} \times \mathcal{U} \rightarrow \mathbb{R}$ is symmetric, and continuous over control space, expressed as $\displaystyle  \mathcal{R}\left(\cdot, \cdot; t; \bm{\mu} \right):= \| \,\bm{u} \left(t; \bm{\mu} \right) \|^2$. This term, known as the \textit{Tikhonov regularization term}~\cite{engl1996,manzo2020tikhonov}, penalizes the magnitude of the control $\bm{u}$ to prevent it from becoming excessively large, and  $\alpha \in \left(0, 1 \right]$ is a regularization parameter balancing the two contributions. A low value of $\alpha$ emphasizes the minimizing the discrepancy between $\bm{v}$ and $\bm{v}_d$, giving less consideration to the regularization term. While, a large value of $\alpha \rightarrow 1$ focuses more on the influence of control. In Eq.~\eqref{eq:2.3}, $p$ represents the pressure, $\bm{v}$ denotes the velocity vector, and the control $\bm{u}$ is sought as an outflow boundary condition. 
More precisely, $\bm{u}$ represents a traction control that includes both viscous and pressure components projected onto the outward normal vector at the boundary $\Gamma_{out}$. The kinematic viscosity is denoted by $\nu$,  $\bm{v}_{in}$ and $\bm{v}_{0}$ represent the inlet and initial velocity profiles, respectively, $\bm{n}$ denotes the unit outward normal vector on the boundary $\Gamma$, and $\bm{f}$ is the source term.

Let us define suitable function spaces for the optimization problem \eqref{eq:2.2} - \eqref{eq:2.3}. The velocity field $ \displaystyle \bm{v}\left(\cdot, \cdot \right)$ belongs to the Hilbert space $\displaystyle \mathbf{V}  := \mathbf{L}^2\left(0, T; \left[\mathbf{H}^{1} \left(\Omega \right) \right]^{n_\mathrm{sd}} \right)$ and the pressure field $p\left(\cdot, \cdot \right)$ belongs to the Hilbert space $ \mathbf{P} := \mathbf{L}^2\left(0, T; \mathbf{L}^{2} \left(\Omega \right) \right)$. Therefore, the full state space is $ \displaystyle \mathcal{V} = \mathbf{V} \times \mathbf{P}$, which can be defined more precisely as, $\displaystyle \mathcal{V} = \Big\{\bm{y} = \left(\bm{v} ,p \right) \in L^2 \left(0, T; \mathcal{V} \right) : \frac{\partial \bm{y}}{\partial t} \in  L^2 \left(0, T; \mathcal{V}^{*} \right)\,\, \mbox{and}\,\, \bm{y}(0) = \bm{y}_0 \Big\}$, where we denoted by $\mathcal{V}^{*}$ the dual-space of $\mathcal{V}$. 
The outflow control variable $\bm{u}$ , influencing the fluid flow physics, belongs to the Hilbert space $\mathcal{U} = L^{2}\left(0, T; \mathbf{U} \right)$ with $\mathbf{U} = L^{2}\left( \Gamma_\mathrm{out} \right).$ 

The weak formulation of parametric state Eq.~\eqref{eq:2.3} is given as: find the pair $\left( \bm{v} \left(t; \bm{\mu} \right), p \left(t; \bm{\mu} \right) \right) \in \mathcal{V}$  such that:
\begin{align}\label{eq:2.4}
\begin{cases}
\begin{aligned}
\displaystyle 
 &   m\left(\bm{v}, \bm{z_s}; t; \bm{\mu} \right) + a \left(\bm{v}, \bm{z_s}; t; \bm{\mu} \right) + e \left(\bm{v}, \bm{v}, \bm{z_s}; t; \bm{\mu} \right) + b\left(p, \bm{z_s}; t; \bm{\mu} \right) + c \left(\bm{u}, \bm{z_s}; t; \bm{\mu} \right) =  F \left(\bm{z_s}, t;\bm{\mu} \right),    &&\forall \bm{z_s}  \in \mathbf{V}, \\
 &   b \left(l_s, \bm{v}; t; \bm{\mu} \right) = 0        &&\forall l_s \in \mathbf{P}. 
\end{aligned}
\end{cases}
\end{align}
Eq.~\eqref{eq:2.4} consists of several terms: the inertial term $ \displaystyle m\left(\bm{v}, \bm{z_s}; t; \bm{\mu} \right),$ which represents the time evolution of the velocity field, the nonlinear convective term $\displaystyle e \left(\bm{v}, \bm{v}, \bm{z_s}; t; \bm{\mu} \right),$ the viscous term $a \left(\bm{v}, \bm{z_s}; t; \bm{\mu}\right),$ the pressure term $\displaystyle b \left(p, \bm{z_s}; t; \bm{\mu} \right),$ the control term $\displaystyle c \left(\bm{u}, \bm{z_s}; t; \bm{\mu} \right),$ and the external force term $F\left(\bm{z_s}; t;\bm{\mu} \right).$  These terms are expressed as follows:
\begin{align}\label{eq:2.4.1}
\begin{gathered}
    m\left(\bm{v}, \bm{z_s}; t; \bm{\mu} \right) = \int_{\Omega} \frac{\partial \bm{v}\left(t; \bm{\mu} \right)}{\partial t} \cdot \bm{z_s}\left(t; \bm{\mu} \right) \, d\Omega, \qquad\qquad 
    e \left(\bm{v}, \bm{v}, \bm{z_s}; t; \bm{\mu} \right) = \int_{\Omega} \left(\bm{v} \left(t; \bm{\mu} \right) \cdot \nabla \right) \bm{v} \left(t; \bm{\mu} \right) \cdot \bm{z_s} \left(t; \bm{\mu} \right) \, d\Omega, \\
    a\left(\bm{v}, \bm{z_s}; t; \bm{\mu}\right) = \nu \int_\Omega \nabla \bm{v}\left(t; \bm{\mu} \right) : \nabla \bm{z_s} \left(t; \bm{\mu} \right) \, d\Omega, \qquad\qquad
    b \left(p, \bm{z_s}; t; \bm{\mu} \right) = -\int_{\Omega} p\left(t; \bm{\mu} \right) \left(\nabla \cdot \bm{z_s}\left(t; \bm{\mu} \right)\right) \, d\Omega, \\
    c \left(\bm{u}, \bm{z_s}; t; \bm{\mu} \right) = - \int_{\Gamma_\mathrm{out}} \bm{u}\left(t; \bm{\mu} \right) \cdot \bm{z_s}\left(t; \bm{\mu} \right) \, d\Gamma, \qquad\qquad 
    F\left(\bm{z_s}; t;\bm{\mu} \right) = \int_{\Omega} \bm{f}\left(t; \bm{\mu} \right) \cdot \bm{z_s}\left(t; \bm{\mu} \right) \, d\Omega.
\end{gathered}
\end{align}
For the \textit{existence and uniqueness} of solution of the state equations, the bilinear operator $a \left(\bm{v}, \bm{z_s}; t; \bm{\mu} \right): \mathbf{V} \times \mathbf{V} \rightarrow \mathbb{R}$, have to be coercive and bounded; and the operator $b\left(p, \bm{z_s}; t; \bm{\mu} \right): \mathbf{P} \times \mathbf{V} \rightarrow \mathbb{R}$, have to satisfy the \textit{Ladyzhenskaya-Babuška-Brezzi (LBB) inf-sup condition}~\cite{temam1977navier}, expressed as 
\begin{equation}\label{eq:inf}
   \displaystyle \exists\,\,\,\, \gamma > 0 \,\,  \text{ such that  }  \inf_{p \neq 0 \in \mathbf{P}} \sup_{\bm{z_s} \neq 0\in  \mathbf{V}} \frac{b \left(p, \bm{z_s}; t; \bm{\mu} \right)}{\| p \|_\mathbf{P} \| \bm{z_s} \|_\mathbf{V}} \, \geq \,\, \gamma. 
\end{equation}
%%
%%----------------------------------------------------%%
\subsection{Lagrangian Formulation}
\label{sec2.2:Lagrange}

To solve the minimization problem, we apply the Lagrangian formulation to the OCP$_{(\bm{\mu})}$ governed by the unsteady N-S equations. 
In practice, we introduce the adjoint variables and derive the so-called first-order optimality conditions. 
We adopt a monolithic approach, where the state equations, adjoint equations, and optimality conditions are solved simultaneously as a coupled system, which is based on the Karush-Kuhn-Tucker (KKT) conditions ~\cite{gunzburger2002,protas2008adjoint,gunzburger1999sensitivities}, ensuring consistency between them. The Lagrangian operator $\mathcal{L}$ combines the cost functional $\mathcal{J}$ with constraints from the unsteady N-S equations, ensuring optimal solutions for fluid flow. In this framework, the optimum is a \textit{saddle-point} of the Lagrangian functional, as discussed in \cite{brezzi1974existence} with KKT optimality conditions.

For any given $\bm{\mu} \in \mathcal{D}$, the Lagrangian functional 
 $ \displaystyle \mathcal{L} \left( \bm{X}\left(t; \bm{\mu} \right); t; \bm{\mu} \right): \mathcal{X} \rightarrow \mathbb{R}$, where $\mathcal{X} := \mathcal{V} \times \mathcal{U} \times \mathcal{V}$, is defined as:
\begin{equation}
\begin{split}\label{eq:2.6}
    \mathcal{L} \left( \bm{v}, p, \bm{u}, \bm{w}, q; t; \bm{\mu} \right) \ =\ &  \frac{1}{2} \int_{0}^{T} \int_{\Omega} \| \left( \bm{v}\left(t; \bm{\mu} \right) - \bm{v}_d; t; \bm{\mu} \right) \|^2 \, d\Omega \, dt + \frac{\alpha}{2} \int_{0}^{T} \int_{\Gamma_\mathrm{out}} \| \bm{u}\left(t; \bm{\mu} \right) \|^2 \, d\Gamma \, dt  \\
   & -\int_0^T \int_{\Omega} \left( \bm{w}\left(t; \bm{\mu} \right) \cdot \left(\frac{\partial \bm{v} \left(t; \bm{\mu} \right)}{\partial t} +\left( \bm{v}\left(t; \bm{\mu} \right) \cdot \nabla \right)\, \bm{v}\left(t; \bm{\mu} \right)  - \nu\, \Delta \bm{v}\left(t; \bm{\mu} \right) + \nabla p\left(t; \bm{\mu} \right) - f\left(t; \bm{\mu} \right) \right) \right) \, d\Omega \, dt 
   \\
   & -\int_0^T \int_{\Omega} q\left(t;  \bm{\mu} \right) \nabla \cdot \bm{v}{\left(t; \bm{\mu} \right)} \, d\Omega \, dt. 
\end{split}
\end{equation}
\\
If we denote by $ \displaystyle \bm{X} \left(t; \bm{\mu} \right) \doteq \left(\bm{v}\left(t; \bm{\mu} \right), p\left(t;  \bm{\mu} \right), \bm{u}\left(t; \bm{\mu} \right), \bm{w}\left(t; \bm{\mu} \right), q\left(t; \bm{\mu} \right) \right),$ the unknowns, where $\displaystyle \bm{z}\left(t; \bm{\mu} \right) = \left( \bm{w}\left(t; \bm{\mu} \right), q\left(t; \bm{\mu} \right) \right) \in \mathcal{V}$ is the Lagrangian multiplier, we aim at computing the optimal solution $\bm{X}^{*}(t; \bm{\mu}) := \left( \bm{y}^{*}\left(t; \bm{\mu} \right),  \bm{u}^{*}\left(t; \bm{\mu} \right), \bm{z}^{*}\left(t; \bm{\mu} \right) \right) \in \mathcal{X}$ for the OCP$_{(\bm{\mu})}$ defined in Eqs.~\eqref{eq:2.2}-\eqref{eq:2.3}, that satisfies the \textit{first order necessary conditions} for any given $\bm{\mu} \in \mathcal{D}$:
\begin{alignat}{3}
\mathcal{L}_{\bm{z}}(\bm{X}^{*}\left(t; \bm{\mu} \right); t; \bm{\mu})[\phi\left(t; \bm{\mu} \right)] &= \quad \langle \mathcal{L}_{\bm{z}}(\bm{X}^{*}\left(t; \bm{\mu} \right); t; \bm{\mu}), \phi\left(t; \bm{\mu} \right) \rangle  \quad &= 0 & \quad & \forall \phi \in \mathcal{V},\label{eq:2.8}
\\
\mathcal{L}_{y}(\bm{X}^{*}\left(t; \bm{\mu} \right); t; \bm{\mu})[\psi\left(t; \bm{\mu} \right)] &= \quad
\langle \mathcal{L}_{\bm{y}}(\bm{X}^{*}\left(t; \bm{\mu} \right); t; \bm{\mu}), \psi\left(t; \bm{\mu} \right) \rangle \quad &= 0 & \quad & \forall \psi \in \mathcal{V},  \label{eq:2.9}
\\
\mathcal{L}_{u}(\bm{X}^{*}\left(t; \bm{\mu} \right); t; \bm{\mu})[\xi\left(t; \bm{\mu} \right)] &= \quad
\langle \mathcal{L}_{\bm{u}}(\bm{X}^{*}\left(t; \bm{\mu} \right); t; \bm{\mu}), \xi\left(t; \bm{\mu} \right) \rangle \quad &= 0 & \qquad & \forall \xi \in \mathcal{U}, \label{eq:2.10}
\end{alignat}
%%
% with any $ \displaystyle \bar{\bm{X}}\left(t; \bm{\mu} \right)  = \left( \bar{\bm{y}}\left(t; \bm{\mu} \right), \bar{\bm{u}}\left(t; \bm{\mu} \right), \bar{\bm{z}}\left(t; \bm{\mu} \right) \right) \in \mathcal{X},$ 
therein, $\phi$, $\psi$ and $\xi$ are the test functions in the corresponding function spaces $\mathcal{V}$ and $\mathcal{U}$, respectively. Here,  Fr\`echet differentiation is used for the Lagrangian functional~\eqref{eq:2.6}, allowing us to obtain a system of PDEs, as detailed in the report \cite{manzoni2019saddle}. More precisely, Eq.~\eqref{eq:2.8} leads to the state equations presented in Eq.~\eqref{eq:2.3}, while Eq.~\eqref{eq:2.9} and Eq.~\eqref{eq:2.10} lead, respectively, to the \textit{adjoint equations}:
\begin{align}\label{eq:2.12}
%\text{(Adjoint Equations)} \quad 
\begin{cases}
\begin{aligned}
\displaystyle 
&-\frac{\partial \bm{w} \left(t; \bm{\mu} \right)}{\partial t} - \left(\bm{v} \left(t; \bm{\mu} \right) \cdot \nabla \right)\, \bm{w} \left(t; \bm{\mu} \right)
 +  \left( \nabla \bm{v}\left(t; \bm{\mu} \right) \right)^{T} \bm{w} \left(t; \bm{\mu} \right) \\ & \qquad\qquad\qquad\qquad\qquad + \nu\, \Delta \bm{w}\left(t; \bm{\mu} \right)   - \nabla q\left(t; \bm{\mu} \right) = \left( \bm{v}\left(t; \bm{\mu} \right) - \bm{v}_d \right),  && \quad    \mbox{in }  \Omega\times \left(0, T\right), \\
&\nabla \cdot \bm{w}{\left(t; \bm{\mu} \right)} = 0,  && \quad  \mbox{in } \Omega\times \left(0, T\right), \\
&\bm{w}{\left(t; \bm{\mu} \right)} = 0, && \quad \mbox{on } \Gamma_{in} \times \left(0, T\right), \\
&\bm{w}\left(T; \bm{\mu} \right) = 0, && \quad  \mbox{at } \Omega \times \{T\}, \\
&\nabla \bm{w}\left(t; \bm{\mu} \right) \cdot \bm{n} = 0,  && \quad   \mbox{on } \Gamma_{out }\times \left(0, T\right).
\end{aligned}
\end{cases}
\end{align}
and the following \textit{optimality condition}
\begin{equation}\label{eq:2.13}
   \alpha \, \bm{u}\left(t; \bm{\mu}\right) + \bm{w}\left(t; \bm{\mu}\right) =0, \quad \mbox{on } \Gamma_{out }\times \left(0, T\right).
\end{equation} 
%%
%%%%%%%%%%%%%%%%%%%%%%%%%%%%%%%%%%%%%%%%%%%%%%%%%%%%%%%%%%%%%%%%%%%%%%%%%%%%%%%%%%%%%%%%%%

%%-------------------------------------------------------%%
\section{Methods} %{Numerical Approximation}
\label{sec3:Numerical}
In this section, we delve into the numerical solution of OCP$_{(\bm{\mu})}$s comprising of the cost functional, nonlinear state equations, adjoint equations, and optimality conditions, by adopting the \textit{optimize-then-discretize} methodology, as detailed in \cite{manzoni2012reduced}. This approach consists in optimizing at the continuous level, and thus deriving the first order optimality condition detailed in the previous section, and then discretizing the obtained optimality system in a \textit{one-shot} manner, meaning that we aim at solving the state, control, and adjoint equations as a unique block system. Given the computational expense of such parametric problems, we introduce a projection-based ROM tailored for such optimization problems, that involves compressing the temporal behaviour, extract the dominant information, and subsequently compressing the parametric dependence of high-fidelity solutions towards efficient simulations of the dynamics \cite{ballarin2016fast,kadeethum2022non}.

%%-------------------------------------------------------%%
\subsection{Galerkin Finite Element Formulation}
\label{sec3.1:FEM}

The complete weak formulation of the coupled optimality system for the optimization problem~\eqref{eq:2.2}--\eqref{eq:2.3} reads as: find $\left( \bm{v}, p, \bm{u}, \bm{w}, q \right) \in \mathcal{X}$  such that: 
\begin{align}\label{eq:2.14}
\begin{cases}
\begin{aligned}
\displaystyle
 & m\left(\bm{v}, \bm{z_s}; t; \bm{\mu} \right) 
 + a \left(\bm{v}, \bm{z_s}; t; \bm{\mu} \right) 
 + e \left(\bm{v}, \bm{v}, \bm{z_s}; t; \bm{\mu} \right) \\
 & \qquad\qquad\qquad + b \left(p, \bm{z_s}; t; \bm{\mu} \right) 
 + c \left(\bm{u}, \bm{z_s}; t; \bm{\mu} \right) 
 = F\left(\bm{z_s}; t;\bm{\mu} \right), 
 && \forall \bm{z_s} \in \mathbf{V} \\
 & b \left(l_{\bm{s}}, \bm{v}; t; \bm{\mu} \right) =  0, 
 && \forall l_{\bm{s}} \in \mathbf{P} \\
 & - m\left(\bm{w}, \bm{z_a}; t; \bm{\mu} \right) 
 + a \left(\bm{w}, \bm{z_a}; t; \bm{\mu} \right) + e \left(\bm{v}, \bm{w}, \bm{z_a}; t; \bm{\mu} \right) \\
 & \qquad\qquad\qquad + e \left(\bm{w}, \bm{v}, \bm{z_a}; t; \bm{\mu} \right) + b \left(q, \bm{z_a}; t; \bm{\mu} \right) 
 = G\left(\bm{z_a}; t;\bm{\mu} \right), 
 && \forall \bm{z_a} \in \mathbf{V} \\
 & b \left(l_{\bm{a}} , \bm{w}; t; \bm{\mu} \right) = 0, 
 && \forall l_{\bm{a}} \in \mathbf{P}, \\
 & \left( \alpha \, \bm{u} \left(t; \bm{\mu} \right) 
 + \bm{w} \left(t; \bm{\mu} \right) \right) 
 \cdot \bm{\eta} = 0,  
 && \forall  \bm{\eta} \in \mathcal{U},
\end{aligned}
\end{cases}
\end{align}
whereas $F$ and $G$ denote the source term and the target flow profile, respectively.
For computing the solution of Eq.~\eqref{eq:2.14} utilizing the Galerkin FE formulation, suitable finite-dimensional approximation spaces are required. Let us define the discrete FE spaces for state $\mathcal{V}_h \subset \mathcal{V}$, and control $\mathcal{U}_h \subset \mathcal{U}$. The finite-dimensional optimal flow control problem is then defined on the space $\mathcal{X}_h \equiv \mathcal{V}_h \times \mathcal{U}_h  \times \mathcal{V}_h  \subset \mathcal{X}$.
Let us rewrite the Eq.~\eqref{eq:2.14} in the discrete setting as: 
\begin{align}\label{eq:2.15}
\begin{cases}
\begin{aligned}
 & m_h\left(\bm{v}_h, \bm{z}_{s_h}; t; \bm{\mu} \right) + a_h \left(\bm{v}_h, \bm{z}_{s_h}; t; \bm{\mu} \right) + e_h \left(\bm{v}_h, \bm{v}_h, \bm{z}_{s_h}; t; \bm{\mu} \right)\\  & \qquad\qquad\qquad + b_h \left(p_h, \bm{z}_{s_h}; t; \bm{\mu} \right) + c_h \left(\bm{u}_h, \bm{z}_{s_h}; t; \bm{\mu} \right) = F_h\left(\bm{z}_{s_h}; t; \bm{\mu} \right), && \forall \bm{z}_{s_h} \in \mathbf{V}_h, \\
 & b_h \left(l_{\bm{s_h}}, \bm{v}_h; t; \bm{\mu} \right) = 0, && \forall l_{\bm{s_h}} \in \mathbf{P}_h, \\
 & -m_h \left(\bm{w}_h, \bm{z}_{a_h}; t; \bm{\mu} \right) + a_h \left(\bm{w}_h, \bm{z}_{a_h}; t; \bm{\mu} \right) + e_h \left(\bm{v}_h, \bm{w}_h, \bm{z}_{a_h}; t; \bm{\mu} \right)\\  & \qquad\qquad\qquad + e_h \left(\bm{w}_h, \bm{v}_h, \bm{z_{a_h}}; t; \bm{\mu} \right)  + b_h \left(q_h, \bm{z}_{a_h}; t; \bm{\mu} \right) = G_h \left(\bm{z}_{a_h}; t; \bm{\mu} \right), && \forall \bm{z}_{a_h} \in \mathbf{V}_h, \\
 & b_h \left(l_{\bm{a_h}}, \bm{w}_h; t; \bm{\mu} \right) = 0, && \forall l_{\bm{a_h}} \in \mathbf{P}_h, \\
 & \left( \alpha \, \bm{u}_h \left(t; \bm{\mu} \right) + \bm{w}_h \left(t; \bm{\mu} \right) \right) \cdot \bm{\eta}_h = 0, && \forall \bm{\eta}_h \in \mathcal{U}_h.
\end{aligned}
\end{cases}
\end{align}
The subscript `$h$' denotes the spatially discretized form of the terms presented in Eq.~\eqref{eq:2.4.1}.
We can write the above nonlinear system in a more compact form as follows: 
\begin{align}\label{eq:2.16}
\begin{cases}
\begin{aligned}
 &  M\left(t; \bm{\mu} \right)  + A \left(t; \bm{\mu} \right)\,\bm{v}_h + E\left(\bm{w}_h \right)\left( t; \bm{\mu} \right)\,\bm{v}_h + B \left(t; \bm{\mu} \right) p_h + C \left(t; \bm{\mu} \right) \,\bm{u}_h = F\left(t; \bm{\mu} \right), \\
 &   B^\top \left( t; \bm{\mu} \right)\bm{v}_h =  0, \\
 &  - M^\top \left(t; \bm{\mu} \right)  +  A_{ad}\left(t; \bm{\mu} \right)\bm{w}_h  + E_{ad}\left( \bm{v}_h \right)\left( t; \bm{\mu} \right) \bm{w}_h  + E_{ad}^\top  \left(\bm{w}_h \right) \left(t; \bm{\mu} \right) \bm{v}_h  + B_{ad} \left(t; \bm{\mu} \right) q_h =   G \left(t; \bm{\mu} \right),  \\ 
 &   B_{ad}^\top \left(t; \bm{\mu} \right)\bm{w}_h =  0, \\
 &  \alpha \, \bm{u}_h \left(t; \bm{\mu} \right) + \bm{w}_h \left(t; \bm{\mu} \right) =0,
\end{aligned}
\end{cases}
\end{align}
where the quantities $M$, $A$, $E$, $B$, and $C$ correspond to the inertial, diffusion, advection, pressure, and control terms, respectively. The subscript `$ad$' corresponds to the adjoint quantities. 
Within the Galerkin approximation, we express the velocity, pressure and control variables as expansions over the respective basis functions\footnote{For simplicity, we have only expressed the state and control variables here. The adjoint variables can be expressed in a similar manner.}: 
$$  \bm{v}_h = \sum_{i=1}^{\mathcal{N}_v}\bm{v}_{h}^{\left(i \right)} \,\phi_i  \in \mathbf{V}_h,    \hspace{10mm}
    {p}_h = \sum_{j=1}^{\mathcal{N}_p}p_{h}^{\left(j \right)} \, \xi_j, \in \mathbf{P}_h,  \hspace{10mm}
\text{ and }  \hspace{5mm}   \bm{u}_h = \sum_{k=1}^{\mathcal{N}_u}\bm{u}_{h}^{\left(k \right)} \, \zeta_k, \in \mathcal{U}_h. $$ 
Here, $\mathcal{N}_v$, $\mathcal{N}_p$, and $\mathcal{N}_u$ represent the dimensions of FE subspaces for the velocity, pressure and control variables, respectively, and the resulting dimension of the high-fidelity system is $\mathcal{N}_h = 2\left( \mathcal{N}_v + \mathcal{N}_p \right) + \mathcal{N}_u$ (the state and adjoint variables belong to the same space, and are spanned by the same basis functions). The vector coefficients for velocity, pressure, and control variables, respectively, are as $$\displaystyle \bm{v}_h = \left( v_h^{\left(1 \right)}, v_h^{\left(2 \right)},  \cdots, v_h^{\left(N_v \right)} \right)^\top \in \mathbb{R}^{\mathcal{N}_v}, \,\,\,\,\,\,\,\,\,\,\,\, \bm{p}_h = \left( p_h^{\left(1 \right)}, p_h^{\left(2 \right)},  \cdots, p_h^{\left(N_p \right)} \right)^\top \in \mathbb{R}^{\mathcal{N}_p},$$ and $$\displaystyle \bm{u}_h = \left(u_h^{\left(1 \right)}, u_h^{\left(2 \right)},  \cdots, u_h^{ \left(N_u \right)} \right)^\top \in \mathbb{R}^{\mathcal{N}_u}.$$ In particular, for the numerical discretization of the system, we exploited the Newton method to handle the nonlinear convective term, and the Implicit Euler scheme for time-discretization, partitioning the interval $[0, T]$ into $N_t$ equally spaced time intervals with time-step $\Delta t$. 
\\
Finally, the optimization problem can be written in algebraic form as follows:
\begin{equation*}%\label{eq:2.17}
\resizebox{\textwidth}{!}{$ % Start resizebox
\begin{aligned}
\begin{bmatrix}
\displaystyle
\left( \frac{1}{\Delta t}M + M_d  + E_{ad}^\top \left( \bm{w}_h \right) \right) \left(t^{\left(n+1\right)}; \bm{\mu} \right)  & 0 & 0 & \left( A_{ad} + E_{ad} \left(\bm{v}_h \right) \right) \left(t^{\left(n+1 \right)}; \bm{\mu} \right)  & B_{ad} \left(t^{ \left(n+1 \right)}; \bm{\mu} \right)\\
0 & 0 & 0  & B_{ad}^\top \left( t^{ \left(n+1 \right)}; \bm{\mu} \right) & 0 \\
0 & 0 & \mathcal{R} \left( t^{ \left(n+1 \right)}; \bm{\mu} \right)  &  0 & 0 \\
\left( A + E \left(\bm{w}_h \right)  \right) \left(t^{\left(n+1 \right)}; \bm{\mu} \right) & B \left(t^{\left(n+1 \right)}; \bm{\mu} \right) &  C \left( t^{\left(n+1 \right)}; \bm{\mu} \right) & \frac{-1}{\Delta t}M^\top \left(t^{ \left(n+1 \right)}; \bm{\mu} \right) & 0 \\
B^\top \left(t^{ \left(n+1 \right)}; \bm{\mu} \right) & 0 & 0  & 0 & 0 \\
\end{bmatrix}
\begin{bmatrix}
\bm{v} \left(t^{\left(n+1 \right)}; \bm{\mu} \right)\\
p \left(t^{ \left(n+1 \right)}; \bm{\mu} \right) \\
\bm{u} \left(t^{ \left(n+1 \right)}; \bm{\mu} \right)\\
\bm{w} \left(t^{ \left(n+1 \right)}; \bm{\mu} \right)\\
q \left(t^{ \left(n+1 \right)}; \bm{\mu} \right)\\
\end{bmatrix}
\\  =
\begin{bmatrix}
\bm{G} \left(t^{\left(n+1 \right)};\bm{\mu} \right) \\
0 \\
\bm{H} \left(t^{ \left(n+1 \right)};\bm{\mu} \right) \\
\bm{F} \left(t^{\left(n+1 \right)};\bm{\mu} \right) \\
0
\end{bmatrix}
+
\begin{bmatrix}
\frac{1}{\Delta t}M \left(t^n; \bm{\mu} \right)  & 0 & 0 & 0 & 0\\
0 & 0 & 0  & 0 & 0 \\
0 & 0 & 0  & 0 & 0 \\
0 & 0 & 0  & \frac{-1}{\Delta t}M^\top \left(t^n; \bm{\mu} \right) & 0 \\
0 & 0 & 0  & 0 & 0 \\
\end{bmatrix}
\begin{bmatrix}
\bm{v}\left(t^n; \bm{\mu} \right)  \\
p\left(t^n; \bm{\mu} \right) \\
\bm{u}\left(t^n; \bm{\mu} \right)  \\
\bm{w}\left(t^n; \bm{\mu} \right) \\
q \left(t^n; \bm{\mu} \right) 
\end{bmatrix}
\end{aligned}
$} % End resizebox
\end{equation*}
where, and we are defining the matrices as, 
\begin{align}\label{eq:2.18}
\resizebox{\textwidth}{!}{$ % 
\begin{gathered}
\displaystyle
  \left( M \left(t;\bm{\mu} \right) \right)_{ij} = m \left(\phi_j, \phi_i; t; \bm{\mu} \right),  
 \qquad \left( M_d \left(t;\bm{\mu} \right) \right)_{ij} = m_d \left( \phi_j, \phi_i; t; \bm{\mu} \right),
 \qquad \left( A \left(t;\bm{\mu}\right) \right)_{ij} = a \left( \phi_j, \phi_i; t; \bm{\mu} \right),\\ 
  \left( E \left(\bm{v}(t; \bm{\mu} \right) \right)_{ij} = \sum_{m=1}^{\mathcal{N}_v} \bm{v}_h^{m} \left(t; \bm{\mu} \right) \, e \left( \phi_m, \phi_j, \phi_i; t; \bm{\mu} \right),
 \qquad \left( C \left(t;\bm{\mu}\right) \right)_{ij} = c \left( \phi_j, \zeta_i; t; \bm{\mu} \right),   
 \qquad \left( B \left(t;\bm{\mu} \right) \right)_{ij} = b \left( \phi_j, \xi_i; t; \bm{\mu} \right), \\
  \left( F \left(t; \bm{\mu} \right) \right)_{i} = \bm{f} \left( \phi_i; t; \bm{\mu} \right),
 \qquad \left( G \left(t; \bm{\mu} \right) \right)_{i} = \bm{g} \left( \phi_i; t; \bm{\mu} \right), 
 \qquad \left( H \left(t; \bm{\mu} \right) \right)_{i} = h \left( \zeta_i; t; \bm{\mu} \right).
\end{gathered}
$}
\end{align}
Similarly, we can define the matrices for the adjoint operators. 
%%-------------------------------------------------------%%
\subsection{Projection-based Reduced Framework for Parametrized OCPs}
\label{sec3.2:ROM}
In this section, we present a ROM for parametrized OCPs for unsteady N--S equations based on a POD-Galerkin technique. %technique and a Galerkin projection. 
We are utilizing the affine decomposition assumption, and expressing the bi(linear) operators, which are mentioned in Eq.~\eqref{eq:2.18}, as follows: 
\begin{align}\label{eq:2.19}
\resizebox{\textwidth}{!}{$ % Start resizebox
\begin{gathered}
\displaystyle
  A\left(t; \bm{\mu} \right) := \sum_{q=1}^{Q_a} \Theta_{q}^{Q_a}\left(t; \bm{\mu} \right)\, A^q\left(\cdot, \cdot \right),
 \qquad M\left(t; \bm{\mu} \right) := \sum_{q=1}^{Q_m} \Theta_{q}^{Q_m}\left(t; \bm{\mu} \right)\, m^q\left(\cdot, \cdot \right), 
 \qquad E\left(\bm{v}\right)\left(t; \bm{\mu}\right) := \sum_{q=1}^{Q_e} \Theta_{q}^{Q_e}\left(t; \bm{\mu} \right)\, e\left(\bm{v}\right)^q\left(\cdot, \cdot \right), \\
  M_d \left(t; \bm{\mu} \right) :=  \sum_{q=1}^{Q_{m_d}} \Theta_{q}^{Q_{m_d}}(t; \bm{\mu})\, m_d^d\left(\cdot, \cdot \right), 
 \qquad C\left(t; \bm{\mu} \right) := \sum_{q=1}^{Q_c} \Theta_{q}^{Q_c} \left(t; \bm{\mu} \right)\, C^q\left(\cdot, \cdot \right), 
 \qquad B\left(t; \bm{\mu} \right) := \sum_{q=1}^{Q_b} \Theta_{q}^{Q_b} \left(t; \bm{\mu} \right)\, B^q\left(\cdot, \cdot \right), \\
  F\left(t; \bm{\mu} \right) := \sum_{q=1}^{Q_f} \Theta_{q}^{Q_f} \left(t; \bm{\mu} \right)\, F^q\left(\cdot, \cdot \right), 
 \qquad G\left(t; \bm{\mu} \right) := \sum_{q=1}^{Q_g} \Theta_{q}^{Q_g} \left(t; \bm{\mu} \right)\, G^q\left(\cdot, \cdot \right), 
 \qquad H\left(t; \bm{\mu} \right) := \sum_{q=1}^{Q_h} \Theta_{q}^{Q_h} \left(t; \bm{\mu} \right)\, H^q \left(\cdot, \cdot \right), 
\end{gathered}
$}
\end{align}
similarly for the other terms. Affine decomposition assumption expresses the operators as sums of separable functions $\Theta_{q}^{Q_*} \left(t; \bm{\mu} \right)$  of parameter $\bm{\mu}$ and time $t$, multiplied by computed spatial components. This decomposition is crucial for implementing the projection-based ROM and significantly reducing the computational complexity of the  OCP$_{\left( \bm{\mu} \right)}$s.

Let us consider a training set $\mathbb{E}_\mathrm{train} = \{ \bm{\mu}^1, \bm{\mu}^2, \cdots, \bm{\mu}^{N_\mathrm{train}} \} \in \mathcal{D}$ of dimension $N_\mathrm{train}$, with randomly chosen values in a specified range. 
To explore the parametric dependence of the OCPs, the offline phase consists in computing the high-fidelity solutions for each parameter value in the training set $\mathbb{E}_\mathrm{train}$, while storing the temporal evolution of the high-fidelity solutions in snapshots in proper matrices along the temporal trajectory. 
For any $\bm{\mu}^{i} \in \mathbb{E}_\mathrm{train}$ where $ i = 1, 2, \dots, N_\mathrm{train}$ we have: 
\begin{align}\label{eq:3.1}
\displaystyle
 S_{\bm{v}}^{i} =& \left[ \bm{v} \left(t^0; \bm{\mu}^{i} \right) \,\,\,|\,\,\, \bm{v} \left(t^1; \bm{\mu}^{i} \right) \,\,\,| \,\,\, \cdots \,\,\,|\,\,\, \bm{v} \left(t^{N_t -1}; \bm{\mu}^{i} \right) \right] \in \mathbb{R}^{{\mathcal{N}^{v}_{h}} \times N_t},  \nonumber \\
 S_{\bm{p}}^{i} =& \left[ \bm{p} \left(t^0; \bm{\mu}^{i} \right) \,\,\,|\,\,\, \bm{p} \left(t^1; \bm{\mu}^{i} \right) \,\,\,| \,\,\, \cdots \,\,\,|\,\,\, \bm{p} \left(t^{N_t -1}; \bm{\mu}^{i} \right) \right] \in \mathbb{R}^{{\mathcal{N}^{p}_{h}} \times N_t},  \\
 S_{\bm{u}}^{i} =& \left[ \bm{u} \left(t^0; \bm{\mu}^{i} \right) \,\,\,|\,\,\, \bm{u} \left(t^1; \bm{\mu}^{i} \right) \,\,\,| \,\,\, \cdots \,\,\,|\,\,\, \bm{u} \left(t^{N_t -1}; \bm{\mu}^{i} \right) \right] \in \mathbb{R}^{{\mathcal{N}^{u}_{h}} \times N_t},  \nonumber 
\end{align}
where $\mathcal{N}^{h}$ and $N_t$ represent the number of spatial degrees of freedom and the number of time snapshots, respectively.
Similarly for the adjoint variables $\left( \bm{w}, \bm{q} \right) \in \mathcal{V}_N$.  In incompressible flow simulations, the \textit{inf-sup (LBB) condition} is critical for ensuring the stability of the pressure-velocity coupling. This condition is crucial to prevent spurious pressure oscillations and numerical instabilities.  For the parametrized bilinear operator $b \left(p_h, \bm{z_{s_h}}; t; \bm{\mu} \right),$  it is expressed as, $\exists \text{  a constant  }  \tilde{\gamma_h} > 0   \text{ such that   }$
\begin{equation}\label{eq:3.2}
   \displaystyle  \gamma_h \left(\bm{\mu} \right) =   \inf_{\bm{p}_h \neq 0\in \mathbf{P}_h} \sup_{\bm{z_{s_h}} \neq 0 \in  \mathbf{V}_h} \frac{ b_h \left(p_h, \bm{z_{s_h}}; t; \bm{\mu} \right)}{\| \bm{p}_h \|_\mathbf{P} \| \bm{z_{s_h}} \|_\mathbf{V}} \geq \tilde{\gamma_h} > 0 \,\,\,\,\,\, \forall \bm{\mu} \in \mathcal{D}, 
\end{equation}
To address this issue in reduced-order models, we introduce the supremizer~\cite{ballarin2015supremizer,quarteroni2014reduced}. \\
We define the \textit{supremizer} operator $ \mathcal{T}_h^{\bm{\mu}} : \mathbf{P}_h  \rightarrow \mathbf{V}_h$  as follows: 
\begin{equation}\label{eq:3.3}
   \left( \mathcal{T}_h^{\bm{\mu}} \bm{p}_h, \bm{z_{s_h}} \right)_{\mathbf{V}_h} = b_h \left( p_h, \bm{z_{s_h}}; t; \bm{\mu} \right),  \hspace{10mm}   \forall \bm{z_{s_h}} \in \mathbf{V}_h,
\end{equation}
enriches the velocity space to maintain the inf-sup stability. Thus, we introduced the state supremizer $\mathbf{s}^{\bm{\mu}}\left(\bm{p} \left(t; \bm{\mu} \right) \right)$ and the adjoint supremizer $\mathbf{r}^{\bm{\mu}}\left(\bm{q}\left(t; \bm{\mu} \right) \right),$ respectively. These supremizer ensures that the velocities space have sufficient degrees of freedom to couple correctly with the pressure space, maintaining stability in the reduced model. To incorporate the supremizer into the ROM, we construct snapshot matrices for these supremizers along the temporal trajectory for each parameter $\bm{\mu}^{i} \in \mathbb{E}_\mathrm{train}, \, \, i = 1, 2, \dots, N_\mathrm{train},$ are defined as follows:
\begin{align}\label{eq:3.4}
\displaystyle
S_{\bm{s}}^{i} =& \left[\mathbf{s}^{\bm{\mu}} \left(\bm{p} \left(t^0; \bm{\mu}^{i} \right) \right)\,\,\,|\,\,\, \mathbf{s}^{\bm{\mu}}\left(\bm{p}\left(t^1; \bm{\mu}^{i} \right)\right) \,\,\,| \,\,\, \cdots \,\,\,|\,\,\, \mathbf{s}^{\bm{\mu}} \left(\bm{p}\left(t^{N_t -1}; \bm{\mu}^{i} \right) \right) \right] \in \mathbb{R}^{{\mathcal{N}^{v}_{h}} \times N_t},  \nonumber \\
S_{\bm{r}}^{i} =& \left[ \mathbf{r}^{\bm{\mu}} \left(\bm{q} \left(t^0; \bm{\mu}^{i} \right) \right)\,\,\,|\,\,\, \mathbf{r}^{\bm{\mu}} \left(\bm{q}\left(t^1; \bm{\mu}^{i} \right) \right) \,\,\,| \,\,\, \cdots \,\,\,|\,\,\, \mathbf{r}^{\bm{\mu}} \left(\bm{q} \left(t^{N_t -1}; \bm{\mu}^{i} \right) \right)\right] \in \mathbb{R}^{{\mathcal{N}^{v}_{h}} \times N_t}.  
\end{align}

Addressing the unsteadiness of high-dimensional systems in OCPs demands substantial computational resources and memory allocation. The standard POD~\cite{rozza2012reduction,ballarin2017numerical,rozza2022advanced} is computationally expensive because it must process large snapshot matrices, representing the state, adjoint and control variables of the system. This challenge is worsened when using finer mesh and smaller time steps, as the number and size of snapshots increase, leading to unbearable computational and memory requirements. Therefore, this study employs a \textit{nested-POD} approach to efficiently manage the computational resources in solving OCP$_{(\bm{\mu})}$s, encompassing the following two steps.

Let us consider the reduced order state space $\mathcal{V}_{rd} \subset \mathcal{V}_h$, and the control space $\mathcal{U}_{rd}\subset \mathcal{U}_h$. Assuming these reduced spaces have been constructed, we can represent the velocity, pressure, and control variables using their respective reduced basis functions as follows: 
$$  \bm{v}_{rd} = \sum_{i=1}^{N_v^{rd}}\, \bm{v}_{rd}^{\left(i \right)} \,\phi_i \in \mathbf{V}_{rd},    \hspace{10mm}
    {p}_{rd} = \sum_{j=1}^{N_p^{rd}}\,p_{rd}^{\left(j \right)} \, \xi_j, \in \mathbf{P}_{rd},  \hspace{10mm}
\text{ and }  \hspace{10mm}   \bm{u}_{rd} = \sum_{k=1}^{N_u^{rd}}\,\bm{u}_{rd}^{\left(k \right)} \, \zeta_k, \in \mathcal{U}_{rd},$$
and their respective vector coefficients are defined as;
$$\displaystyle \underline{\bm{v}}_{rd} = \left(v_{rd}^{\left(1 \right)}, v_{rd}^{\left(2 \right)},  \cdots, v_{rd}^{\left(N_v^{rd} \right)} \right)^\top \in \mathbb{R}^{N_v^{rd}},\,\,\,\,\,\,\,\, \underline{\bm{p}}_{rd} = \left(p_{rd}^{\left(1 \right)}, p_{rd}^{\left(2 \right)},  \cdots, p_{rd}^{\left(N_p^{rd} \right)}\right)^\top \in \mathbb{R}^{N_p^{rd}},$$
and 
$$ \underline{\bm{u}}_{rd} = \left( u_{rd}^{\left(1 \right)}, u_{rd}^{\left(2 \right)},  \cdots, u_{rd}^{\left(N_u^{rd} \right)} \right)^\top \in \mathbb{R}^{N_u^{rd}},$$
where $N_v^{rd}$, $N_p^{rd}$, and $N_u^{rd}$ denote the dimensions of the reduced order spaces for velocity, pressure, and control, respectively. Similarly, the adjoint variables and supremizers can be expressed using their corresponding reduced basis functions.

%%==============================================================================%%
\subsubsection{Temporal compression}
\label{sec3.2.1:temporal}

In this step, our goal is to reduce the number of time steps by extracting the most significant temporal modes from the snapshots. 
For each fixed parameter $\bm{\mu}^{i} \in \mathbb{E}_\mathrm{train}$, where $i = 1, 2, \dots, N_\mathrm{train}$, we compute the Galerkin FE solutions at each time step, storing them in snapshot matrices as described in Eqs.~\eqref{eq:3.1} and \eqref{eq:3.4}. We then apply POD to extract the dominant temporal modes, efficiently capturing the evolution of system for each parameter.

Let us consider the snapshot matrix for state velocity as; 
\begin{equation}\label{eq:3.5}
 S_{\bm{v}}^{i} = \left[ \bm{v}(t^0; \bm{\mu}^{i}) \,\,\,|\,\,\, \bm{v}(t^1; \bm{\mu}^{i}) \,\,\,| \,\,\, \cdots \,\,\,|\,\,\, \bm{v}(t^{N_{t} -1}; \bm{\mu}^{i})\right] \in \mathbb{R}^{{\mathcal{N}^{v}_{h}} \times N_{t}}. 
\end{equation}
The POD basis for the state velocity is obtained by performing the singular value decomposition (SVD) of the above snapshot matrix, given by:
\begin{equation}\label{eq:3.6}
\displaystyle
    S_{\bm{v}}^{i} = L_{\bm{v}}^{i} \, \Sigma_{\bm{v}}^{i}\,{R_{\bm{v}}^{i}}^\top,
\end{equation}
where,
\begin{itemize}
    \item $L_{\bm{v}}^{i} \in \mathbb{R}^{\mathcal{N}_{v}^{h} \times N_t}:$ The matrix of left singular vectors representing the spatial modes,
    \item $\Sigma_{\bm{v}}^{i} \in \mathbb{R}^{N_t \times N_t}:$ Diagonal matrix of singular values $\left(\lambda_{\bm{v}}^{i} \right)$ in descending order, representing the energy of each mode,
    \item $R_{\bm{v}}^{i} \in \mathbb{R}^{N_t \times N_t}:$ The matrix of right singular vectors represents the temporal modes.
\end{itemize}
% $L_{\bm{v}}^{i} \in \mathbb{R}^{\mathcal{N}_{v}^{h} \times N_t}$ contains the first $N_t$ left singular vectors,
% $R_{\bm{v}}^{i} \in \mathbb{R}^{N_t \times N_t}$ is the orthogonal matrix of right singular vectors, and
% $\Sigma_{\bm{v}}^{i} \in \mathbb{R}^{N_t \times N_t}$ is the diagonal matrix of singular 
% values.

% = \left( \sigma_{\bm{v}}^{i} \right)_j^2$.
% % with 
% % The singular values also called the \textit{weighted singular values} for \textit{nested-POD}, are given by  
% \begin{equation}\label{eq:3.7}
%     \displaystyle \left( \Sigma_{\bm{v}}^{i} \right)_{kk} = \left( \sigma_{\bm{v}}^{i} \sigma \right)_k,  %\hspace{5mm} \text{ such that } \hspace{2mm} \left( \sigma_{\bm{v}}^{i} \right)_1 \geq  \left( \sigma_{\bm{v}}^{i} \right)_2 \geq \cdots \geq \left( \sigma_{\bm{v}}^{i} \right)_{N_t} \geq 0.
% \end{equation}
% \end{itemize}  

Typically, we require $N_{t}^\mathrm{POD} \ll N_t$,  and the POD basis is given by the first $N_{t}^\mathrm{POD}$ columns of $L_{\bm{v}}^{i},$ therefore, we can define the basis matrix for velocity  as, for each fixed  $\bm{\mu}^{i} \in \mathbb{E}_\mathrm{train}$, where $i = 1, 2, \dots, N_\mathrm{train}$:
\begin{equation}\label{eq:3.8}
\displaystyle
Z_{\bm{v}}^{i} = \left[ \varphi_1^{i} |\,\,\,  \varphi_2^{i} \,\,\,|\,\,\, \cdots \,\,\,|\,\,\, \varphi_{N_{t}^\mathrm{POD}}^{i} \right] \in \mathbb{R}^{\mathcal{N}^{v}_{h} \times N_{t}^\mathrm{POD}}.
\end{equation}
%%
% and the corresponding the singular values $\Sigma_{\bm{v}}^{i} \in \mathbb{R}^{ N_{t}^\mathrm{POD} \times  N_{t}^\mathrm{POD}}$ also called as \textit{weighted singular values} for \textit{nested-POD}. 
% We are relying on the solution of the so-called \textit{method of snapshots}~\cite{sirovich1987turbulence}, we define the correlation matrix as, for each fixed $\bm{\mu}^{i} \in \mathbb{E}_\mathrm{train}$,  where $i = 1, 2, \dots, N_\mathrm{train}$:
% %%
% \begin{equation}\label{eq:3.9}
% \displaystyle
%     \mathcal{C}_{\bm{v}}^{i} =  {S_{\bm{v}}^{i}}^\top S_{\bm{v}}^{i} \in \mathbb{R}^{N_t  \times N_t},
% \end{equation}
% %%
% and defines the POD bases for velocity by their first $N_{t}^\mathrm{POD}$ singular vectors of the correlation matrix as:
% %%
% \begin{equation}\label{eq:3.10}
% \displaystyle
%     \mathcal{C}_{\bm{v}}^{i}  \left(\varphi_{\bm{v}}^{i} \right)_j = \left( \lambda_{\bm{v}}^{i} \right)_j \, \left(\varphi_{\bm{v}}^{i} \right)_j,  \hspace{10mm} \text{for} \hspace{2mm} j = 1, 2, \cdots, N_{t}^\mathrm{POD}
% \end{equation}
% %%
During the temporal compression, the reduced space dimension $N_{t}^\mathrm{POD}$ is selected as the smallest integer for which the  ``energy" of retained modes, 
\begin{equation}\label{eq:3.11}
\displaystyle
    E_{\bm{v}}^{i} \left(\varphi_1, \varphi_2, \cdots, \varphi_{N_{t}^\mathrm{POD}} \right) = \frac{\Sigma_{k=1}^{N_{t}^\mathrm{POD}} \left(\lambda_{\bm{v}}^{i} \right)_k}{\Sigma_{k=1}^{N_t} \left( \lambda_{\bm{v}}^{i} \right)_k}.
\end{equation}
is greater than $1 - \epsilon_\mathrm{tol},$ for some prescribed tolerance $\epsilon_\mathrm{tol}.$ 
Similarly, we can perform the POD on the snapshot matrices for the other variables, including supremizers.  The resulting reduced snapshot matrix for the state velocity is 
\begin{equation}\label{eq:3.12}
\displaystyle
\tilde{S}_{\bm{v}}^{i} \equiv {Z_{\bm{v}}^{i}}\, {\overline{S}_{\bm{v}}^{i}},
%\in \mathbb{R}^{\mathcal{N}^{v}_{h} \times N_{t}^\mathrm{POD}}.
\end{equation}
where, ${Z_{\bm{v}}^{i}} \in \mathbb{R}^{\mathcal{N}^{v}_{h} \times N_{t}^\mathrm{POD}}$ denotes the reduced spatial basis, capturing the dominant spatial modes, and ${\overline{S}_{\bm{v}}^{i} } \in \mathbb{R}^{ N_{t}^\mathrm{POD} \times N_{t}}$ represents the temporal coefficients that encode the time evolution of these spatial modes.

Similarly, we can rewrite the compressed snapshot matrices for the other variables, including supremizers, as follows: 
\begin{equation}\label{eq:3.13}
\displaystyle
\tilde{S}_{\bm{\square}}^{i} \equiv {Z_{\bm{\square}}^{i}} \,  {\overline{S}_{\bm{\square}}^{i}},
% \in \mathbb{R}^{{\mathcal{N}^{\square}_{h}} \times N_{t}^\mathrm{POD}},
\end{equation}
where $\square \in \{\bm{p}, \bm{u}, \bm{w}, \bm{q},  \bm{s}, \bm{r}\}$, and we can compute the POD singular values and vectors as before. The critical aspect of temporal compression is effectively capturing temporal dynamics with a chosen $N_{t}^\mathrm{POD} \ll N_t$ POD bases, thereby enhancing computational efficiency.
%%====================================================================%%
\subsubsection{Parametric-space compression} 
\label{sec3.2.2:Parameteric}
Once we have compressed the temporal behaviour with the first POD, we aim at reducing the dimensionality also with respect to the parametric evolution. Thus, we stack together all the compressed matrices $\tilde{S}_{\bm{v}}^{i},$ $i = 1, 2, \dots, N_\mathrm{train}$, obtained in the previous step \eqref{eq:3.12} and \eqref{eq:3.13}, obtaining the new snapshots matrix for state velocity as:
\begin{equation}\label{eq:3.14}
\displaystyle
 \mathcal{S}_{\bm{v}}^{\bm{\mu}} = \left[  \tilde{S}_{\bm{v}}^{1}, \,\,\, \tilde{S}_{\bm{v}}^{2}, \,\,\, \cdots \,\,\,, \tilde{S}_{\bm{v}}^{N_\mathrm{train}}\right] \in \mathbb{R}^{{\mathcal{N}^{v}_{h}} \times \left( N_{t}^\mathrm{POD} \times N_\mathrm{train} \right)}. 
\end{equation}
Then, we perform an additional POD on the matrix to obtain the first $N_{\bm{v}}^{\mathrm{POD}}$ (say, $N_\mathrm{max} = N_{\bm{v}}^{\mathrm{POD}}$) singular values as basis functions for $\mathcal{S}_{\bm{v}}^{\bm{\mu}}$, following the same procedure described above. 
Similarly, the compressed matrices for the other variables, including the supremizers, are obtained as follows:
\begin{equation}\label{eq:3.15}
\displaystyle
\mathcal{S}_{\bm{\square}}^{\bm{\mu}} = \left[ \tilde{S}_{\bm{\square}}^{1}, \,\,\, \tilde{S}_{\bm{\square}}^{2}, \,\,\, \cdots \,\,\, \tilde{S}_{\bm{\square}}^{N_\mathrm{train}} \right] \in \mathbb{R}^{{\mathcal{N}^{\square}_{h}} \times \left(N_{t}^\mathrm{POD} \times N_\mathrm{train} \right)},
\end{equation}
where $\square \in \{p, \bm{u}, \bm{w}, q,  \bm{s}, \bm{r}\}.$ The POD singular values and vectors for these variables can be computed in the same way as for the state velocity.

%%====================================================================%%
\subsubsection{Algebraic formulation of projection-based POD}
\label{sec3.2.3:rom_prob}
A reduced-order approximation of the velocity, pressure, and control variables is obtained by performing a Galerkin projection onto the reduced spaces, $\mathbf{V}_{rd},$ $\mathbf{P}_{rd},$ and $\mathcal{U}_{rd}$. The approximations are sought in the following form:
\begin{equation}\label{eq:3.3.1}
\centering
\begin{gathered}
\displaystyle
 \bm{v}_h \left(t;\bm{\mu} \right)  \approx Z_{\bm{v}, \bm{s}}\,\bm{v}_{rd} \, \left(t;\bm{\mu} \right) ,
 \quad \bm{p}_h \left(t;\bm{\mu} \right)  \approx Z_{\bm{p}}\,\bm{p}_{rd} \, \left(t;\bm{\mu} \right) , 
 \quad \bm{u}_h \left(t;\bm{\mu} \right)  \approx Z_{\bm{u}}\,\bm{u}_{rd} \,\left(t;\bm{\mu} \right) , \\
\bm{w}_h \left(t;\bm{\mu} \right)  \approx Z_{\bm{w}, \bm{r}}\,\bm{w}_{rd} \, \left(t;\bm{\mu} \right), 
 \quad \bm{q}_h \left(t;\bm{\mu} \right) \approx Z_{\bm{q}}\,\bm{q}_{rd} \, \left(t;\bm{\mu} \right) , 
 \end{gathered}
\end{equation}
where, $Z_{\bm{v}, \bm{s}}$, $Z_{\bm{p}}$, $Z_{\bm{w}, \bm{r}}$, $Z_{\bm{q}}$,and $Z_{\bm{u}}$, respectively, represent the reduced basis matrices for state velocity, state pressure, adjoint velocity, adjoint pressure, and control along with the supremizers.  The resulting algebraic formulation of the reduced-order system for any $\bm{\mu} \in \mathcal{D}$ is given by:
\begin{equation*}%\label{eq:3.3.2}
\resizebox{\textwidth}{!}{$ % Start resizebox
\begin{aligned}
\begin{bmatrix}
\displaystyle
\left( \frac{1}{\Delta t}M^{rd} + M_d^{rd} \right) \left(t^{\left(n+1 \right)}; \bm{\mu} \right)  & 0 & 0 & \left(A_{ad}^{rd} + E_{ad}^{rd} \left(\bm{v}_{rd} \right) \right) \left(t^{\left(n+1 \right)}; \bm{\mu} \right)  & B_{ad}^{rd} \left(t^{\left(n+1 \right)}; \bm{\mu} \right)\\
0 & 0 & 0  & {B_{ad}^{rd}}^\top \left(t^{\left(n+1 \right)}; \bm{\mu} \right) & 0 \\
0 & 0 & \mathcal{R}^{rd}\left(t^{\left(n+1\right)}; \bm{\mu}\right)  &  0  & 0 \\
\left(A^{rd} + E^{rd} \left(\bm{w}_{rd} \right) + {E_{ad}^{rd}}^\top \left(\bm{w}_{rd} \right) \right) \left(t^{\left(n+1 \right)}; \bm{\mu} \right) & B^{rd}\left(t^{\left(n+1 \right)}; \bm{\mu} \right) &  C^{rd} \left(t^{\left(n+1 \right)}; \bm{\mu} \right) & \frac{-1}{\Delta t} {M^{rd}}^\top \left(t^{\left(n+1 \right)}; \bm{\mu} \right) & 0 \\
{B^{rd}}^\top(t^{(n+1)}; \bm{\mu}) & 0 & 0  & 0 & 0 \\
\end{bmatrix}
\begin{bmatrix}
\bm{v}_{rd} \left(t^{\left(n+1 \right)}; \bm{\mu} \right)\\
p_{rd} \left(t^{\left(n+1 \right)}; \bm{\mu} \right) \\
\bm{u}_{rd} \left(t^{ \left(n+1 \right)}; \bm{\mu} \right)\\
\bm{w}_{rd} \left(t^{\left(n+1 \right)}; \bm{\mu}\right)\\
q_{rd} \left(t^{ \left(n+1 \right)}; \bm{\mu} \right)\\
\end{bmatrix}
\\  =
\begin{bmatrix}
\bm{F}^{rd} \left(t^{(n+1)};\bm{\mu} \right) \\
0 \\
\bm{H}^{rd} \left(t^{(n+1)};\bm{\mu} \right) \\
\bm{G}^{rd} \left(t^{(n+1)};\bm{\mu} \right) \\
0
\end{bmatrix}
+
\begin{bmatrix}
\frac{1}{\Delta t}M^{rd} \left(t^n; \bm{\mu} \right)  & 0 & 0 & 0 & 0\\
0 & 0 & 0  & 0 & 0 \\
0 & 0 & 0  & 0 & 0 \\
0 & 0 & 0  & \frac{-1}{\Delta t}{M^{rd}}^\top \left(t^n; \bm{\mu} \right)  & 0 \\
0 & 0 & 0  & 0 & 0 \\
\end{bmatrix}
\begin{bmatrix}
\bm{v}_{rd} \left(t^n; \bm{\mu} \right) \\
p_{rd} \left(t^n; \bm{\mu} \right)  \\
\bm{u}_{rd} \left(t^n; \bm{\mu} \right) \\
\bm{w}_{rd}\left(t^n; \bm{\mu} \right)  \\
q_{rd}\left(t^n; \bm{\mu} \right)  
\end{bmatrix}
\end{aligned}
$} % End resizebox
\end{equation*}
Wherein, the reduced order matrices are written as:
\begin{align}\label{eq:3.3.3}
% \begin{cases}
\begin{gathered}
\displaystyle
A^{rd}\left(t; \bm{\mu} \right)  = Z_{\bm{w}, \bm{r}}^\top \, A\left(t; \bm{\mu} \right) \,  Z_{\bm{v}, \bm{s}},
 \qquad M^{rd}\left(t; \bm{\mu} \right)  =  Z_{\bm{v}, \bm{s}}^\top\, M\left(t; \bm{\mu} \right) \, Z_{\bm{v}, \bm{s}},
 \qquad E^{rd}\left(\cdot ; t; \bm{\mu} \right)  =  Z_{\bm{w}, \bm{r}}^\top\, E\left(\cdot ; t; \bm{\mu} \right)\, Z_{\bm{v}, \bm{s}},\\
  M_d^{rd}\left(t; \bm{\mu} \right)  =  Z_{\bm{v}, \bm{s}}^\top M_d\left(t; \bm{\mu} \right)  Z_{\bm{v}, \bm{s}}, 
 \qquad C^{rd}\left(t; \bm{\mu} \right) =  Z_{\bm{u}}^\top \, C\left(t; \bm{\mu} \right) \,  Z_{\bm{u}}, 
 \qquad B^{rd}\left(t; \bm{\mu} \right)  =   Z_{\bm{q}}^\top \, B\left(t; \bm{\mu} \right) \,  Z_{\bm{v}, \bm{s}},\\
  A_{ad}^{rd}\left(t; \bm{\mu} \right) = Z_{\bm{v}, \bm{s}}^\top \, A_{ad} \left(t; \bm{\mu} \right)\,  Z_{\bm{w}, \bm{r}},
 \qquad \mathcal{R}^{rd} \left(t;\bm{\mu} \right) =   Z_{\bm{u}}^\top \, \mathcal{R} \left(t; \bm{\mu} \right)\,  Z_{\bm{u}},
 \qquad E_{ad}^{rd} \left(\cdot ; t; \bm{\mu} \right)  =  Z_{\bm{v}, \bm{s}}^\top\, E_{ad} \left(\cdot ; t; \bm{\mu} \right)\, Z_{\bm{w}, \bm{r}},\\
  B_{ad}^{rd}\left(t;\bm{\mu} \right) =   Z_{\bm{p}}^\top \, B\left(t; \bm{\mu} \right)\,  Z_{\bm{w}, \bm{r}},
 \qquad F^{rd}\left(t;\bm{\mu}\right) =  Z_{\bm{v}, \bm{s}}^\top\, F\left(t;\bm{\mu} \right),
 \qquad G^{rd}\left(t;\bm{\mu} \right) =  Z_{\bm{v}, \bm{s}}^\top\, G\left(t;\bm{\mu}\right), \\
  H^{rd}\left(t;\bm{\mu} \right) = Z_{\bm{u}}^\top\, H\left(t;\bm{\mu} \right).
\end{gathered}
% \end{cases}
\end{align}
The computation of the proposed reduced order optimal flow control are performed using \textit{multiphenics}~\cite{multiphenics} for high-fidelity solutions and RBniCS~\cite{rbnics} for reduced-order solutions. Both libraries are based on FEniCS~\cite{fenics-github}. The \textit{multiphenics} library utilizes PETSc~\cite{petsc} to efficiently solve matrices, employing a block-structured formulation that is specifically designed to address the complexities typical of optimal control problems.
%%-------------------------------------------------------%%
%%%%%%%%%%%%%%%%%%%%%%%%%%%%%%%%%%%%%%%%%%%%%%%%%%%%%%%%%%%%%%%%%%%%%%%%%%%%%%%%%%%%%%%%%%
\section{Results}%{Numerical Results and Discussion}
\label{sec4:Results}

In this section, we present and analyze the results of the numerical test cases for the optimization problem \eqref{eq:2.2}--\eqref{eq:2.3}, with parametrized inlet flow profiles. We show the efficiency of the projection-based reduced order methodology in capturing the essential dynamics of the flow, to generalize to unseen configurations. In particular, we examine the temporal evolution and parametric dependence of the cardiovascular flow, understanding how the Reynolds number $Re$ interacts with the control variable in driving the dynamics towards a desired configuration. 

%%--------------------------------------------------------%%
\subsection{Computational Settings and Flow Models}
\label{sec4.1:comput}
\begin{figure}[htbp]
    \centering
    \begin{subfigure}{0.45\textwidth} % Left column
        \centering
        \includegraphics[width=\linewidth]{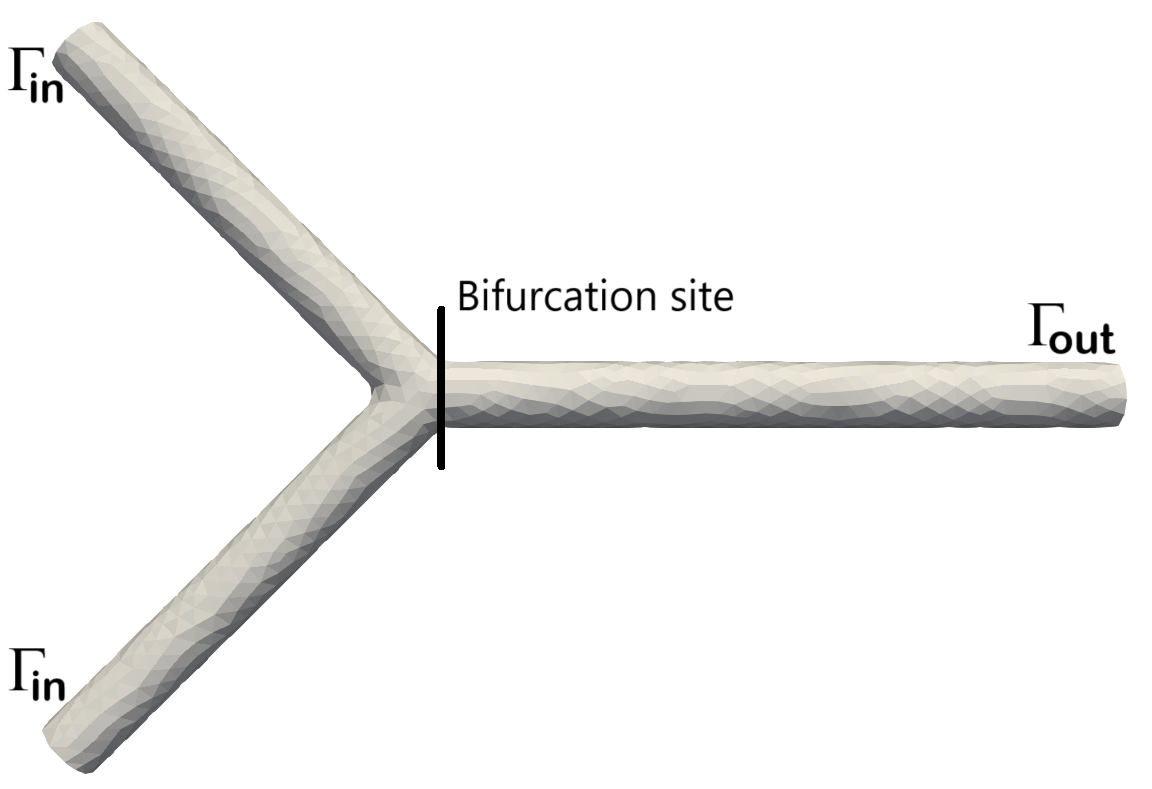}
        \caption{Test case 1: Idealistic bifurcation model}
        \label{fig1:Bif_geom}
    \end{subfigure}
    \hfill
    \begin{subfigure}{0.45\textwidth} % Right column
        \centering
        \includegraphics[width=0.85\linewidth]{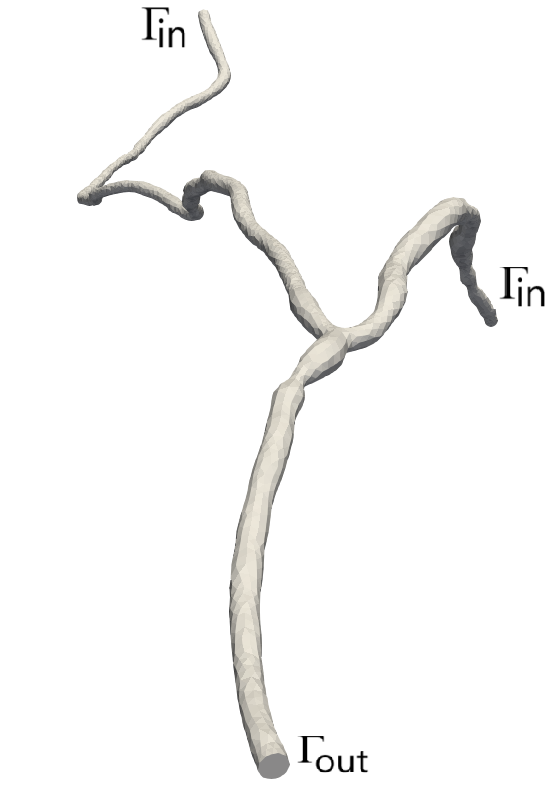}
        \caption{Test case 2: Patient-specific CABG model}
        \label{fig1:cabg}
    \end{subfigure}

    \caption{Three-dimensional vascular flow models.}
    \label{fig1:geom1}
\end{figure}

In this study, we consider two vascular flow models, namely an idealized bifurcation model and a patient-specific CABG model, as shown in Figure~\ref{fig1:geom1}. The bifurcation model consists of a simplified geometry which retains the crucial features for understanding the complex blood flow dynamics within vascular networks. This vascular model has two inlet boundaries $\Gamma_\mathrm{in}$ and one outflow boundary $\Gamma_\mathrm{out}$, where the control can act as stress boundary condition. This way, it plays a key role in both advancing our knowledge of haemodynamic flows as well as in optimizing the design of medical devices, such as vascular grafts, to enhance clinical outcomes. Figure \ref{fig1:Bif_geom} also presents a bifurcation site in the \textit{xz-plane} view, showing the junction where two inlet flows converge and exit through the outlet. 
%This view highlights the spatial arrangement of the inlets and the outlet, which is essential for understanding the flow haemodynamics.
For test case 2, the realistic geometry counterpart\footnote{Data provided by \textit{Sunnybrook Health Sciences Centre}, Toronto (Canada) as in \cite{zainib2021}.}, includes a post-surgery computed tomography (CT) scan of a 62-year-old male with a history of smoking and elevated cholesterol~\cite{zainib2021}. This patient underwent an aortocoronary bypass procedure, with the right internal mammary artery (RIMA) grafted to bypass a blockage in the left anterior descending artery (LAD), resulting in a coronary artery bypass graft (CABG), as depicted in Figure~\ref{fig1:cabg}. For numerical computations, we employed the Taylor-Hood stable $\mathbb{P}^2 - \mathbb{P}^1$ finite element polynomials for the velocity and pressure fields, respectively, and $\mathbb{P}^2$ polynomials for the control variables. Additional details on the computational setup for both test cases are provided in Table \ref{tab1:parameters}.

Our objective is to minimize the discrepancy between a desired outflow velocity profile and the velocity profile achieved through the unsteady N--S equations across various parametrized inlet flow profiles. In particular, for both the idealized bifurcation and patient-specific CABG test cases, we chose the inlet flow profile defined on the inlet boundaries $\Gamma_\mathrm{in}$ as:
\begin{equation}\label{eq:4.1}
    \bm{v}_\mathrm{in} = -\frac{\nu \, Re}{R_{in}}\left( 1 - \frac{r^2}{R_{in}^2}\right)\,f\left(t \right)\,\bm{n}_{in},
\end{equation}
where $r$ represents the radius of cross-sections along the domain computed using centerlines obtained from VMTK~\cite{antiga2008image}, and $R_{in}$ denotes the radius of inlet cross-sections with their respective unit outward normal directions. The negative sign allows for an inflow direction, $f(t) = 0.02 + 0.02 \sin \left(\pi t \right)$ represents a time-dependent function accounting for pulsatile behaviour, and the kinematic viscosity is fixed as $\nu = 3.6$ mm$^2$/s. In blood flow problems, the non-dimensional quantity we use to parametrized the model, \textit{i.e.}\ the Reynolds number $Re$, holds significant importance due to its influence on flow field characteristics across different regimes. For the first test case, we consider $Re$ as the parameter of the model for which we aim at reducing the computational cost. In the following, it will be denoted as $\mu = Re$, and it will vary in the parametric domain  $\mathcal{D} = [50, 80]$. We used a similar inlet flow profile also for the patient-specific CABG case \ref{fig1:cabg}, as real patient-specific data for the inlet conditions was not provided. For this test case, we consider the more complex setting in which each inlet flow boundary is characterized by a Reynolds number, i.e. \ $\bm \left(\mu_1, \mu_2 \right) = \left(Re_1, Re_2 \right)$, where both $Re_1$ and $Re_2$ vary independently within the range $\mathcal{D} = [50, 80]^2$.  These simulations enable a thorough examination of how different inlet flow profiles, defined by varying Reynolds numbers, influence the flow behaviour and dynamics within the vascular models. 

To complete the optimal control framework, the target flow velocity profile $\bm{v}_d$ is chosen as: 
\begin{equation}\label{eq:4.2}
     \bm{v}_d = v_\mathrm{const}\,\left( 1 - \frac{r^2}{R(s)^2}\right)\,\bm{t}\left(s \right),
\end{equation}
where $\bm{t}\left(s \right)$ is the tangent along the centerline $\bm{c}\left(s \right)$, $R\left(s \right)$ is the maximum radius of the vessel around $\bm{c} \left(s \right)$, $r$ is the distance between the mesh nodes and nearest point lying on $\bm{c}\left(s \right)$, $s$ is the curvilinear abscissa, and $v_\mathrm{const} > 0$ represents the maximum magnitude of the  $\bm{v}_d$. As discussed in ~\cite{sankaranarayanan2005computational}, the maximum velocity of the smaller blood beds is much less than the aortic vessel and lies in the range $50-500$ mm/s; therefore, we consider $v_\mathrm{const} = 250$ mm/s for the first test case, and we set $v_\mathrm{const} = 100$ mm/s for the second one.  For this study, we heuristically choose $\alpha = 10^{-3},$ as this value provides a reasonable balance between regularization and accuracy for the computations.  Table \ref{tab1:parameters} contains the details of computation settings and physical parameters. 
\begin{table}[htbp]
\centering
\caption{Physical Parameters and Computational Details}
\label{tab:parameters}
\renewcommand{\arraystretch}{1.5} % Adjust spacing between rows
\begin{tabular}{|l|l|l|}
\hline
Computational Parameters & Test case 1 & Test case 2\\ \hline

Parametric Space & $Re_1 \in [50, 80]$ & $\left(Re_1, Re_2\right) \in [50, 80]^2$ \\ \hline
FE dofs & $\mathcal{N}_h = 225'248$ & $\mathcal{N}_h = 433'288$\\ \hline
Final time & $T=1.0$ s  & $T=0.8$ s\\ \hline
Time step & $\Delta t= 0.01$ s & $\Delta t = 0.01$ s\\ \hline
\hline
No. snapshots w.r.t.\ time $t$ & 21 & 21 \\ \hline
No. snapshots w.r.t.\ parameter $\bm{\mu}$ $\left(= N_\mathrm{train} \right)$ & 21 & 25 \\ \hline
\textit{nested-POD} $\left(= N_t^{\mathrm{POD}} \right)$ & 10 & 10 \\ \hline 
$N_\mathrm{max} \left( = N_{\bm{\square}}^{\mathrm{POD}} \right)$ &  15 & 20 \\ \hline 
Offline CPU time  &  $\approx 4-5$ hours  & $\approx 7-8$ hours \\ \hline
Online CPU time  &  $\approx  30$ minutes  & $\approx  1$ hour\\ \hline
\end{tabular}
\label{tab1:parameters}
\end{table}

Performing the unsteady simulations of OCP$_{\left(\bm{\mu} \right)}$s for the three-dimensional model, which includes state equations, adjoint equations, and cost functionals, requires significant computational resources and memory allocation, even when $\bm{\mu}$ is fixed. Managing these simulations becomes particularly challenging due to the varying computational demands associated with different $\bm{\mu}$ values. This requires efficient allocation and utilization of computational resources to ensure that the simulations run smoothly and within acceptable time frames despite the complexities involved.  To address these computational challenges effectively, we implemented a strategy involving two intermediate steps to capture spatial-temporal information of the variables over the time interval, $[0, T],$ where $T >0$ is the final time. Initially, we saved information at every $k^{th}$ time-step $\left(\Delta t \right)$ out of $N_t$, which helped mitigate memory allocation issues during transient simulations. However, conducting simulations over the entire time interval $[0, T]$ posed challenges due to resource constraints. To address this, we adopted a sequential approach, in which initially, we saved simulation results for the time interval $[0, t_1]$ with $t_1 > 0$. This step allowed us to store information up to $t_1$, effectively managing memory usage. Subsequently, we continued the simulation from the interval $[t_1, t_2]$ with $t_2 > t_1$, and we read the previously saved results at $t=t_1$, then continued the computation from $t_1 + \Delta t$. We repeated this process until reaching the final time $T$, incrementally saving and reusing simulation data between each successive time interval. This dual-step strategy optimizes memory usage and computational efficiency, ensuring that long-duration simulations are conducted effectively without compromising accuracy or computational performance, thus addressing the challenges posed by varying computational demands and resource limitations.

%%------------------
\subsection{Quantifying the Efficacy of OCPs:}
\label{sec:4.1.1:Efficacy}
\begin{figure}[htbp]
    \centering
    % Subfigure 1
    \begin{subfigure}[b]{0.325\textwidth}
        \centering
        \includegraphics[width=\linewidth]{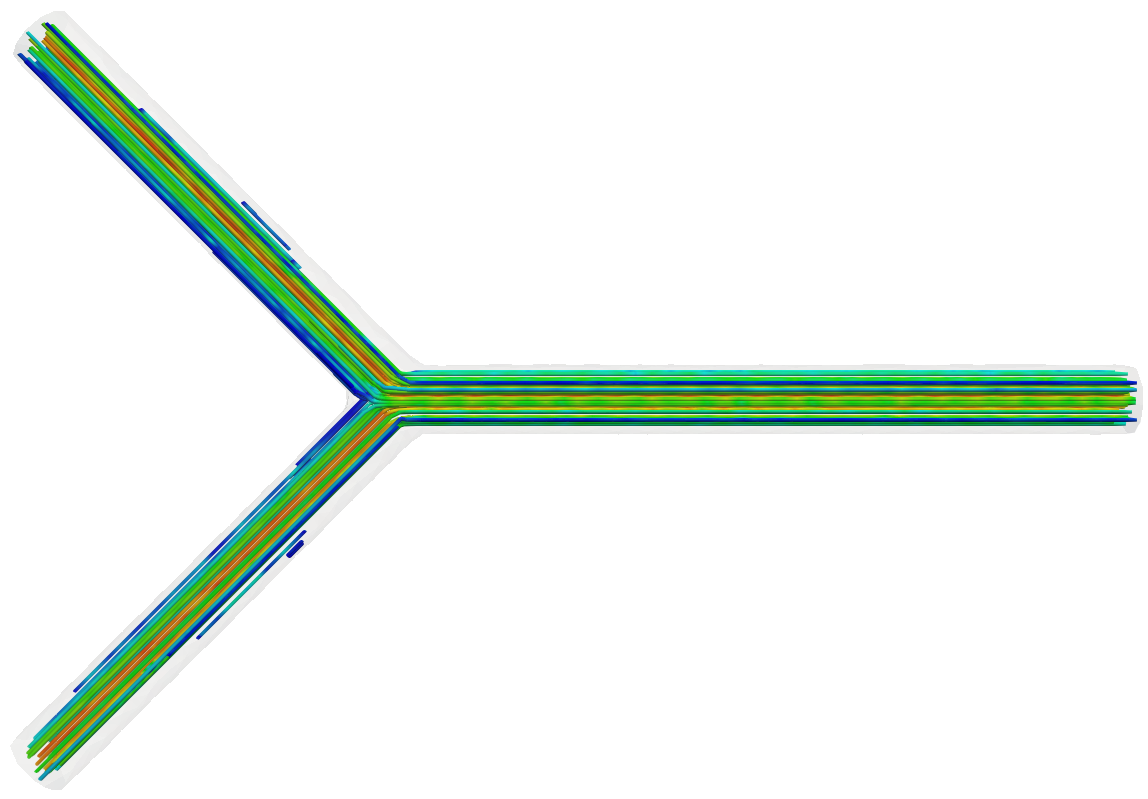}
        \caption{{Target velocity profile}}
        \label{fig1311:subfig1}
    \end{subfigure}
    \hfill
    % Subfigure 2
    \begin{subfigure}[b]{0.325\textwidth}
        \centering
        \includegraphics[width=\linewidth]{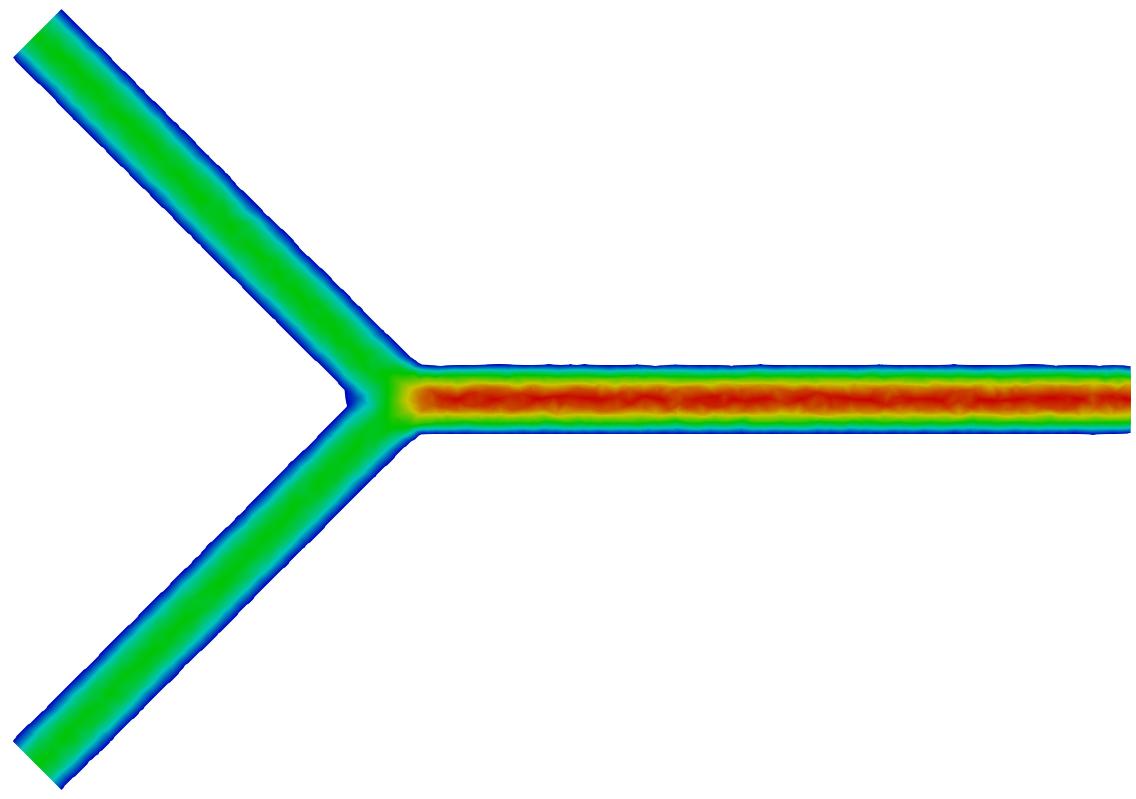}
        \caption{{Uncontrolled velocity profile}}
        \label{fig1311:subfig2}
    \end{subfigure}
    % Subfigure 3
    \begin{subfigure}[b]{0.325\textwidth}
        \centering
        \includegraphics[width=\linewidth]{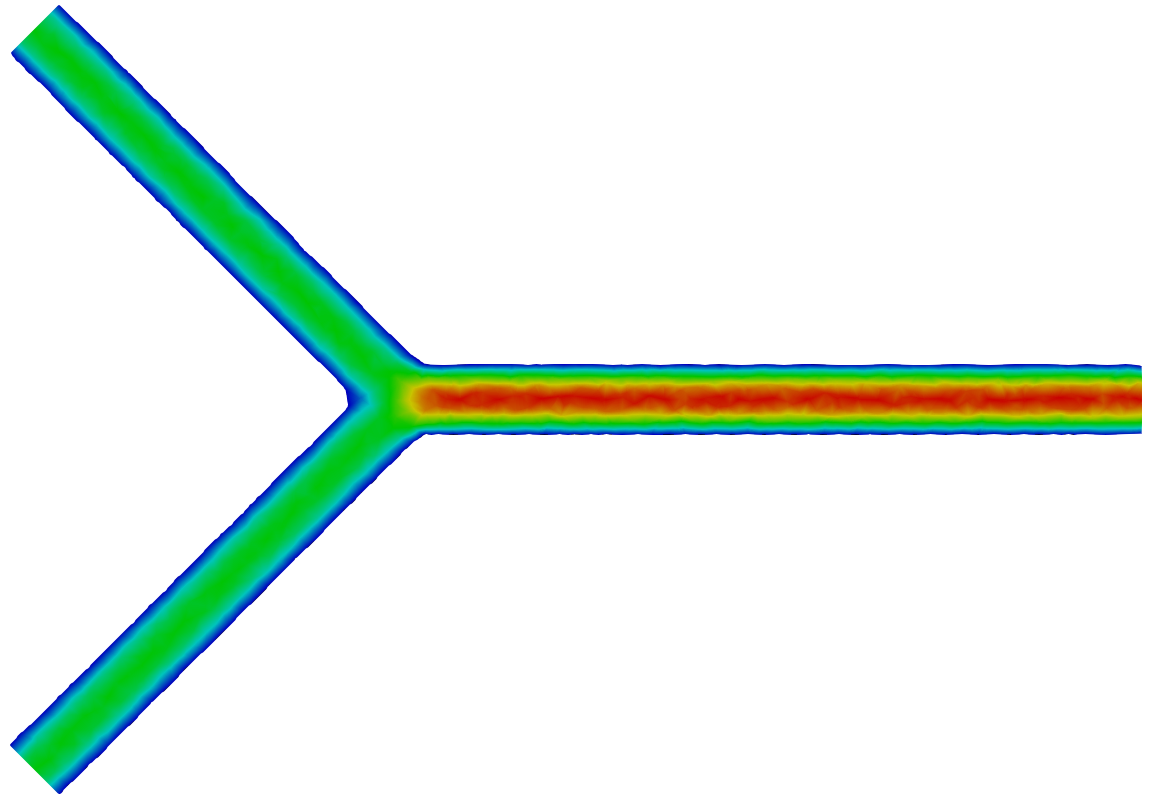}
        \caption{{Controlled velocity profile}}
        \label{fig1311:subfig3}
    \end{subfigure}
    
     %\vskip\baselineskip % Add space between rows
     
    \begin{subfigure}[b]{\textwidth}
        \centering
        \includegraphics[width=0.325\textwidth]{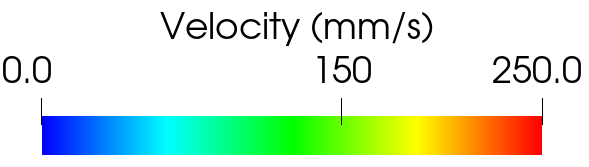}
        \label{fig1311:colorbar}
    \end{subfigure}
  \caption{{Streamlines of target, and comparison of uncontrolled and controlled high-fidelity profiles with $Re=50$ at $t=0.25$ s.}}
  \label{fig1311:discrepencies}
\end{figure}

\begin{figure}[ht]
    \centering
    % Subfigure 1
    \begin{subfigure}[b]{0.425\textwidth}
        \centering
        \includegraphics[width=0.75\linewidth]{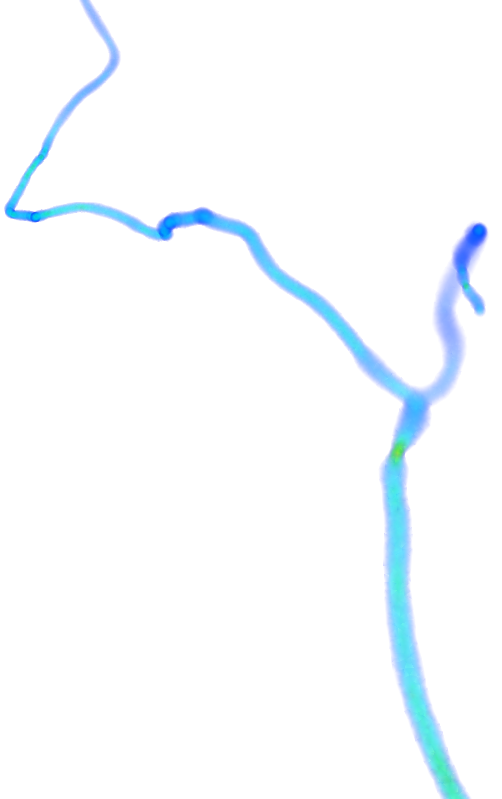}
        \caption{{Uncontrolled velocity profile}}
        \label{fig13:subfig1}
    \end{subfigure}
    \hfill
    % Subfigure 2
    \begin{subfigure}[b]{0.425\textwidth}
        \centering
        \includegraphics[width=0.75\linewidth]{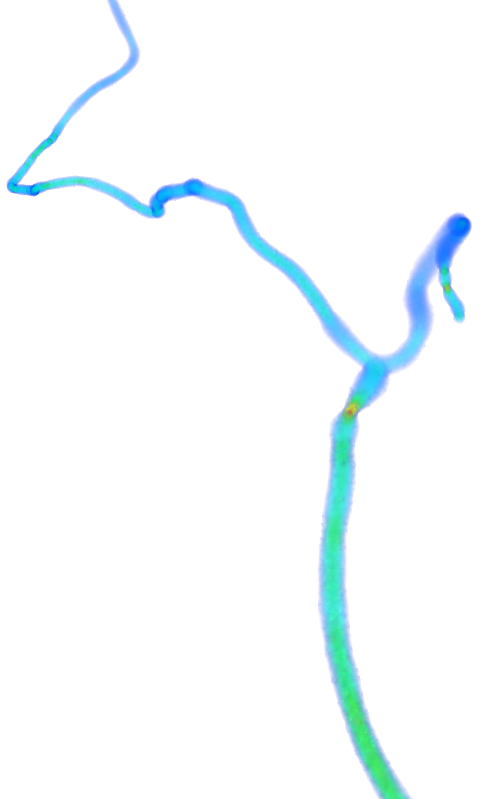}
        \caption{{Controlled velocity profile}}
        \label{fig13:subfig2}
    \end{subfigure}

    \begin{subfigure}[b]{\textwidth}
        \centering
        \includegraphics[width=0.325\textwidth]{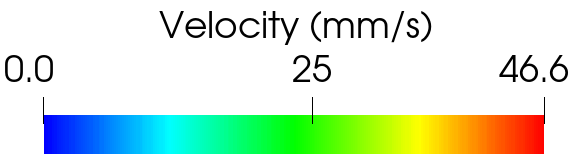}
        \label{fig13:colorbar}
    \end{subfigure}
  \caption{{Comparison of uncontrolled and controlled high-fidelity profiles with $\left(Re_1, Re_2 \right) = \left(50, 50 \right)$ at the maximum inlet.}}
  \label{fig13:discrepencies}
\end{figure}

{Figures~\ref{fig1311:discrepencies} and \ref{fig13:discrepencies} quantify the efficacy of the OCPs by comparing high-fidelity velocity profiles between uncontrolled and controlled cases within both vascular test case models. Uncontrolled flow is computed by solving the unsteady Navier–Stokes equations with a homogeneous Neumann boundary condition, \textit{i.e.} $$ \displaystyle \left( \nu \, \nabla \bm{v} \left(t; \bm{\mu} \right) - p \left(t; \bm{\mu} \right) \right) \cdot \bm{n} = 0,$$ which implies that no stress is applied at the outflow boundary. This often leads to unrealistic flow behaviour and non-physiological velocity patterns, particularly in anatomically complex regions post-surgery.}

{ We simulated the flow dynamics within the idealized bifurcation model, characterized as a symmetric geometry, rigid vessel walls, and a prescribed inlet flow profile as given in~\eqref{eq:4.1} for both the uncontrolled and controlled configurations. In controlled configuration, the outlet pressure control, imposed as a \textit{non-homogeneous Neumann} boundary condition mentioned in the last equation of equation~\eqref{eq:2.3}, indicating that outlet traction has a negligible impact on the flow fields presented in Figure~\ref{fig1311:subfig3}, which appear nearly identical with the uncontrolled one, as shown in Figure~\ref{fig1311:subfig2}. 
These high-fidelity solutions are similar to $Re = 50$ at $t = 0.25$ s, due to the consideration of a symmetric target profile, as illustrated in Figure~\ref{fig1311:subfig1}, which is derived analytically using the steady Poiseuille flow assumption~\eqref{eq:4.2}, providing a reference for assessing control accuracy and model performance \cite{sutera1993history}. This situation is similar to the single-vessel scenario, where, assuming negligible viscosity, the velocity field and pressure drop are uniquely determined by the inlet velocity and rigid vessel walls. In this case, the outlet pressure does not affect the velocity field, as it is already fully determined by the inflow and geometry. Since a control variable must influence the velocity field to reduce the discrepancy with the target profile, the outlet pressure cannot serve as an effective control variable for the symmetric test cases.}
In contrast, the patient-specific CABG geometry, characterized as an asymmetric configuration, exhibits a substantial difference that is observed between the uncontrolled and controlled velocity fields, as shown in Figures~\ref{fig13:subfig1} and~\ref{fig13:subfig2}. This highlights the significance of the control strategy in achieving a more accurate representation of CV flow, ensuring that the simulated dynamics closely align with physiological conditions. This also underscores the critical role of optimal control strategies in enforcing physiologically meaningful, boundary-driven behaviour, thereby improving the accuracy of flow predictions; such strategies enhance the model's ability to reproduce realistic hemodynamics, supporting personalized modelling and optimized surgical decision-making.
%%
%%============================================================================%%
\subsection{Test case 1: Idealistic Bifurcation Model}
\label{sec4.1:testcase1}
In this test case, we consider an idealized bifurcation model to investigate the underlying flow dynamics and assess the performance of the proposed computational framework in a controlled setting. This idealistic geometry enables the analysis of essential hemodynamics, such as velocity gradients and pressure variations, while minimizing anatomical complexity. %The target velocity profile is derived analytically using the steady Poiseuille flow assumption~\eqref{eq:4.2}, providing a reference for assessing control accuracy and model performance \cite{sutera1993history}. % {Figure~\ref{fig1311:subfig1} shows the streamlines of the target profile and high-fidelity solution computed using the state equations, with a controlled outflow boundary and $Re = 50$ at $t = 0.25$ s.} This high-fidelity simulation captures the intricate fluid behaviour at the bifurcation, including the effects of unsteady and viscous forces. %Notably, the outlet flow velocity is higher than the inlet velocities, as mentioned in~\eqref{eq:4.1}, a deviation from the Poiseuille flow assumption. This increase in outlet velocity arises from the convergence of flow at the bifurcation, where the merging streams must accelerate to conserve mass, as described by the continuity equation $Q_\mathrm{inlet} \equiv Q_\mathrm{outlet}$. This highlights the importance of detailed computational models for accurately capturing intricate fluid dynamics within complex vascular flow geometries.

%%--------------------------------------------------------%%
\subsubsection{Assessment of ROM Performance}
\begin{figure}[htbp]
    \centering
    % Subfigure 1
    \begin{subfigure}{0.475\textwidth}
        \centering
        \includegraphics[width=\linewidth]{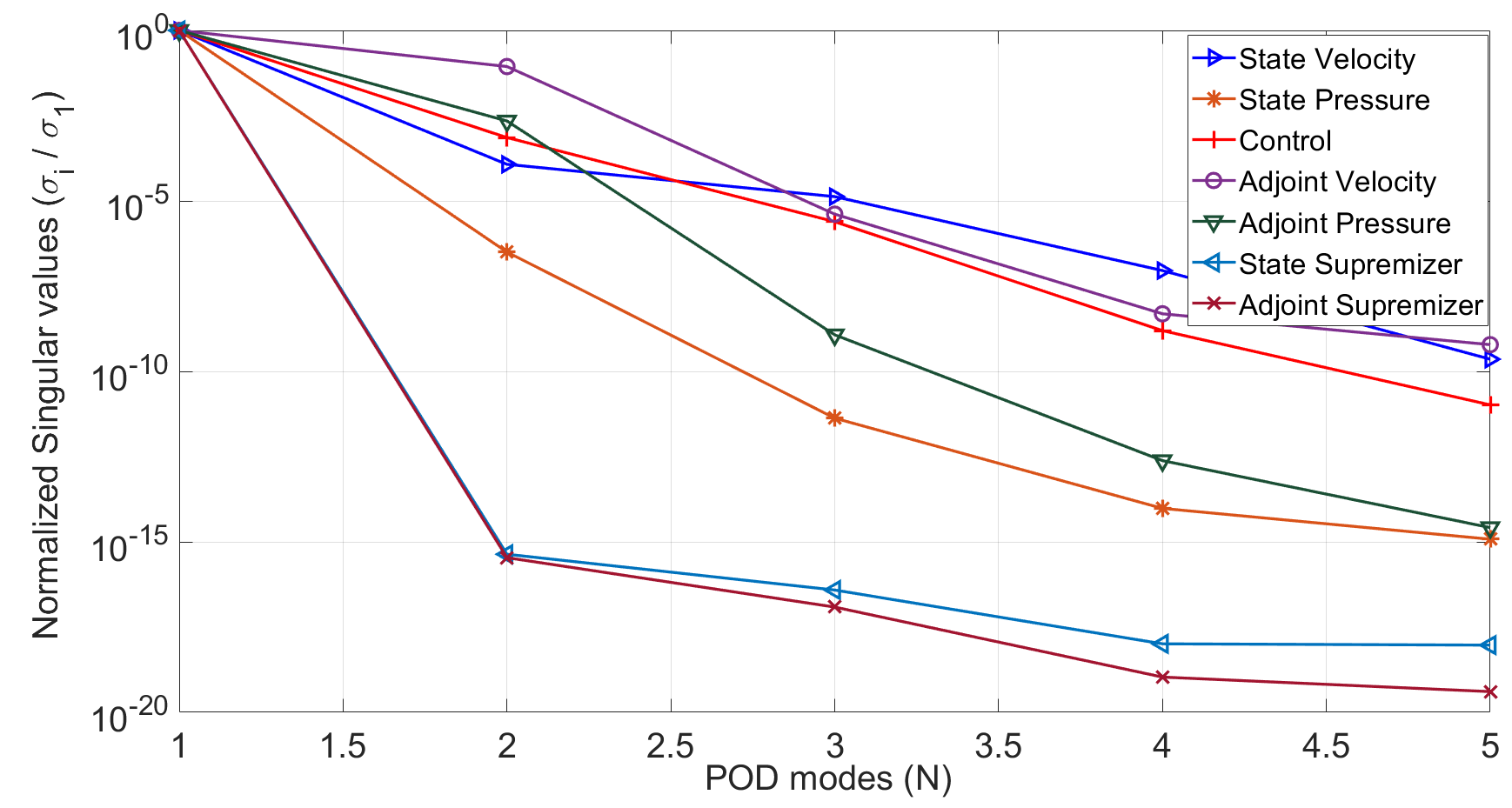}
        \caption{Singular values:  \textit{temporal compression}}
        \label{fig2:subfig1}
    \end{subfigure}
    \hfill
    % Subfigure 2
    \begin{subfigure}{0.475\textwidth}
        \centering
        \includegraphics[width=\linewidth]{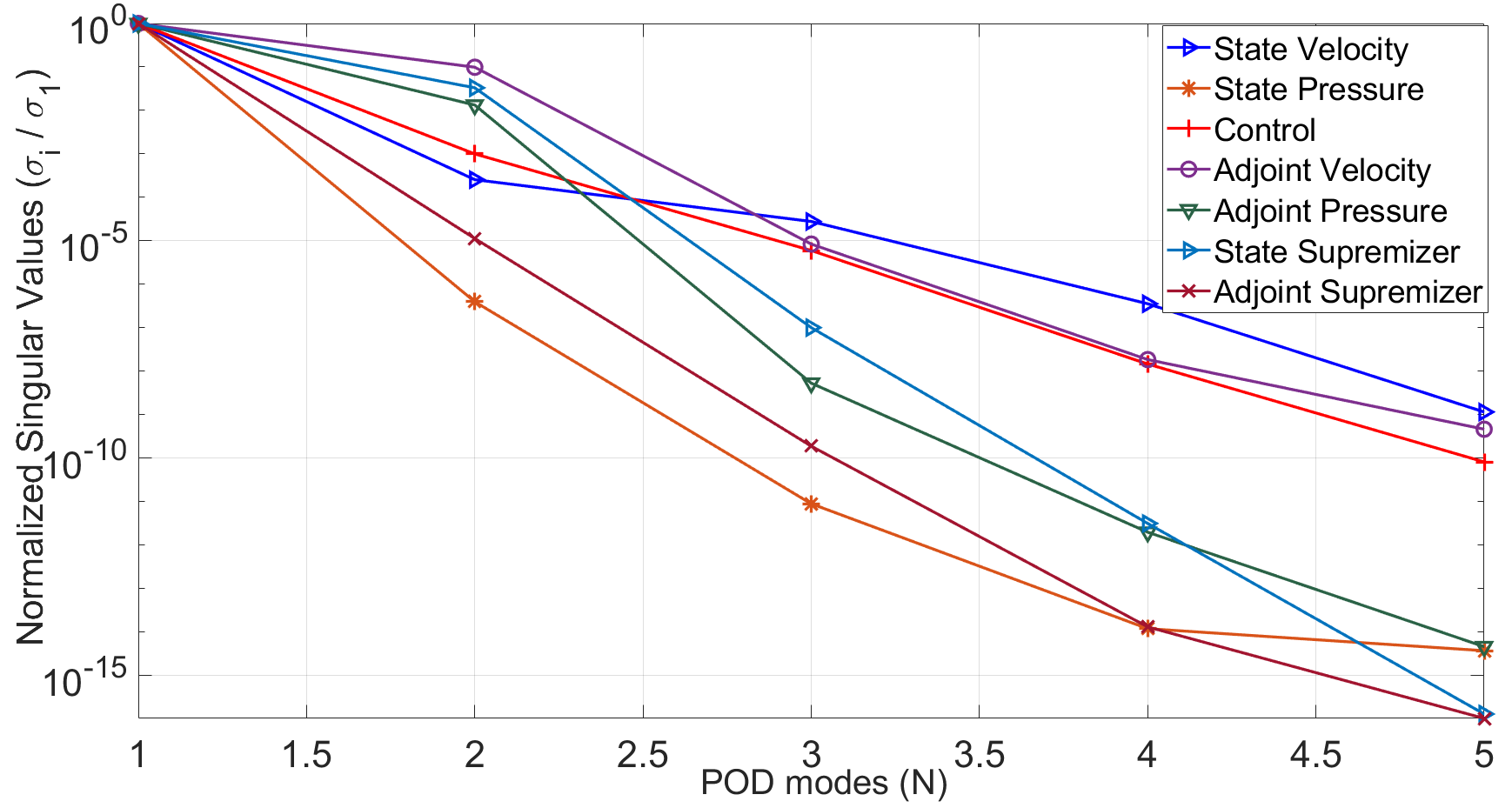}
        \caption{Singular values: \textit{parametric-space compression}}
        \label{fig2:subfig2}
    \end{subfigure}
  \caption{{Normalized POD singular values for velocity, pressure, and supremizers {for $Re = 50$}.}}
  \label{fig2:eigen_bif}
\end{figure}

%\subsection{ Test case 1: Idealistic Vascular Model}
\label{sec4.2.1:bif_rom}
As mentioned previously in Section~\ref{sec4.1:comput}, we collect only $21$ snapshots in time out of the 100 time-steps solutions, for all the variables. Thus, by subsampling the time domain, we are able to ensure the accuracy utilizing a small time-step, while addressing the memory issue. In the context of \textit{nested-POD}, for the \textit{temporal compression}, we selected $N_{t}^{\mathrm{POD}}=10$ modes from these 21 snapshots, for each of the $N_\mathrm{train}=21$ parameter samples. Therefore, the resulting matrix size is $210 \times \mathcal{N}_h$, which is significantly smaller than the matrix size $2100 \times \mathcal{N}_h$ used in classical POD. This $90\%$ reduction not only manages memory but also accelerates computations, particularly in large-scale settings with $\mathcal{N}_h = 225,248$, making the nested-POD approach both efficient and accurate. The offline CPU time for a fixed $\bm{\mu}$ snapshot was approximately $4-5$ hours, leading to a total offline time of around $4$ days for all $21$ snapshots. In contrast, the online CPU time was significantly reduced, taking around $30$ minutes for a snapshot, demonstrating the efficiency over the offline. Figure~\ref{fig2:subfig1} shows a consistent decay of eigenvalues for state, adjoint, and control variables, however, the supremizers exhibit a rapid decay. A decay of $\mathcal{O}\left(10^{-10}\right)$ in the POD singular values is achieved with fewer than $10$ POD modes, indicating a uniform representation of system dynamics during \textit{temporal compression}. From Figure~\ref{fig2:subfig2}, we observe that pressure and supremizers show a rapid decay in normalized singular values compared to the velocity and control; they also show a significant decay, indicating that the considered system dynamics can be efficiently represented by a few modes, $N \leq 10$. Hence, Figure~\ref{fig2:eigen_bif} demonstrates the efficiency of \textit{nested-POD} in capturing the essential spatiotemporal features of the system with a minimal number of modes while preserving its dynamics. The rapid decay of singular values suggests a significant reduction in computational cost without affecting accuracy. 
\begin{figure}[htbp]
    \centering
    \begin{subfigure}[b]{0.325\textwidth}
        \centering
        \includegraphics[width=\textwidth]{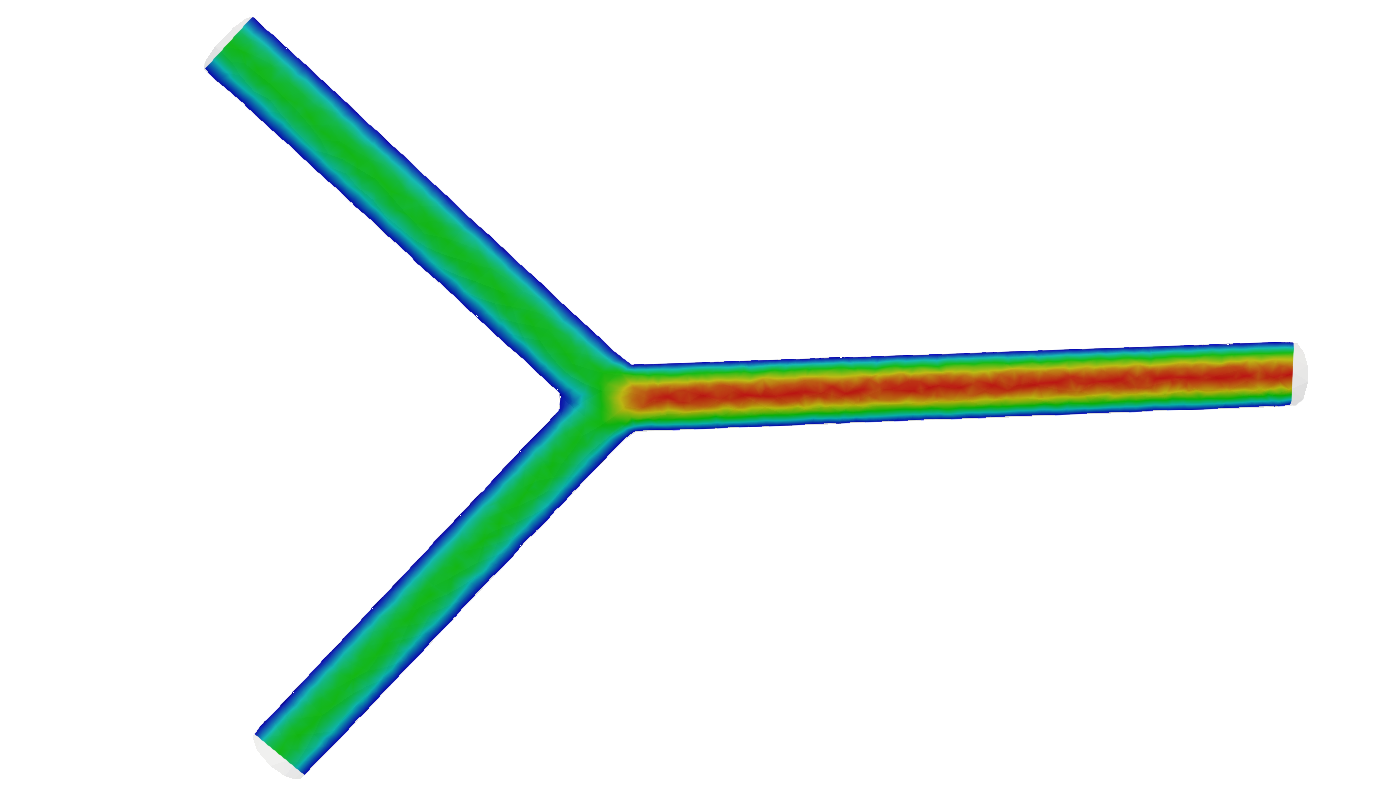}
        \caption{\textit{High-fidelity} solution}
        \label{fig3:fom_vel}
    \end{subfigure}%
    \hfill
    \begin{subfigure}[b]{0.325\textwidth}
        \centering
        \includegraphics[width=\textwidth]{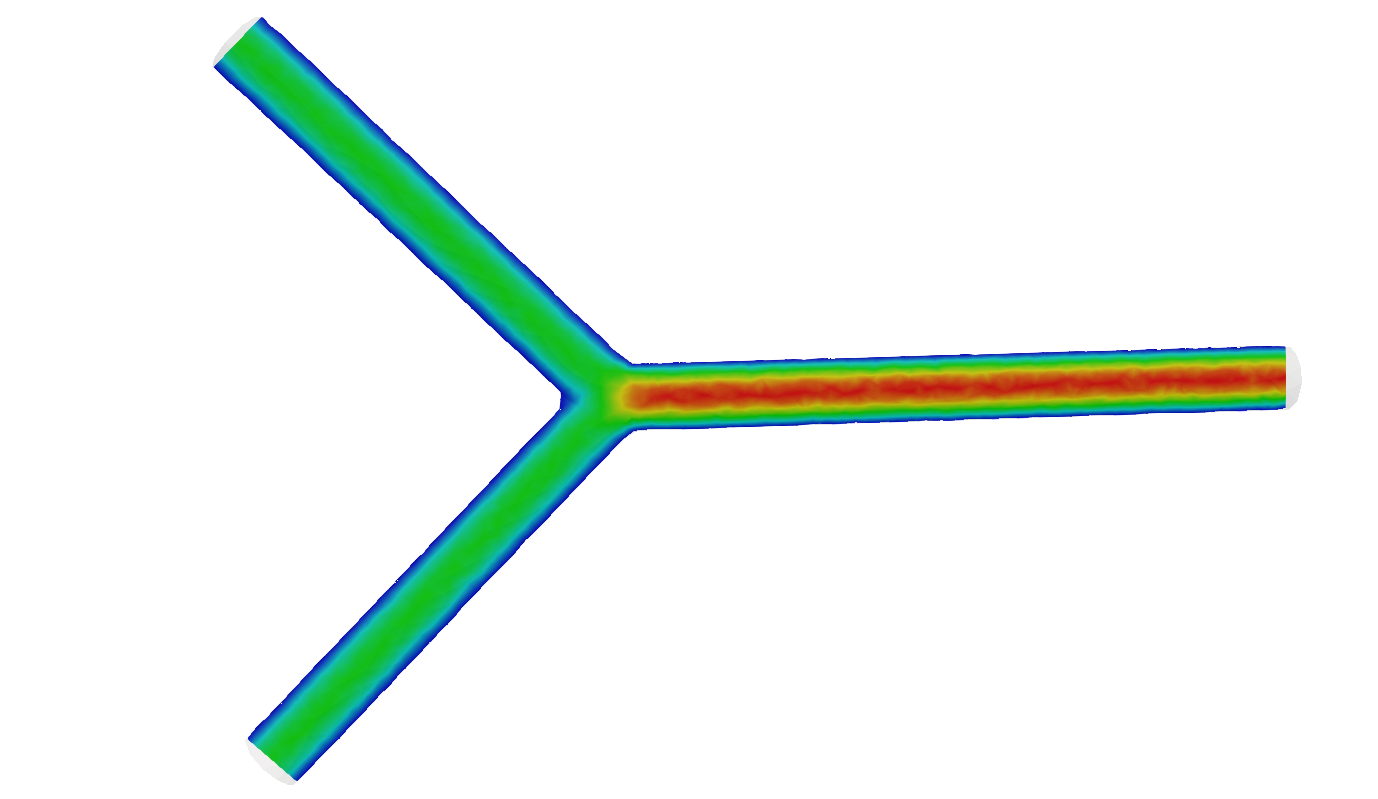}
        \caption{\textit{Reduced-order} solution}
        \label{fig3:rom_vel}
    \end{subfigure}%
    \hfill
    \begin{subfigure}[b]{0.325\textwidth}
        \centering
        \includegraphics[width=\textwidth]{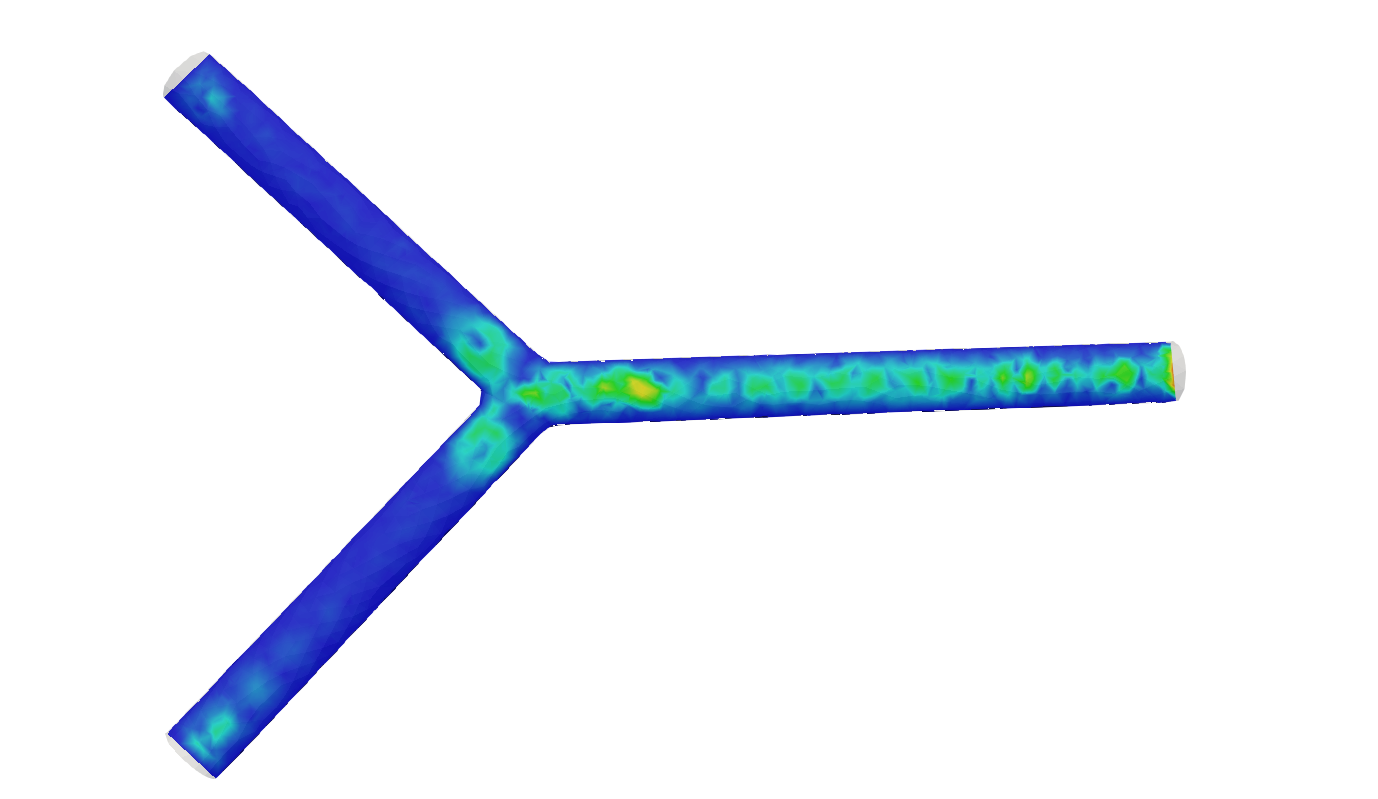}
        \caption{Absolute error}
        \label{fig3:error_vel}
    \end{subfigure}
    %\vskip\baselineskip % Add space between rows
    \begin{subfigure}[b]{0.65\textwidth}
        \centering
        \includegraphics[width=0.5\textwidth]{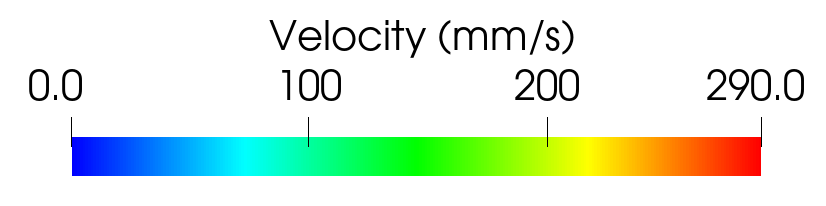}
        \label{fig3:velocity_colorbar}
    \end{subfigure}
    \hfill
    \begin{subfigure}[b]{0.325\textwidth}
        \centering
        \includegraphics[width=0.85\textwidth]{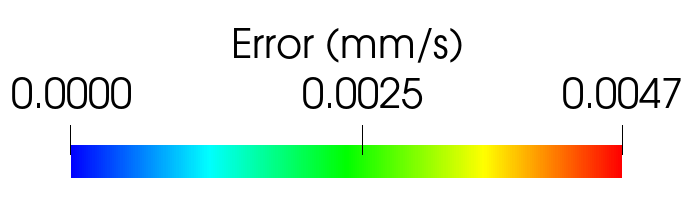}
        \label{fig3:error_colorbar}
    \end{subfigure}
    \caption{Comparison of velocity profiles at $Re = 50$ and $t = 0.5$ s. (a) High-fidelity solution, (b) Reduced-order solution, and (c) Absolute error.}
    \label{fig3:Vel_comp_bif}
\end{figure}

\begin{figure}[htbp]
    \centering
    \begin{subfigure}[b]{0.325\textwidth}
        \centering
        \includegraphics[width=\textwidth]{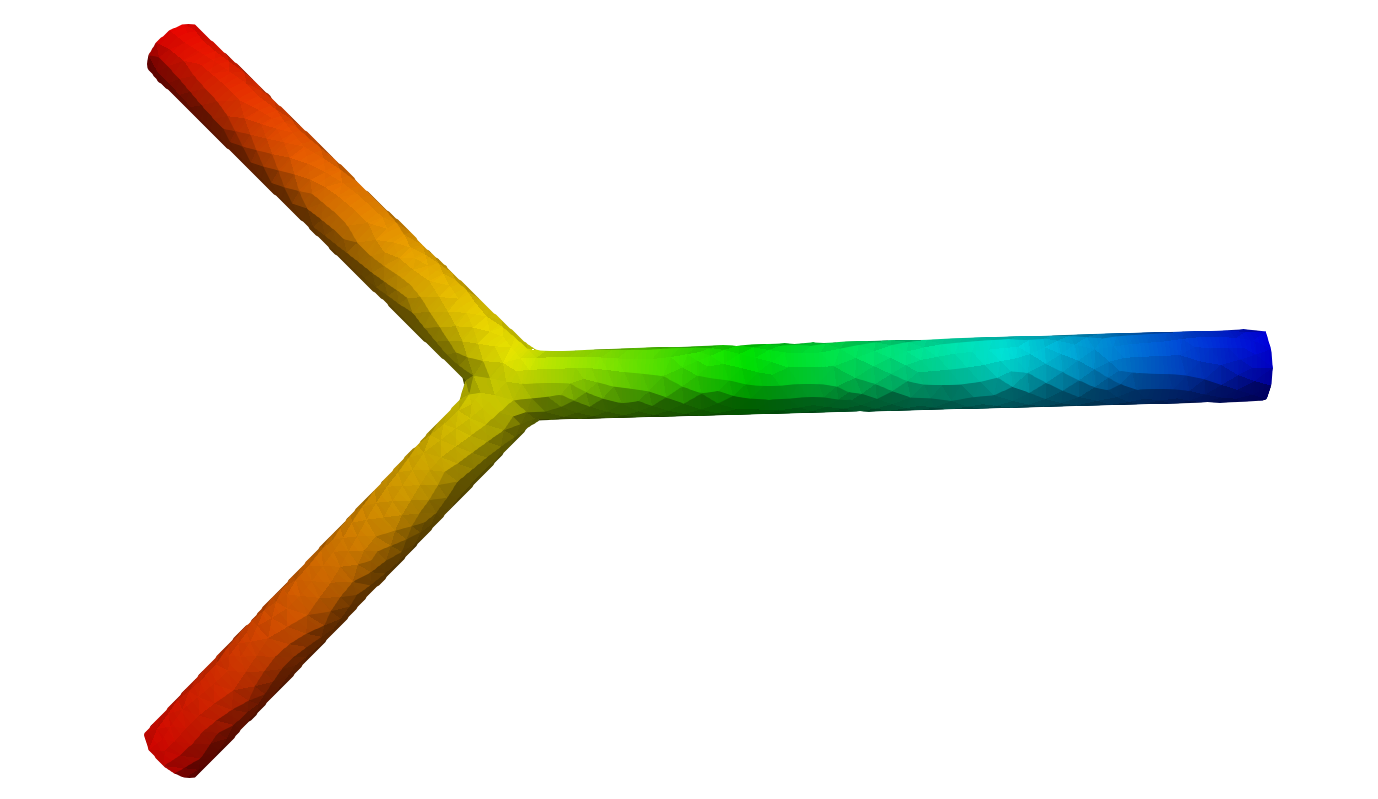}
        \caption{\textit{High-fidelity} solution}
        \label{fig4:fom_press}
    \end{subfigure}%
    \hfill
    \begin{subfigure}[b]{0.325\textwidth}
        \centering
        \includegraphics[width=\textwidth]{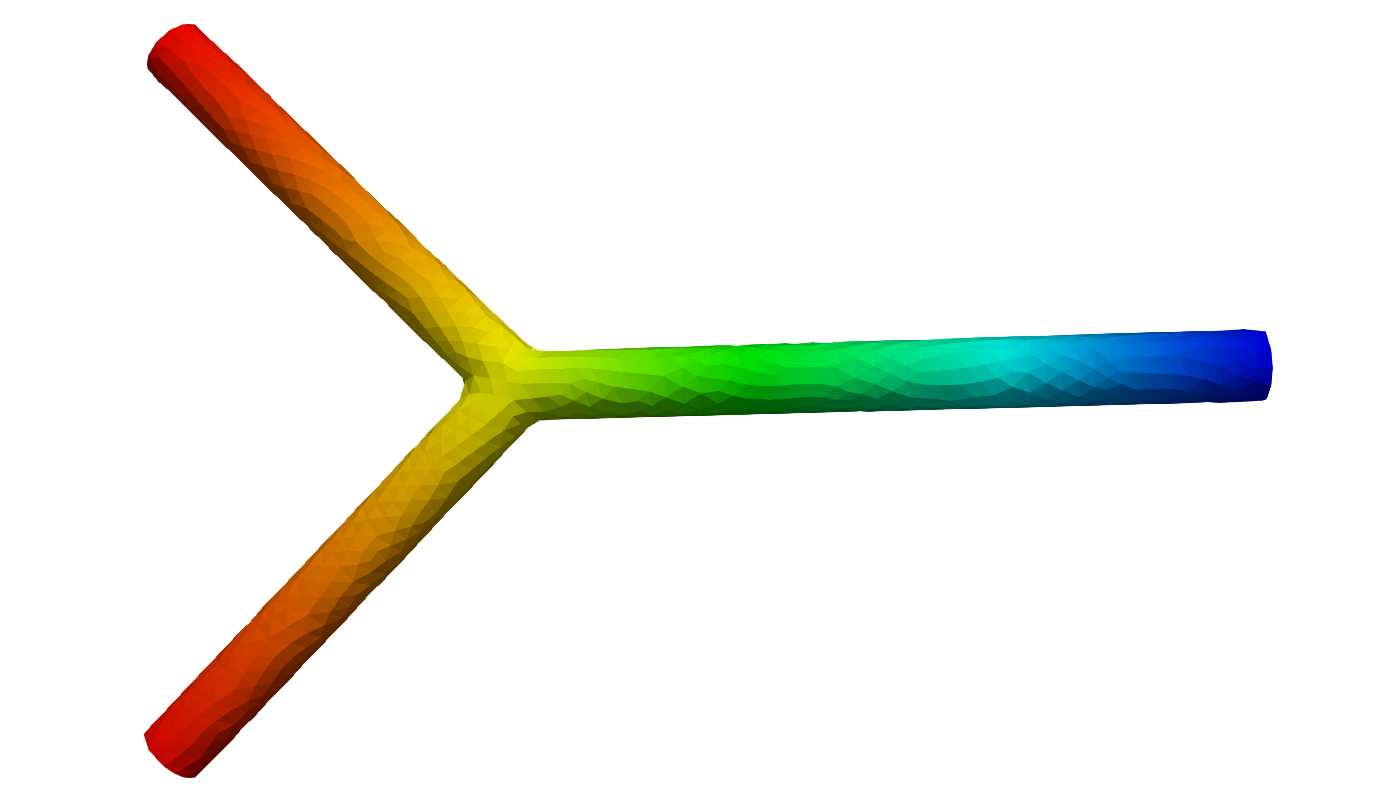}
        \caption{\textit{Reduced-order} solution}
        \label{fig4:rom_press}
    \end{subfigure}%
    \hfill
    \begin{subfigure}[b]{0.3\textwidth}
        \centering
        \includegraphics[width=\textwidth]{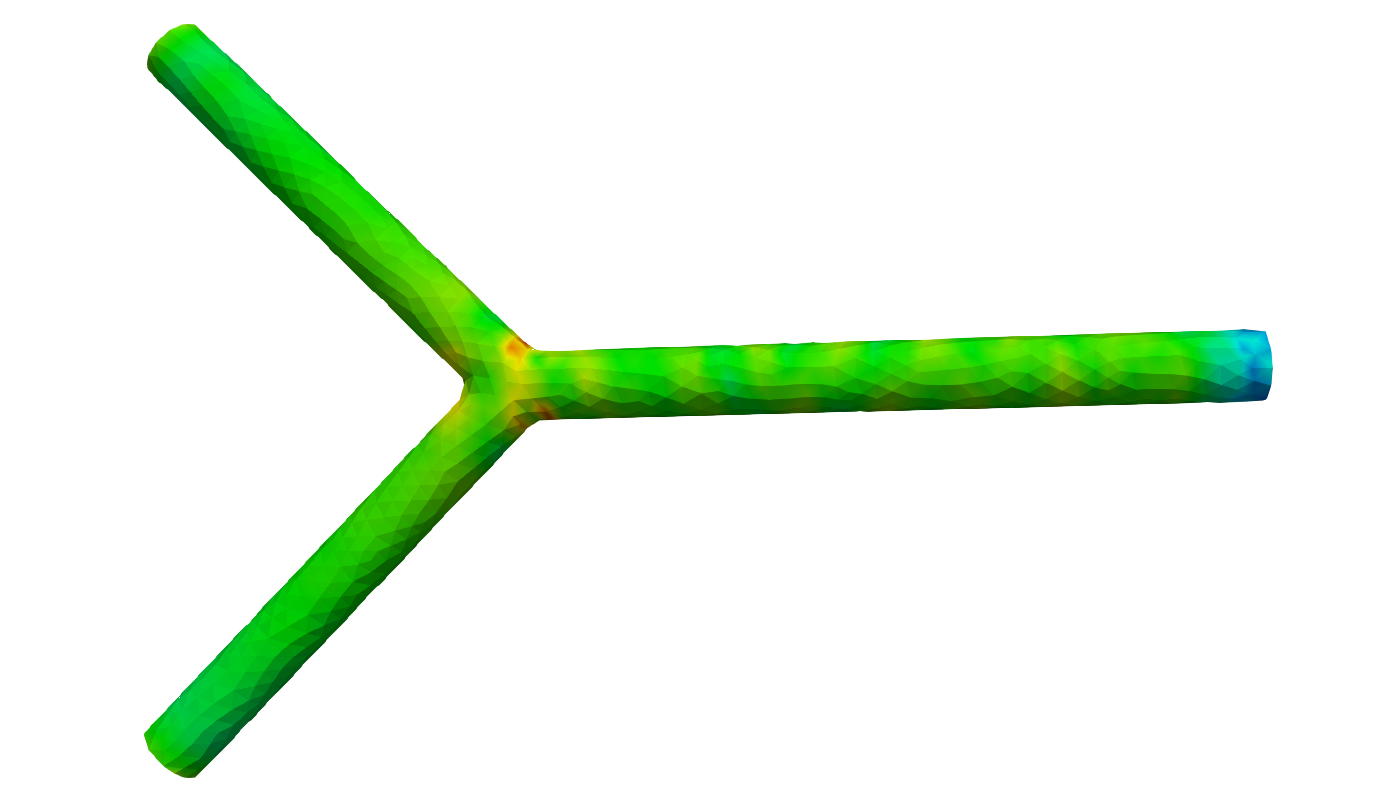}
        \caption{Absolute error}
        \label{fig4:error_press}
    \end{subfigure}
    
    %\vskip\baselineskip % Add space between rows
    
    \begin{subfigure}[b]{0.65\textwidth}
        \centering
        \includegraphics[width=0.52\textwidth]{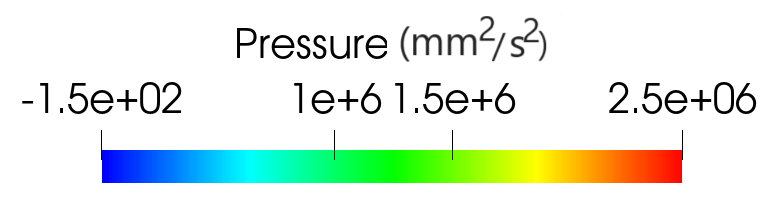}
        \label{fig4:press_colorbar}
    \end{subfigure}
    \hfill
    \begin{subfigure}[b]{0.325\textwidth}
        \centering
        \includegraphics[width=0.95\textwidth]{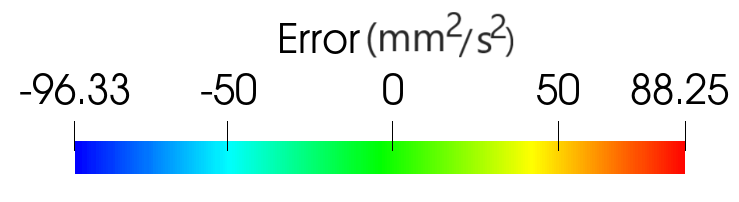}
        \label{fig4:error_colorbar_press}
    \end{subfigure}
    
    \caption{Comparison of pressure distribution at $Re = 50$ and $t = 0.5$ s. (a) High-fidelity solution, (b) Reduced-order solution, and (c) Absolute error.}
    \label{fig4:press_comp_bif}
\end{figure}
Figure~\ref{fig3:Vel_comp_bif} presents the high-fidelity and reduced-order velocity profiles at maximum inlet flow profile with $Re=50$ along with the absolute error between the solutions. We have observed that the inlet flows from the inlet branches merge at the bifurcation site and exit smoothly through the outlet, as depicted in the high-fidelity solution in Figure~\ref{fig3:fom_vel}. The fluid moves faster along the centerline, forming a parabolic velocity profile, while near the walls, the velocity gradually decreases to zero due to the rigid wall condition. The reduced order solution \ref{fig3:rom_vel} closely approximates the high-fidelity solution, achieving an absolute error of order $\mathcal{O}\left(10^{-3} \right)$  with fewer POD modes. Similarly, Figure~\ref{fig4:press_comp_bif} shows comparable results for pressure distribution. We observe a significant pressure drop from the inlet to the outlets, where the higher inlet pressure drives the flow, and the decreasing pressure corresponds to an increase in velocity as the fluid exits. These figures demonstrate that the proposed ROM accurately reproduces the detailed velocity and pressure distributions from the Galerkin FE formulation, providing significant computational efficiency and accuracy.  

%%--------------------------------------------------------%%
\subsubsection{Flow field Characteristics}
\label{sec4.2.2:Bif_flow}
\begin{figure}[htbp]
  \centering
    \begin{subfigure}[b]{0.175\textwidth}
    \centering
    \includegraphics[width=\textwidth]{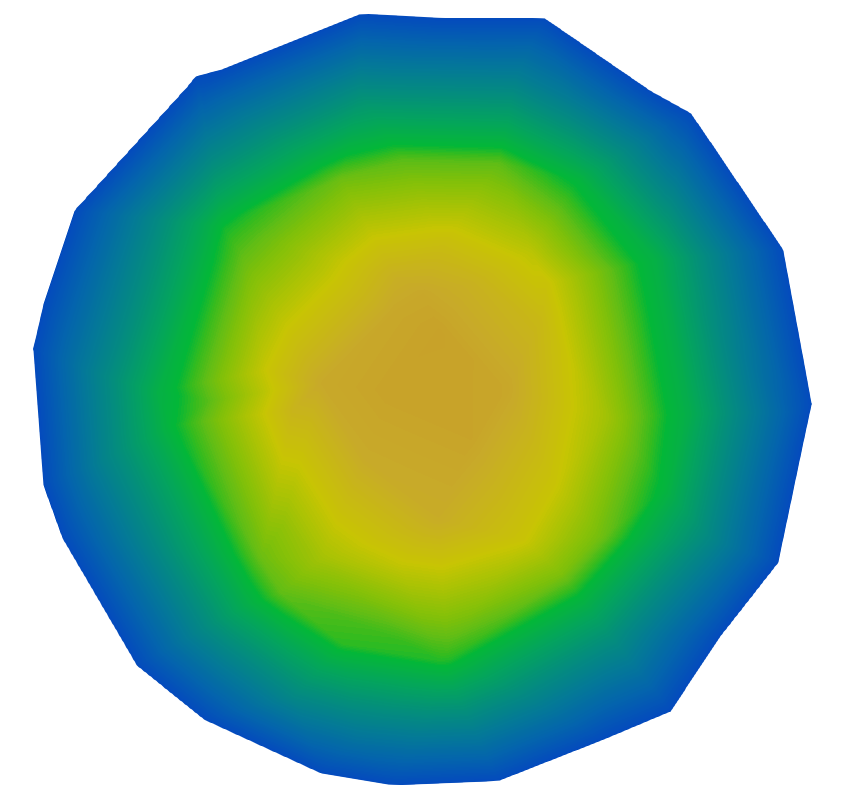}
    \caption{$t=0.05$ s}
    \label{fig5:bif_v_t0.05}
  \end{subfigure}
  \hfill
  \begin{subfigure}[b]{0.175\textwidth}
    \centering
    \includegraphics[width=\textwidth]{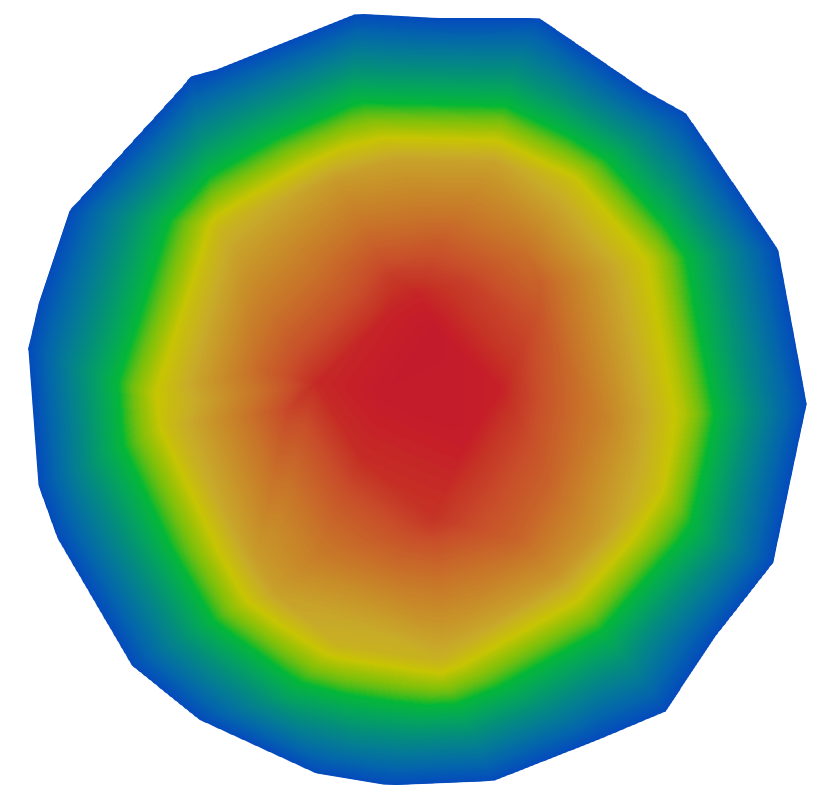}
    \caption{$t=0.25$ s}
    \label{fig5:bif_v_t0.25}
  \end{subfigure}
  \hfill
  \begin{subfigure}[b]{0.175\textwidth}
    \centering
    \includegraphics[width=\textwidth]{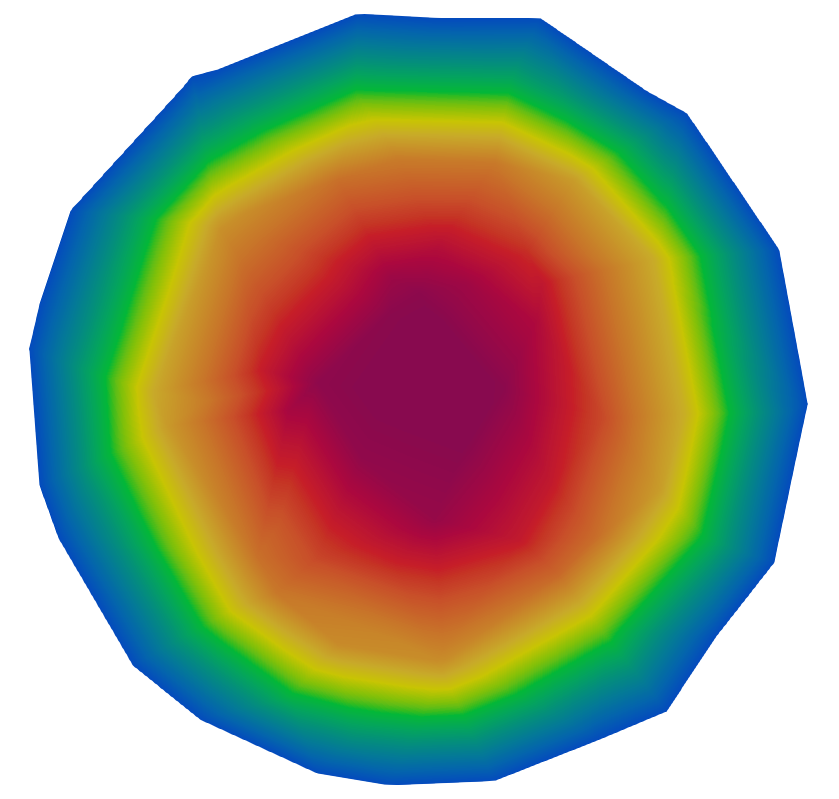}
    \caption{$t=0.5$ s}
    \label{fig5:bif_v_t0.5}
  \end{subfigure}
  \hfill
  \begin{subfigure}[b]{0.175\textwidth}
    \centering
    \includegraphics[width=\textwidth]{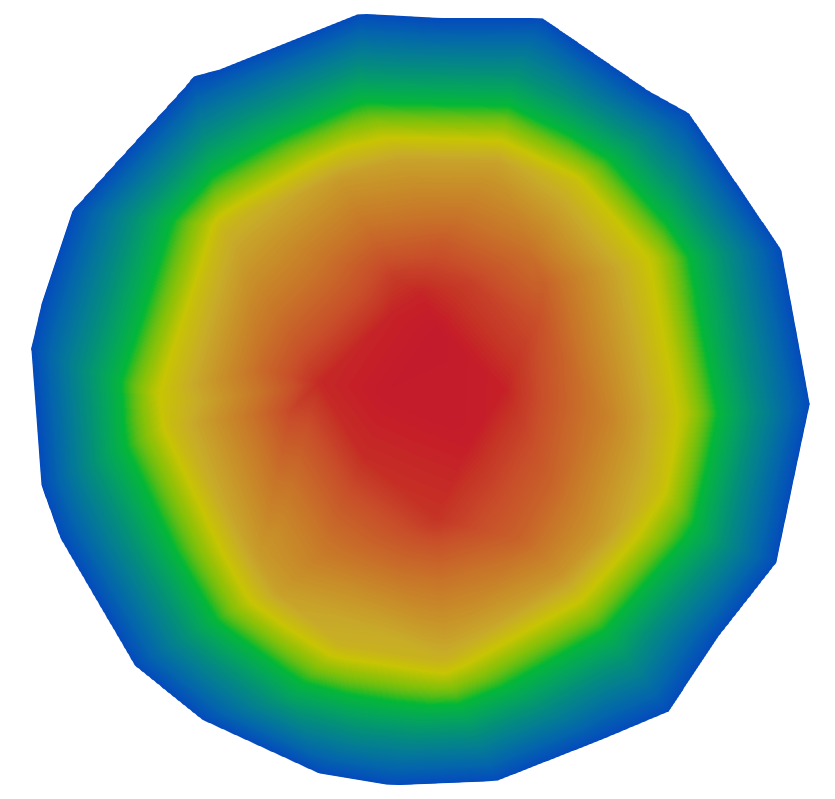}
    \caption{$t=0.75$ s}
    \label{fig5:bif_v_t0.75}
  \end{subfigure}
    \hfill
  \begin{subfigure}[b]{0.175\textwidth}
    \centering
    \includegraphics[width=\textwidth]{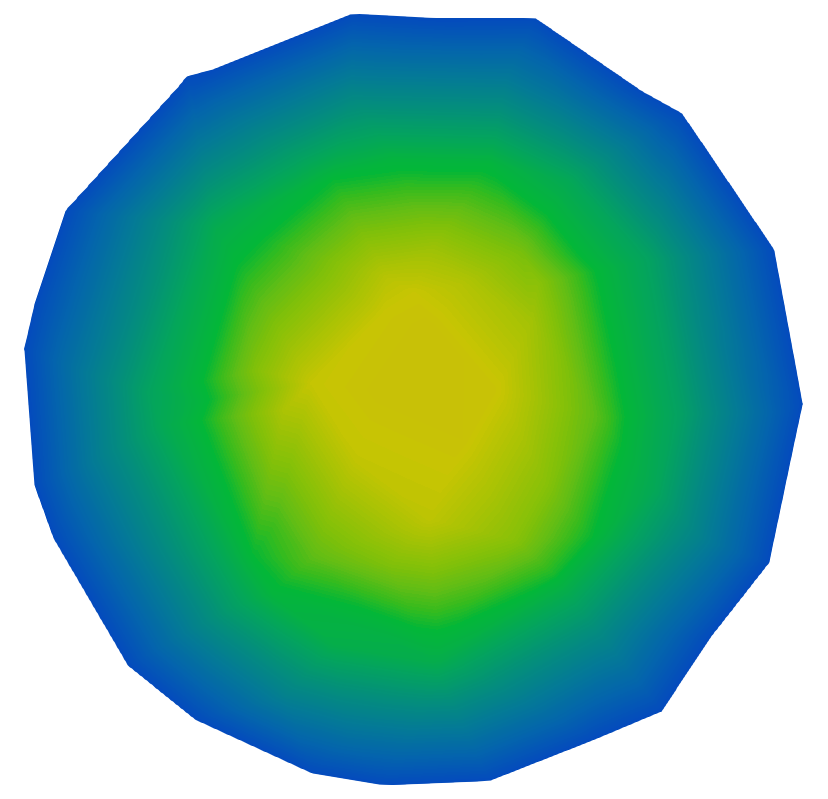}
    \caption{$t=1.0$ s}
    \label{fig5:bif_v_t1.0}
  \end{subfigure}

    %\vskip\baselineskip 

    \begin{subfigure}[b]{0.6\textwidth}
        \centering
        \includegraphics[width=0.55\textwidth]{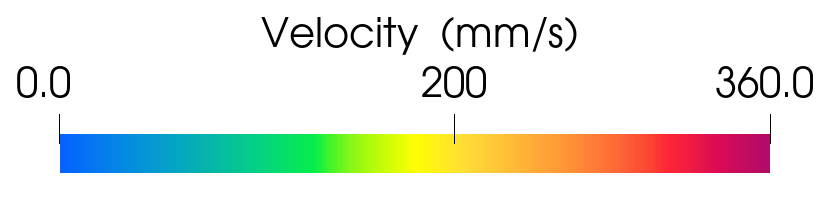}
        \label{fig5:bif_v_t_colorbar}
    \end{subfigure}

  \vskip\baselineskip 
    
      \begin{subfigure}[b]{0.175\textwidth}
    \centering
    \includegraphics[width=\textwidth]{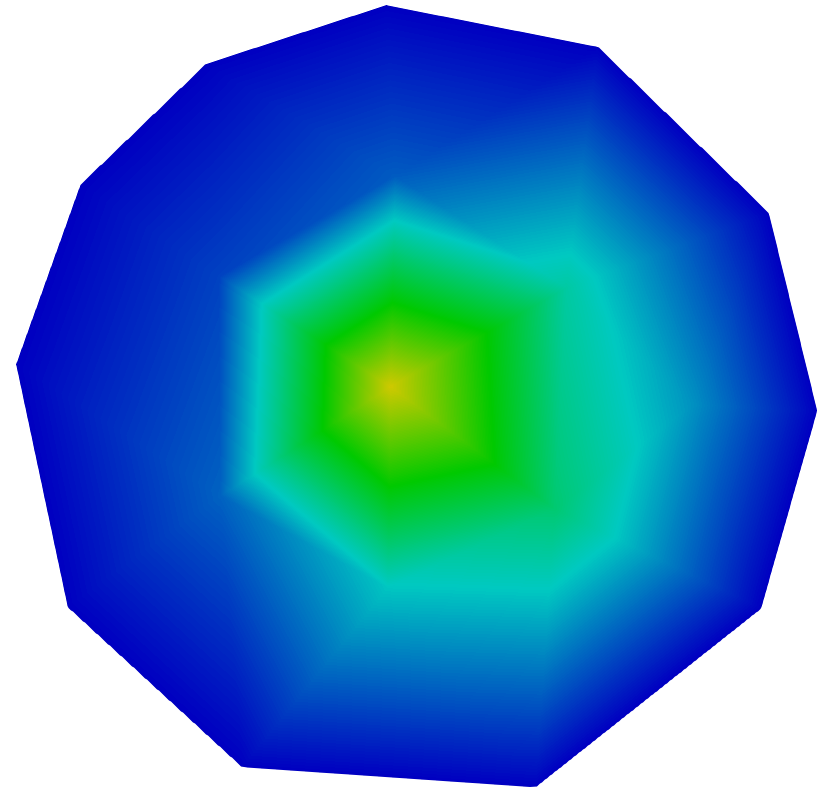}
    \caption{$t=0.05$ s}
    \label{fig5:bif_u_t0.05}
  \end{subfigure}
  \hfill
  \begin{subfigure}[b]{0.175\textwidth}
    \centering
    \includegraphics[width=\textwidth]{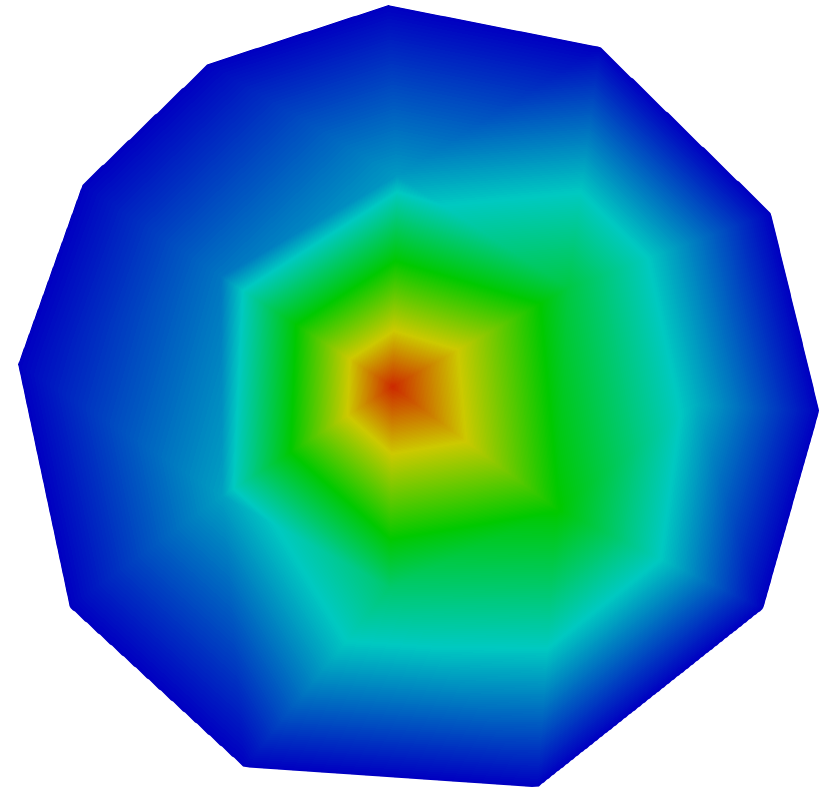}
    \caption{$t=0.25$ s}
    \label{fig5:bif_u_t0.25}
  \end{subfigure}
  \hfill
  \begin{subfigure}[b]{0.175\textwidth}
    \centering
    \includegraphics[width=\textwidth]{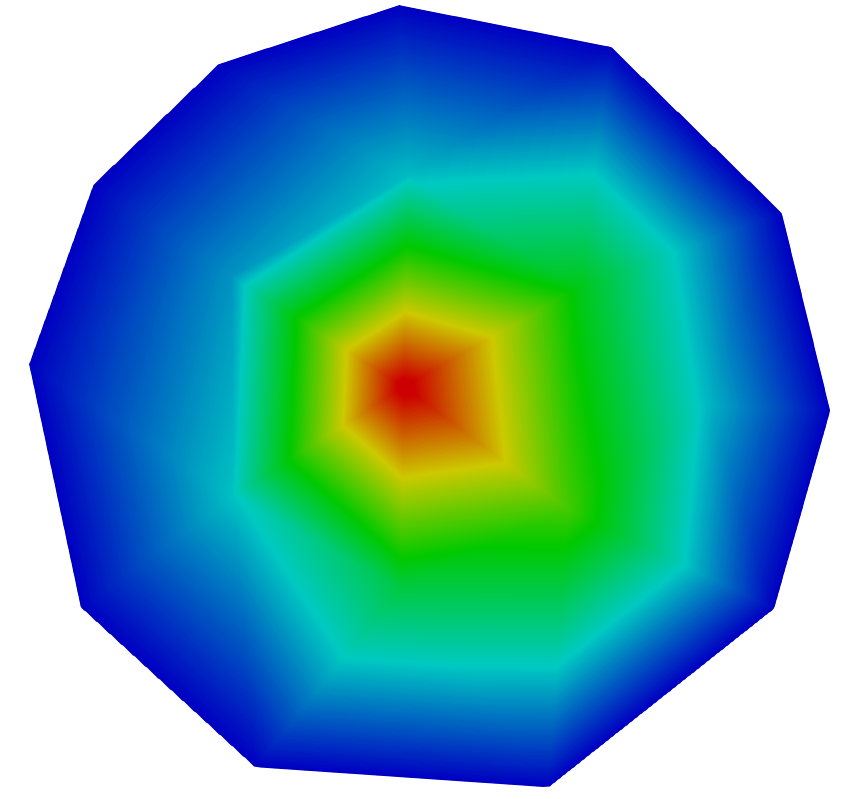}
    \caption{$t=0.5$ s}
    \label{fig5:bif_u_t0.5}
  \end{subfigure}
  \hfill
  \begin{subfigure}[b]{0.175\textwidth}
    \centering
    \includegraphics[width=\textwidth]{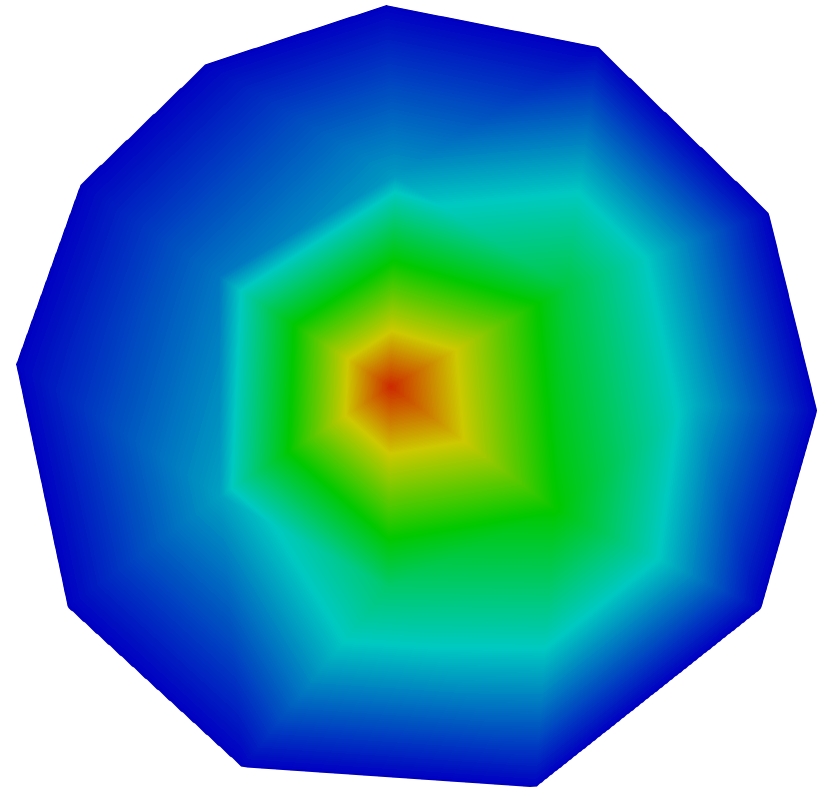}
    \caption{$t=0.75$ s}
    \label{fig5:bif_u_t0.75}
  \end{subfigure}
    \hfill
  \begin{subfigure}[b]{0.175\textwidth}
    \centering
    \includegraphics[width=\textwidth]{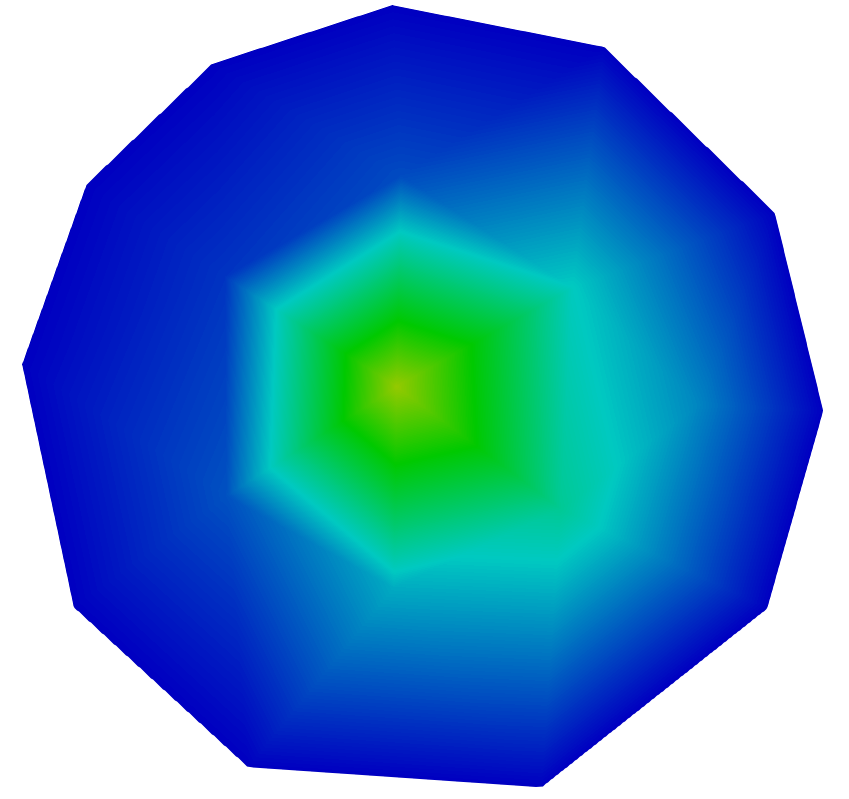}
    \caption{$t=1.0$ s}
    \label{fig5:bif_u_t1.0}
  \end{subfigure}

    %\vskip\baselineskip 

    \begin{subfigure}[b]{0.6\textwidth}
        \centering
        \includegraphics[width=0.55\textwidth]{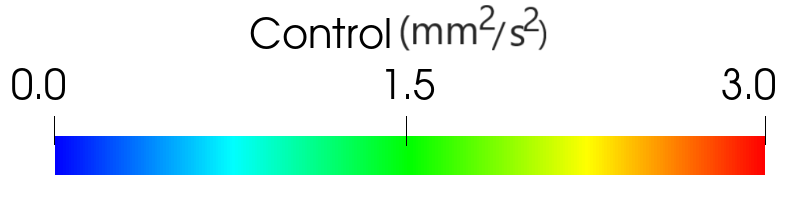}
        \label{fig5:bif_u_t_colorbar}
    \end{subfigure}
  \caption{Temporal variation at $Re = 70$: Velocity distribution at bifurcation site (top row), and Control distribution at $\Gamma_\mathrm{out}$ (bottom row).}
  \label{fig5:timevar_bif}
\end{figure}

Figure~\ref{fig5:timevar_bif} presents the velocity distribution at the bifurcation site, illustrated in the Figure~\ref{fig1:Bif_geom}, specifically on the $xz$-plane of the domain, and the control distribution at the outflow boundary $\Gamma_\mathrm{out}$ at different time instances for $Re=70$. From the figure, we can investigate the flow dynamics where two inlet flows merge into a single outlet, influencing the flow characteristics such as velocity profiles and control distributions. At any fixed time $t$, we observed that the velocity profile exhibits a parabolic shape, with the highest velocity magnitude at the centre of the bifurcation site, which gradually diminishes towards the walls and becomes zero due to the considered assumptions, as shown in the top row of Figure~\ref{fig5:timevar_bif}. As time progresses from $t=0.05$ s to $t=0.5$ s, the velocity significantly increases, indicating the flow is fully developed. The velocity is influenced by the function $f\left(t\right)$, which peaks at $t=0.5$ s; subsequently, the velocity decreases.
\begin{figure}[htbp]
  \centering
  
  \begin{subfigure}[b]{0.18\textwidth}
    \centering
    \includegraphics[width=\textwidth]{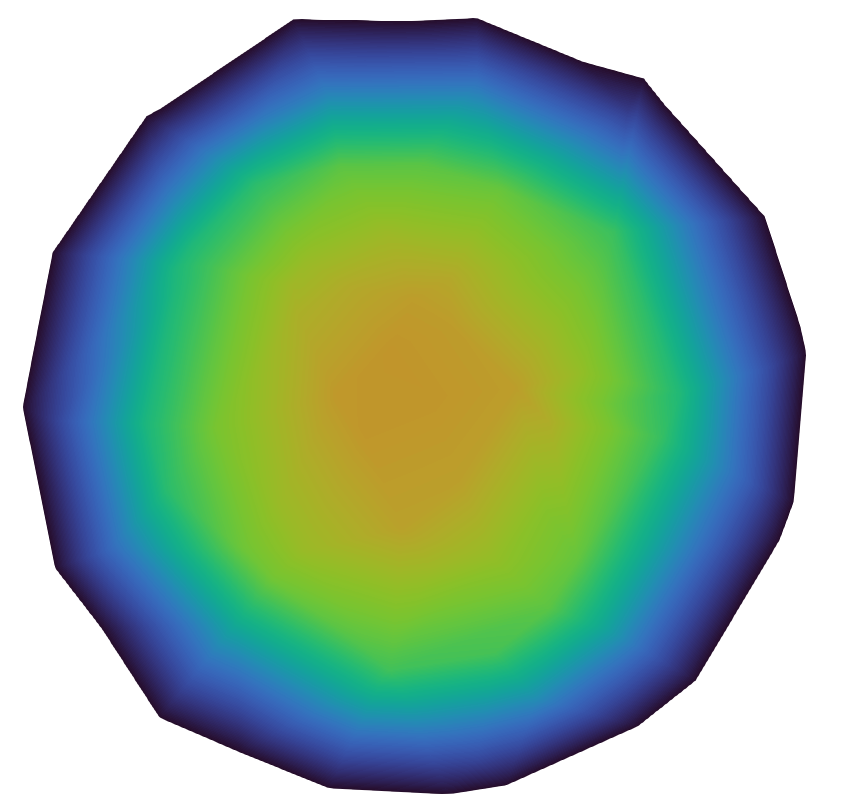}
    \caption{$Re = 50$}
    \label{fig6:bif_v_Re50}
  \end{subfigure}
  \hfill
  \begin{subfigure}[b]{0.18\textwidth}
    \centering
    \includegraphics[width=\textwidth]{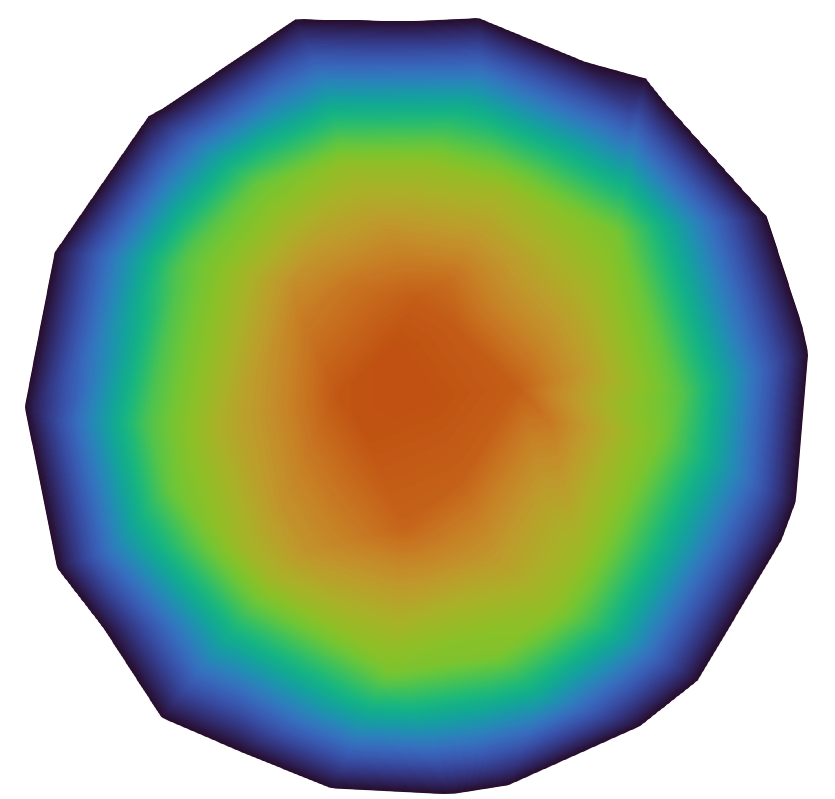}
    \caption{$Re = 60$}
    \label{fig6:bif_v_Re60}
  \end{subfigure}
  \hfill
  \begin{subfigure}[b]{0.18\textwidth}
    \centering
    \includegraphics[width=\textwidth]{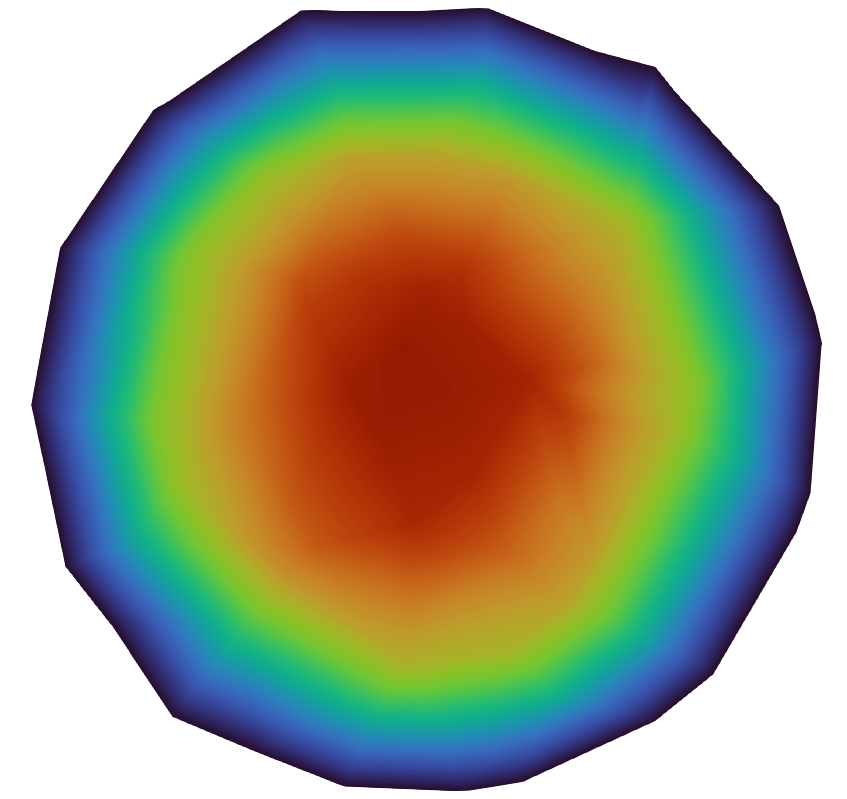}
    \caption{$Re = 70$}
    \label{fig6:bif_v_Re70}
  \end{subfigure}
  \hfill
  \begin{subfigure}[b]{0.18\textwidth}
    \centering
    \includegraphics[width=\textwidth]{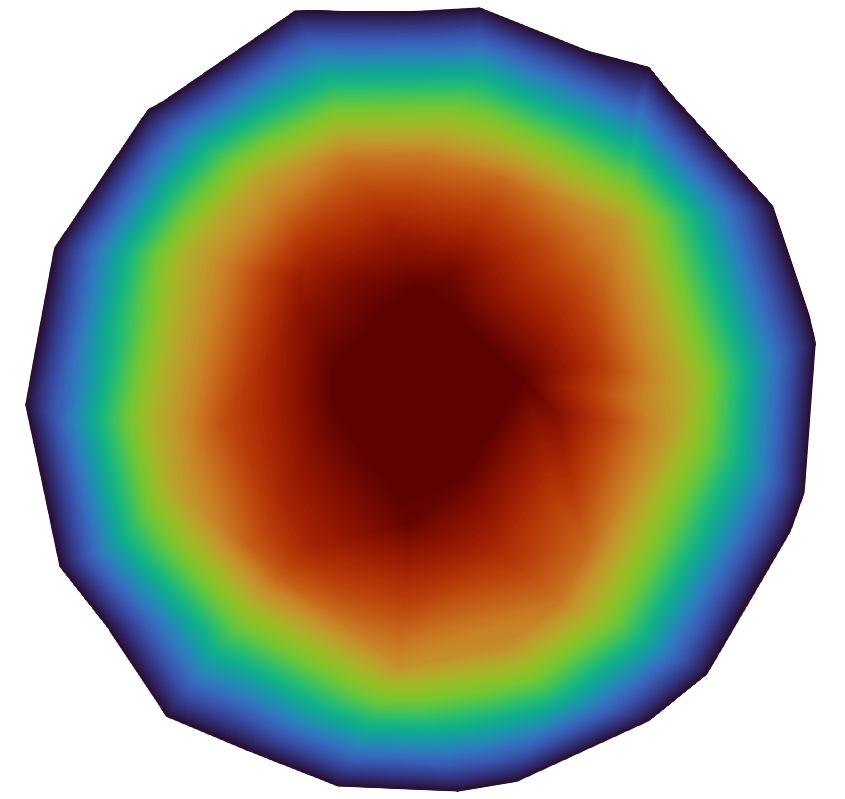}
    \caption{$Re = 80$}
    \label{fig6:bif_v_Re80}
  \end{subfigure}
    %\vskip\baselineskip 
    \begin{subfigure}[b]{0.5\textwidth}
        \centering
        \includegraphics[width=0.8\textwidth]{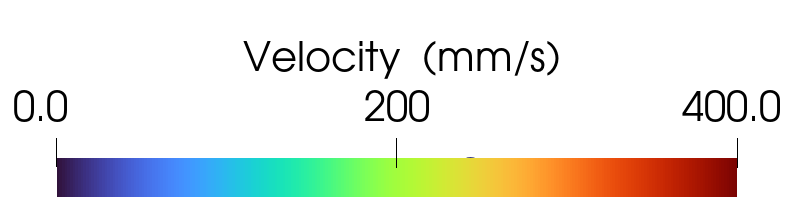}
        \label{fig6:bif_vel_colorbar_t0.5_Re}
    \end{subfigure}
    
  \caption{Parametric dependency of $Re$: Reduced-order velocity distribution at the bifurcation site at $t = 0.5$ s.}
  \label{fig6:Re_vel_bif}
\end{figure}
The control on $\Gamma_\mathrm{out}$, is concentrated at the center and asymmetrically distributed with lower values toward the periphery. At time instant $t=0.05$ s, it is relatively small, effectively regulating the early flow dynamics. As time progresses, the control distribution broadens and highest around $t=0.5$ s, while still peaking at the center. After the peak, the control decreases, allowing the flow to gradually stabilize while still preserving the essential Poiseuille profile. This dynamic adjustment of control ensures smooth regulation of the flow dynamics, as depicted in the bottom row of Figure~\ref{fig5:timevar_bif}.
\begin{figure}[htbp]
  \centering
  
  \begin{subfigure}[b]{0.245\textwidth}
    \centering
    \includegraphics[width=\textwidth]{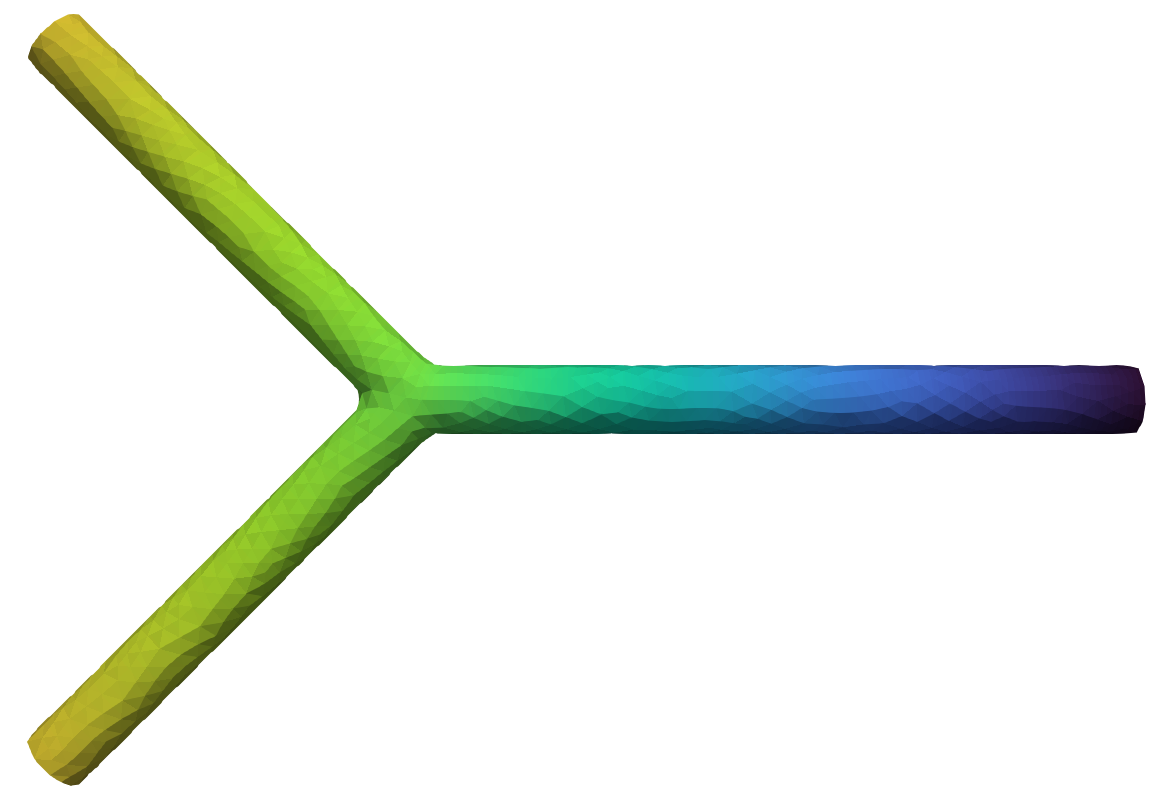}
    \caption{$Re = 50$}
    \label{fig7:bif_Re50_press}
  \end{subfigure}
  \hfill
  \begin{subfigure}[b]{0.245\textwidth}
    \centering
    \includegraphics[width=\textwidth]{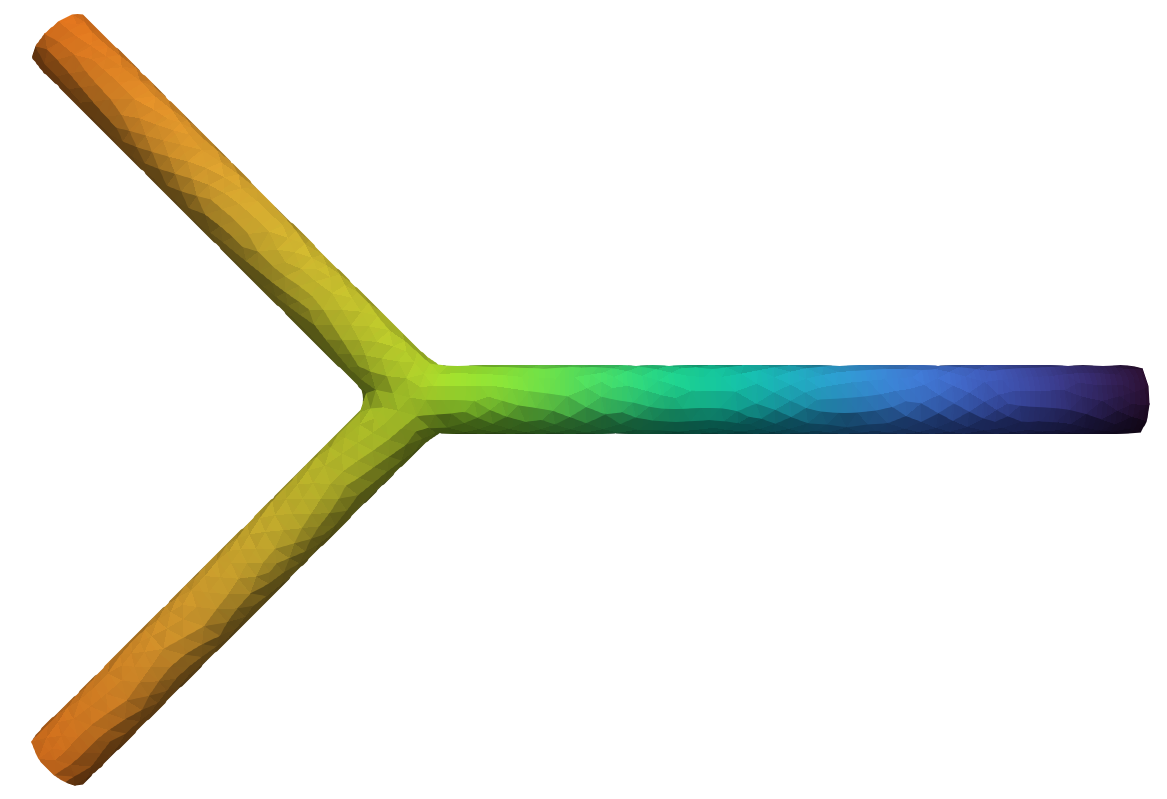}
    \caption{$Re = 60$}
    \label{fig7:bif_Re60_press}
  \end{subfigure}
  \hfill
  \begin{subfigure}[b]{0.245\textwidth}
    \centering
    \includegraphics[width=\textwidth]{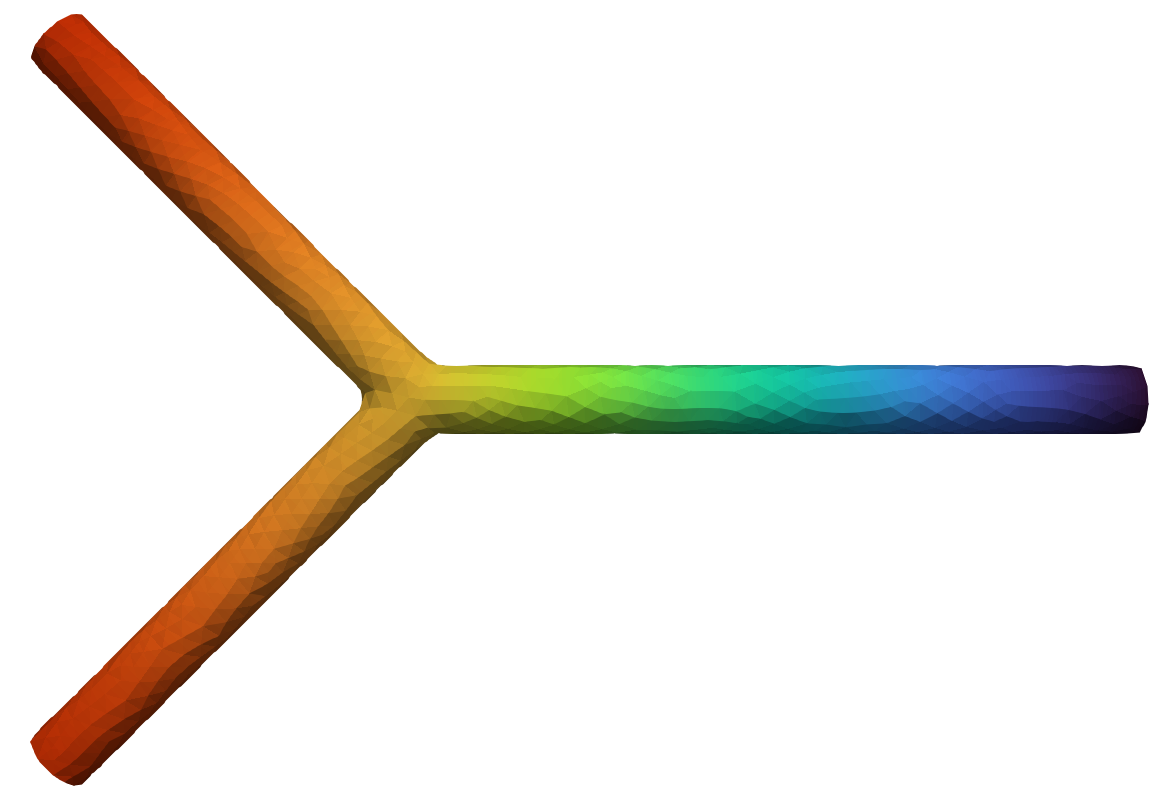}
    \caption{$Re = 70$}
    \label{fig7:bif_Re70_press}
  \end{subfigure}
  \hfill
  \begin{subfigure}[b]{0.245\textwidth}
    \centering
    \includegraphics[width=\textwidth]{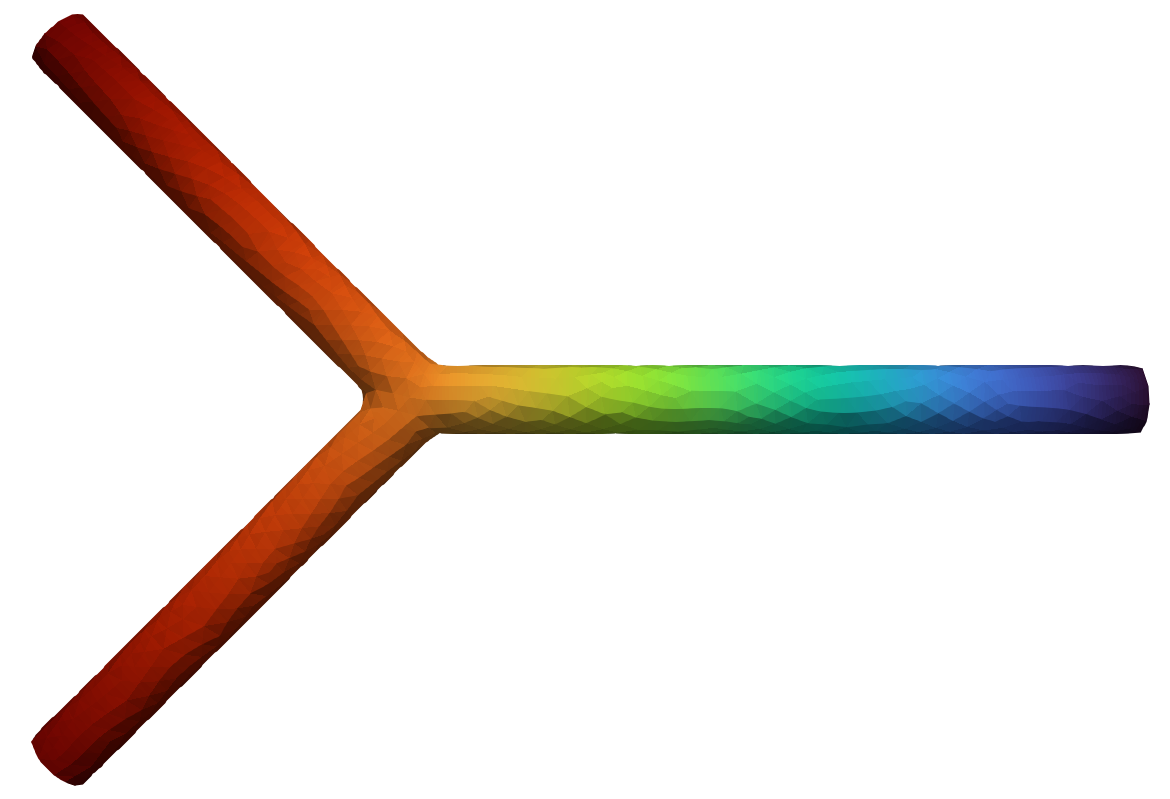}
    \caption{$Re = 80$}
    \label{fig7:bif_Re80_press}
  \end{subfigure}

   % \vskip\baselineskip 

    \begin{subfigure}[b]{0.5\textwidth}
        \centering
        \includegraphics[width=0.9\textwidth]{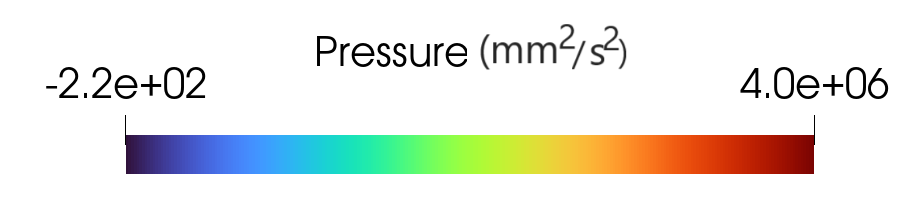}
        \label{fig7:bif_press_colorbar_t0.5_Re}
    \end{subfigure}
  
  \caption{Parametric dependency of $Re$: Reduced order pressure distribution at $t=0.5$ s.}
  \label{fig7:Re_press_bif}
\end{figure}
Moreover, Figures~\ref{fig6:Re_vel_bif} and \ref{fig7:Re_press_bif} present the parametric dependency of velocity and pressure distributions for various $Re$ values at $t=0.5$ s. As $Re$ increases, inertial forces dominate, leading to higher-speed flow patterns, as depicted in Figure~\ref{fig6:Re_vel_bif}. The velocity and pressure distributions are closely related; high-velocity regions correspond to low-pressure zones and vice versa, ensuring the fluid remains incompressible and momentum is balanced, as shown in Figure~\ref{fig7:Re_press_bif}. The control at $\Gamma_\mathrm{out}$ significantly influences the pressure distribution, ensuring a smooth exit of merged inlet flows from the bifurcation site.
% %%
% \begin{figure}[htbp]
%   \centering
  
%   \begin{subfigure}[b]{0.245\textwidth}
%     \centering
%     \includegraphics[width=\textwidth]{Bif_Adj_Press_Re50_t0.5.png}
%     \caption{$Re = 50$}
%     \label{fig71:bif_adj_press_Re50_press}
%   \end{subfigure}
%   \hfill
%   \begin{subfigure}[b]{0.245\textwidth}
%     \centering
%     \includegraphics[width=\textwidth]{Bif_Adj_Press_Re60_t0.5.png}
%     \caption{$Re = 60$}
%     \label{fig71:bif_adj_press_Re60_press}
%   \end{subfigure}
%   \hfill
%   \begin{subfigure}[b]{0.245\textwidth}
%     \centering
%     \includegraphics[width=\textwidth]{Bif_Adj_Press_Re70_t0.5.png}
%     \caption{$Re = 70$}
%     \label{fig71:bif_adj_press_Re70_press}
%   \end{subfigure}
%   \hfill
%   \begin{subfigure}[b]{0.245\textwidth}
%     \centering
%     \includegraphics[width=\textwidth]{Bif_Adj_Press_Re80_t0.5.png}
%     \caption{$Re = 80$}
%     \label{fig71:bif_adj_press_Re80_press}
%   \end{subfigure}

%    % \vskip\baselineskip 

%     \begin{subfigure}[b]{0.5\textwidth}
%         \centering
%         \includegraphics[width=0.8\textwidth]{Bif_Adj_Press_Re_t0.5_colorbar.png}
%         \label{fig71:bif_adj_press_colorbar_t0.5_Re}
%     \end{subfigure}
  
%   \caption{{Parametric dependency of $Re$: Reduced order Adjoint pressure (Lagrangian multiplier) distribution at $t=0.5$ s.}}
%   \label{fig71:Re_press_bif}
% \end{figure}
% %% 
%%%%%%%%%%%%%%%%%%%%%%%%%%%%%%%%%%%%%%%%%%%%%%%%%%%%%%%%%%%%%%%%%%%%%%%%%%%%%%%%%%%%%%%%%%
\subsection{Test case 2: Patient-specific CABG}
\label{sec4.3:cabg_case}
This test case, a patient-specific CABG model, represents a realistic scenario of the idealized bifurcation model presented in Test Case \ref{sec4.1:testcase1}. 
%The idealized bifurcation model simplifies the geometry to understand the fundamental flow dynamics, while the CABG model incorporates the complex morphology of the coronary artery, providing a more accurate and realistic simulation of CV flow dynamics. 
For this case, the final time is $T = 0.8$ s as in \cite{ballarin2016fast,ballarin2017numerical} with a time-step $\Delta t = 0.01$ s, storing the solutions at every $4^{th}$ time-step, resulting in 20 snapshots in time. We considered parametric-dependent inlet flow profiles computed using Eq.~\eqref{eq:4.1} with $\left(Re_1, Re_2\right)$ values independently within the range $[50, 80]$ for each inlet boundary. 

%%--------------------------------------------------------%%
\subsubsection{Assessment of ROM Performance}
\label{sec4.3.1:perfor}
\begin{figure}[H]
  \centering
  
  % First row of subfigures: POD Singular Values
  \begin{subfigure}[b]{0.475\textwidth}
    \centering
    \includegraphics[width=\textwidth]{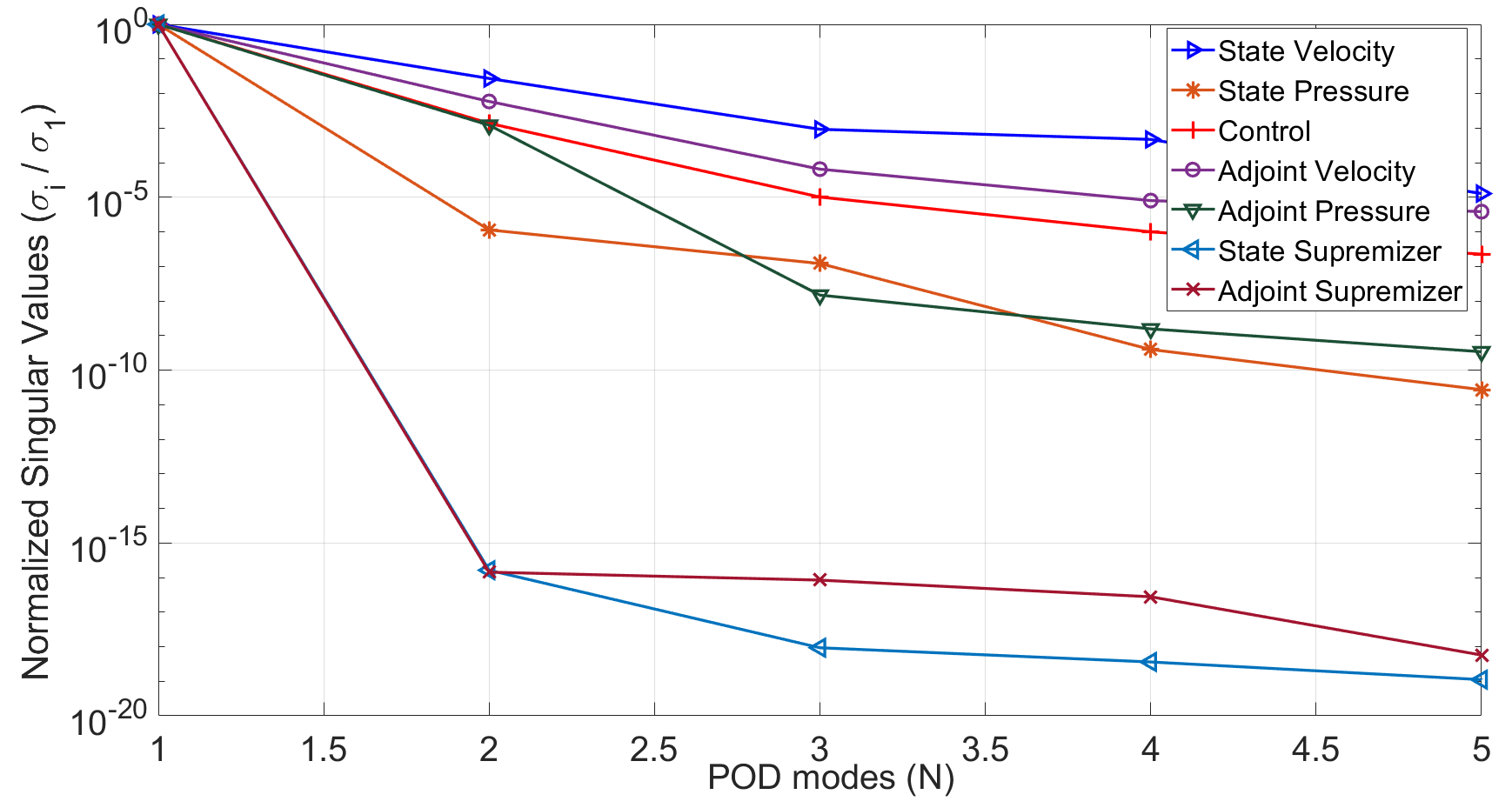}
    \caption{Singular values: \textit{temporal compression}}
    \label{fig8:cabg_pod_eig0}
  \end{subfigure}
  \hfill
  \begin{subfigure}[b]{0.475\textwidth}
    \centering
    \includegraphics[width=\textwidth]{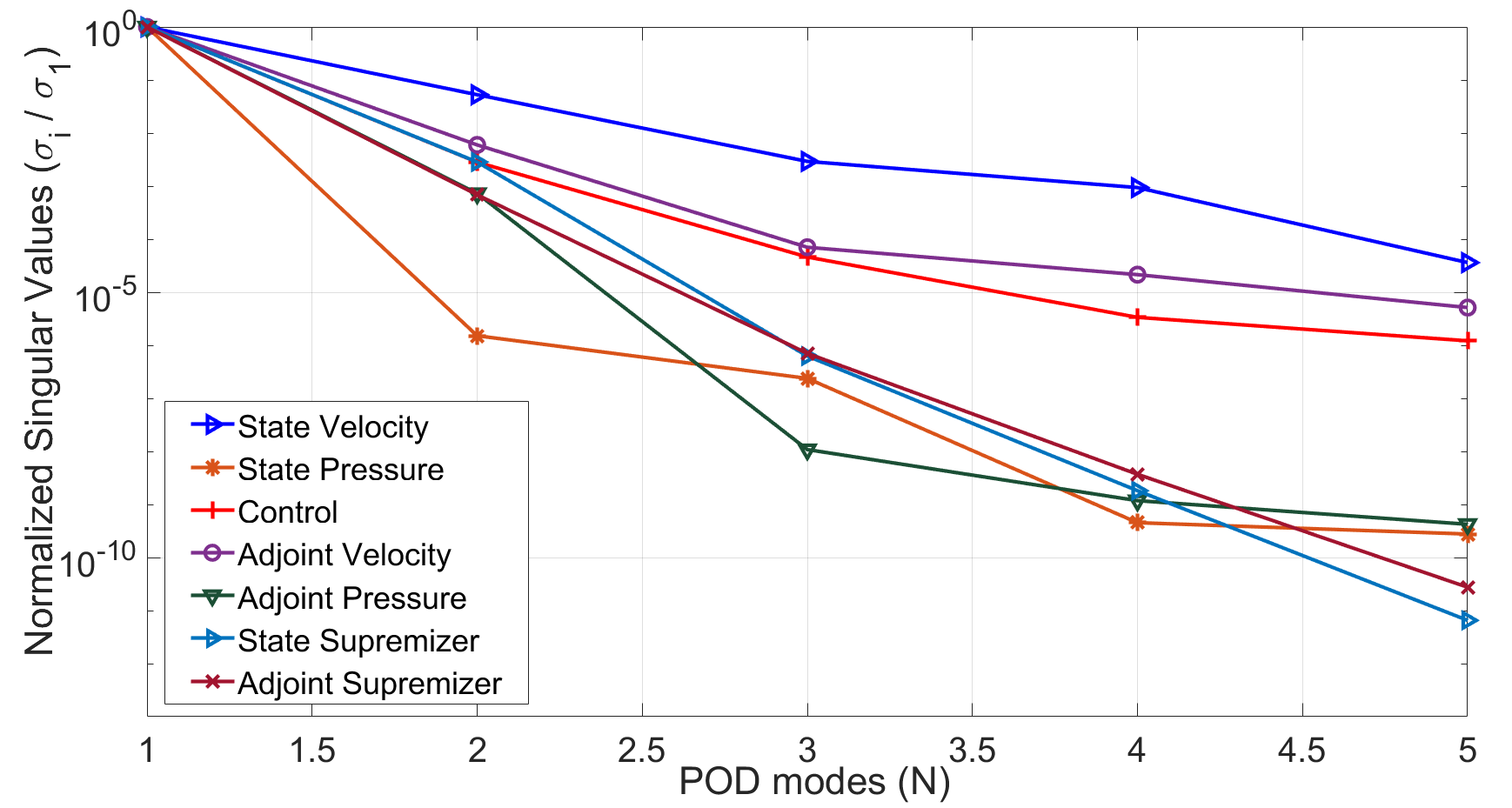}
    \caption{Singular values: \textit{parametric-space compression}}
    \label{fig8:cabg_pod_eig1}
  \end{subfigure}
  
  \vskip\baselineskip % Vertical space between the rows

  % Second row of subfigures: POD Retaining Energy
  \begin{subfigure}[b]{0.46\textwidth}
    \centering
    \includegraphics[width=\textwidth]{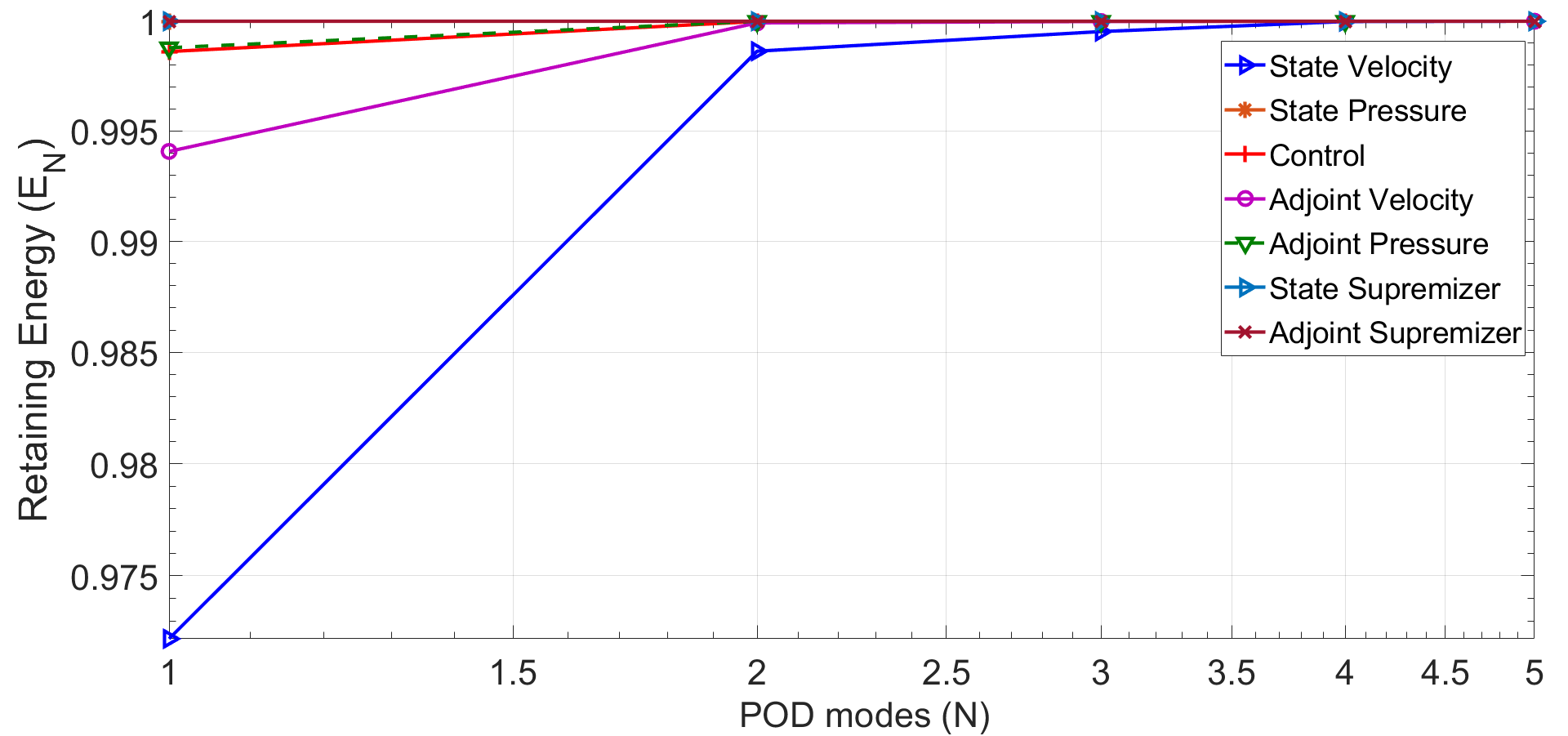}
    \caption{Retaining energy: \textit{temporal compression}}
    \label{fig8:cabg_pod_eng0}
  \end{subfigure}
  \hfill
  \begin{subfigure}[b]{0.48\textwidth}
    \centering
    \includegraphics[width=\textwidth]{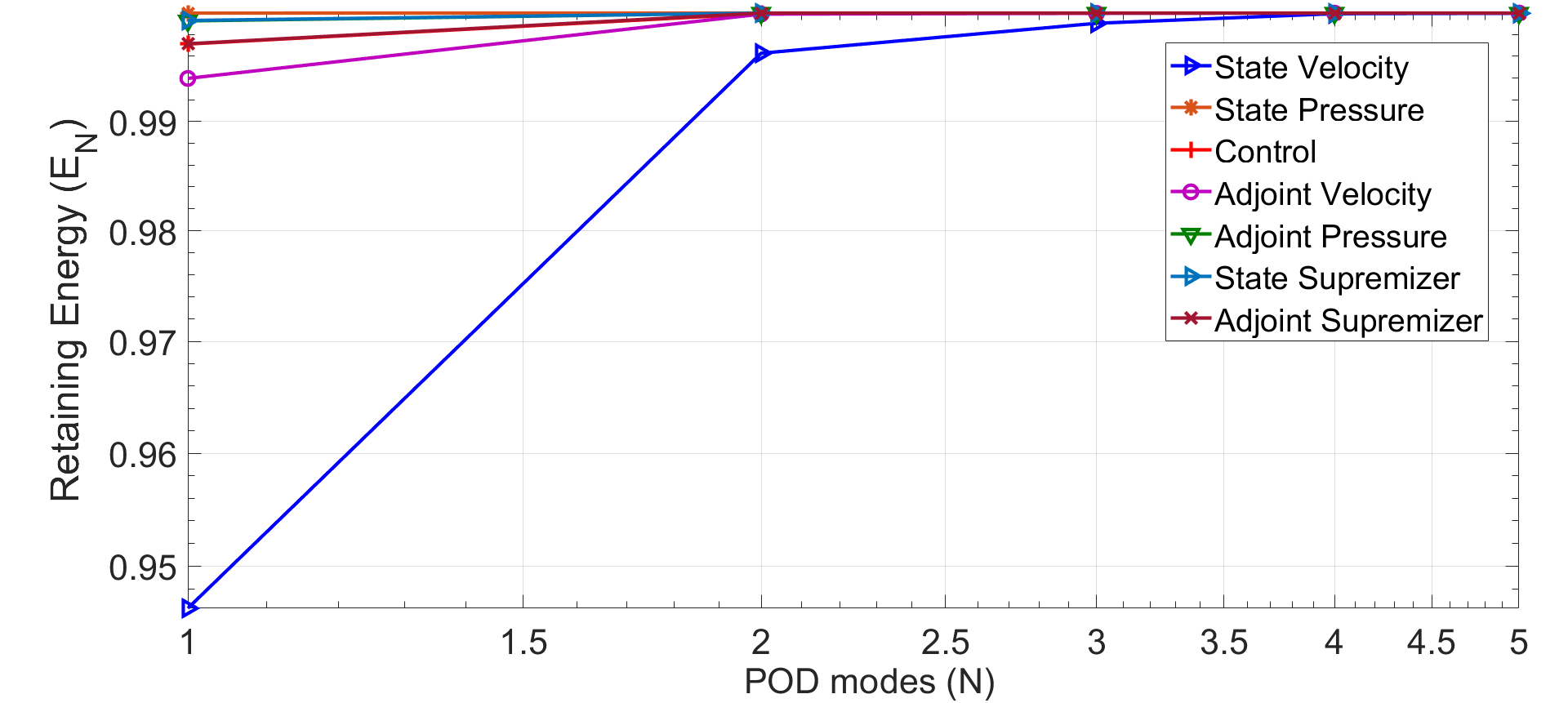}
    \caption{Retaining energy: \textit{parametric-space compression}}
    \label{fig8:cabg_pod_eng1}
  \end{subfigure}

  \caption{{Normalized POD singular values (top row) and retaining energy (bottom row) for velocity, pressure, and supremizers  {for $\left(Re_1, Re_2\right) = \left(50, 50\right)$}.}}
  \label{fig8:eigen_cabg}
\end{figure}

In this case study, we consider a patient-specific CABG model characterized by its complex morphology and a large-scale mesh comprising $\mathcal{N}_h = 433,288$ dofs, posing significant computational challenges. Using the \textit{nested-POD} approach, extracting  $N_t^{\mathrm{POD}} = 10$ POD-modes significantly reduces the dimensionality from $2500 \times \mathcal{N}_h$ to $250 \times \mathcal{N}_h$  with $N_\mathrm{train} = 25$ snapshots, for capturing the dominant dynamics of the system. Each $\bm{\mu}$ snapshot in the offline phase required about $7-8$ hours, adding up to approximately $8$ days for all snapshots. In contrast, the online phase significantly enhanced efficiency, reducing the CPU time to around  $1$ hour per snapshot, resulting in substantial time savings. Figure \ref{fig8:eigen_cabg}  presents POD singular values and retained energy for all variables during the temporal and parametric-space compression.
We observed that state velocity, adjoint velocity, and control exhibit significant decay in normalized singular values, while pressure and supremizers show a more rapid decay. This suggests the system dynamics can be captured with a reduced number of POD modes $N \leq 5$. Additionally, we noticed that state variables decay more slowly compared to adjoint variables and retain the $95\%$ energy of the system, as shown in Figures~\ref{fig8:cabg_pod_eng0} and \ref{fig8:cabg_pod_eng1}. Subsequently, the essential spatio-temporal features of the cardiovascular system are accurately captured using a minimal number of POD modes, highlighting the effectiveness of the proposed approach.

\begin{figure}[htbp]
    \centering
    \begin{subfigure}[b]{0.325\textwidth}
        \centering
        \includegraphics[width=\textwidth]{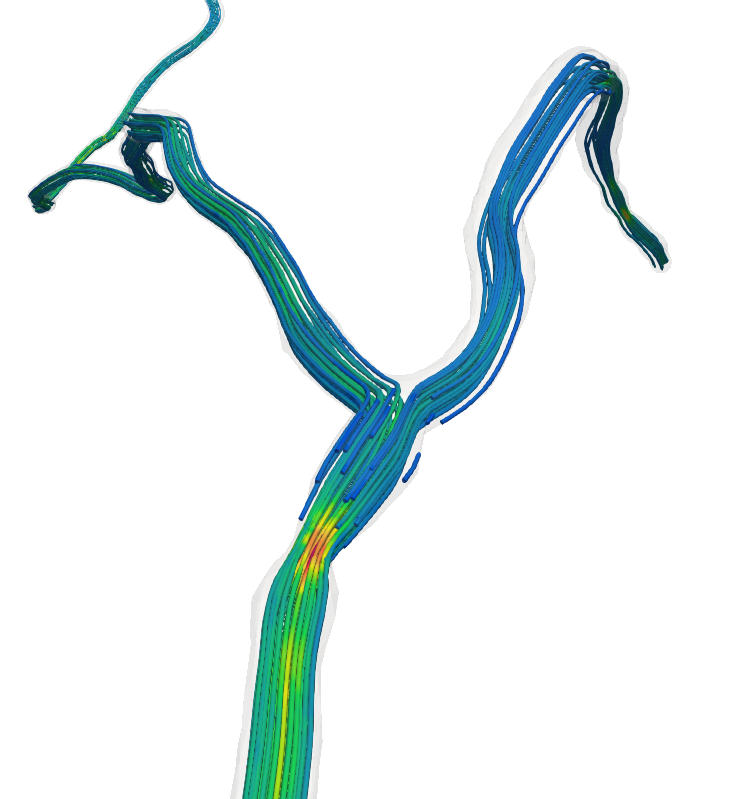}
        \caption{\textit{High-fidelity} solution}
        \label{fig9:cabg_fom_vel}
    \end{subfigure}%
    \hfill
    \begin{subfigure}[b]{0.325\textwidth}
        \centering
        \includegraphics[width=\textwidth]{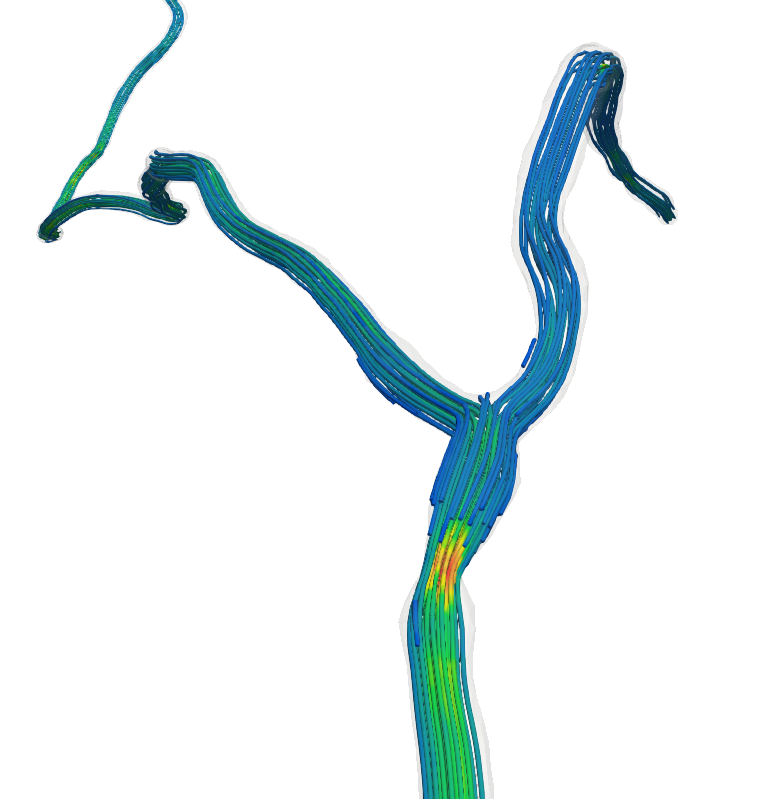}
        \caption{\textit{Reduced-order} solution}
        \label{fig9:cabg_rom_vel}
    \end{subfigure}%
    \hfill
    \begin{subfigure}[b]{0.325\textwidth}
        \centering
        \includegraphics[width=\textwidth]{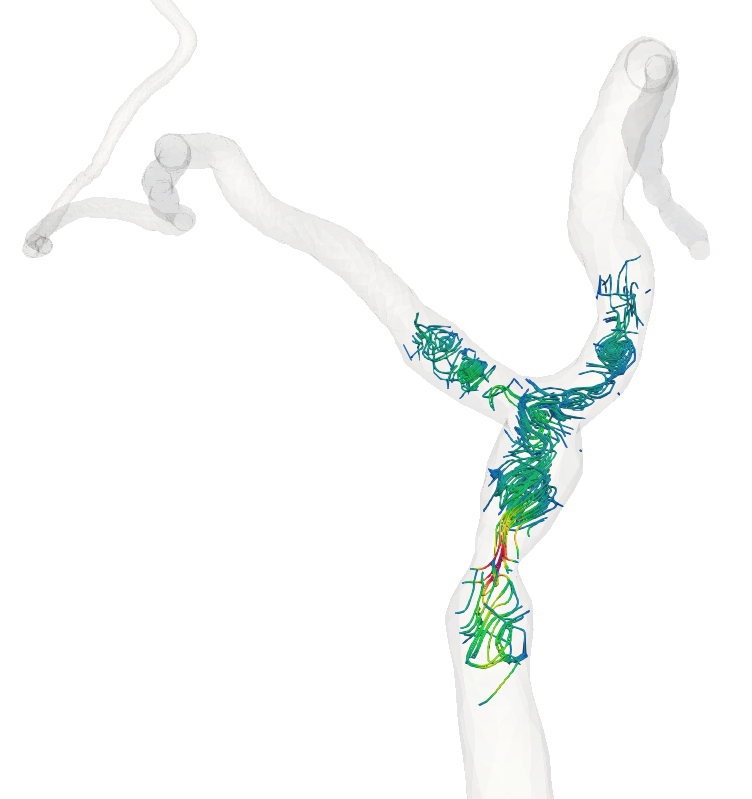}
        \caption{Absolute error}
        \label{fig9:cabg_error_vel}
    \end{subfigure}
    
    \vskip\baselineskip % Add space between rows
    
    \begin{subfigure}[b]{0.65\textwidth}
        \centering
        \includegraphics[width=0.51\textwidth]{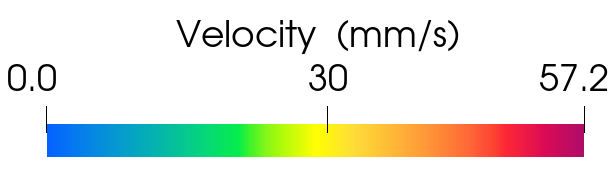}
        \label{fig9:cabg_velocity_colorbar}
    \end{subfigure}
    \hfill
    \begin{subfigure}[b]{0.3\textwidth}
        \centering
        \includegraphics[width=0.95\textwidth]{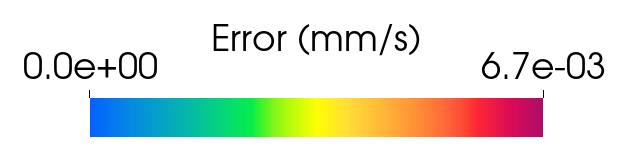}
        \label{fig9:cabg_v_e}
    \end{subfigure}
    
    \caption{Comparison of the streamlines at $\left(Re_1, Re_2\right) = \left(80, 80\right)$ and $t = 0.4$ s. (a) High-fidelity solution, (b) Reduced-order solution, and (c) Absolute error.}
    \label{fig9:vel_comp_cabg}
\end{figure}
\begin{figure}[htbp]
    \centering
    \begin{subfigure}[b]{0.325\textwidth}
        \centering
        \includegraphics[width=0.8\textwidth]{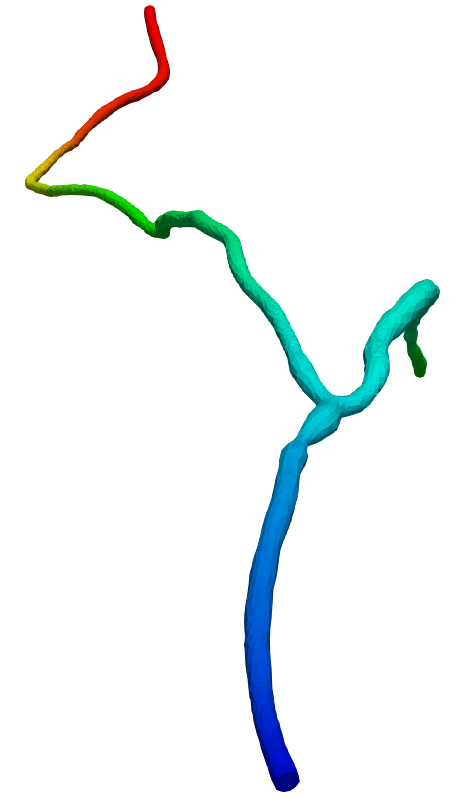}%
        \caption{\textit{High-fidelity} solution}
        \label{fig10:cabg_fom_press}
    \end{subfigure}%
    \hfill
    \begin{subfigure}[b]{0.325\textwidth}
        \centering
        \includegraphics[width=0.8\textwidth]{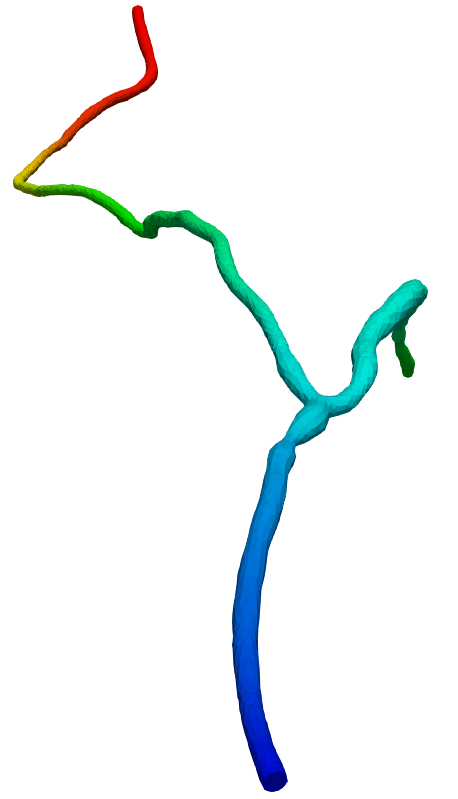}
        \caption{\textit{Reduced-order} solution}
        \label{fig10:cabg_rom_press}
    \end{subfigure}%
    \hfill
    \begin{subfigure}[b]{0.325\textwidth}
        \centering
        \includegraphics[width=0.8\textwidth]{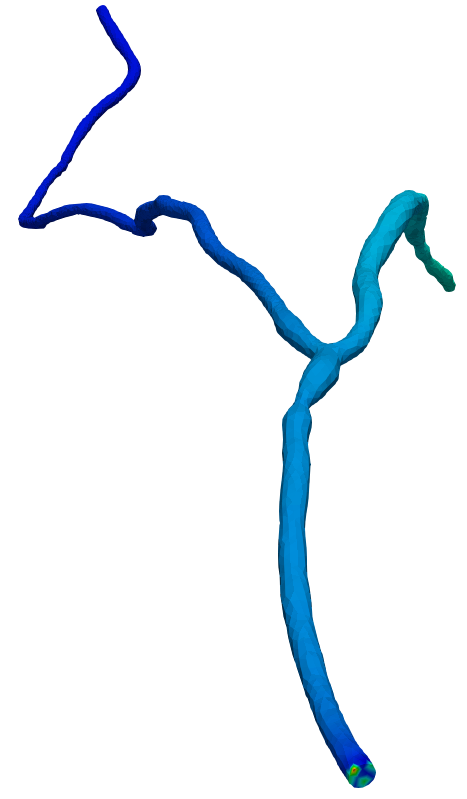}
        \caption{Relative error}
        \label{fig10:cabg_error_press}
    \end{subfigure}
    %\vskip\baselineskip % Add space between rows
    \begin{subfigure}[b]{0.65\textwidth}
        \centering
        \includegraphics[width=0.5\textwidth]{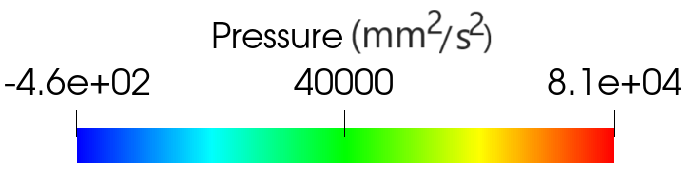}
        \label{fig10:cabg_press_colorbar}
    \end{subfigure}
    \hfill
    \begin{subfigure}[b]{0.3\textwidth}
        \centering
        \includegraphics[width=0.9\textwidth]{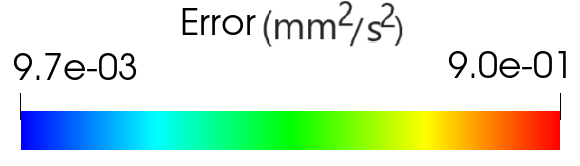}
        \label{fig10:cabg_error_press_colorbar}
    \end{subfigure}
    
    \caption{Comparison of pressure distributions for the unseen Reynolds numbers $\left(Re_1, Re_2\right) = \left(50, 70\right)$ and $t = 0.4$ s. (a) High-fidelity solution, (b) Reduced-order solution, and (c) Relative error.}
    \label{fig10:press_comp_cabg}
\end{figure}
Figure~\ref{fig9:vel_comp_cabg} illustrates the comparison between the high-fidelity and reduced-order solutions for the streamlines, along with the corresponding absolute errors between the two methods at $\left(Re_1, Re_2 \right) = \left(80, 80 \right)$ and $t=0.4$ s. The streamlines in Figure~\ref{fig9:cabg_fom_vel} show the flow behaviour in a patient-specific CABG model, where fluid enters through the RIMA and LAD inlets and exits through the outlet, with high velocities observed in the outlet and region of stenosis. Figure~\ref{fig10:press_comp_cabg} compares the high-fidelity and reduced-order pressure distributions for unseen Reynolds numbers,  specifically $\left(Re_1, Re_2 \right) = \left(50, 70 \right)$ and $t=0.4$ s, along with relative error between the solutions.  The high-fidelity solutions capture the intricate flow dynamics and detail pressure variations, while the reduced-order solution~\ref{fig10:cabg_rom_press}, which closely approximates the overall pressure distribution while requiring significantly less computational effort.
The relative error between the ROM and high-fidelity solutions remains low across most domains, with slightly elevated values near bifurcations and outlet boundaries. The proposed ROM performs effectively even for unseen Reynolds numbers, showcasing its robustness while efficiently approximating high-fidelity solutions, making it a valuable tool for real-time, patient-specific haemodynamics simulations and clinical decision support.

\begin{figure}[htbp]
    \centering
    \begin{subfigure}[b]{0.325\textwidth}
        \centering
        \includegraphics[width=\textwidth]{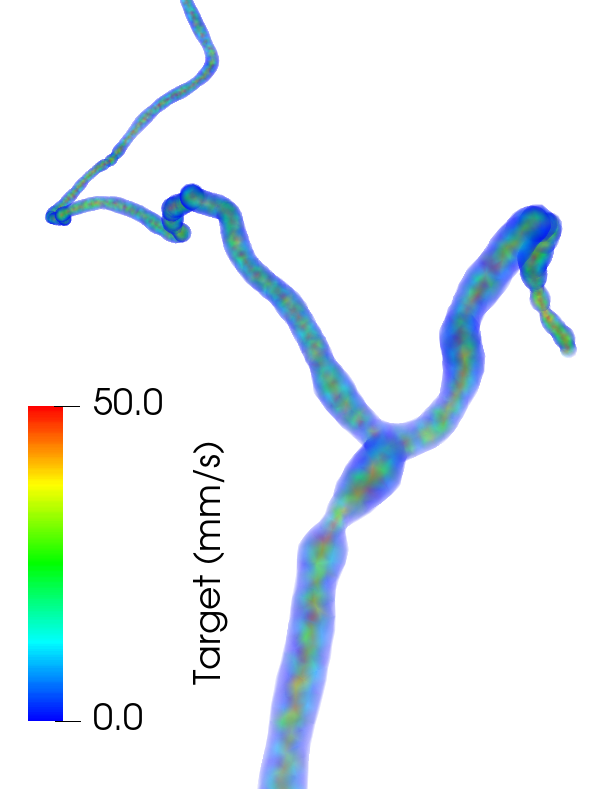}
        \caption{Target Profile}
        \label{fig91:cabg_desired}
    \end{subfigure}%
    \hfill
    \begin{subfigure}[b]{0.325\textwidth}
        \centering
        \includegraphics[width=\textwidth]{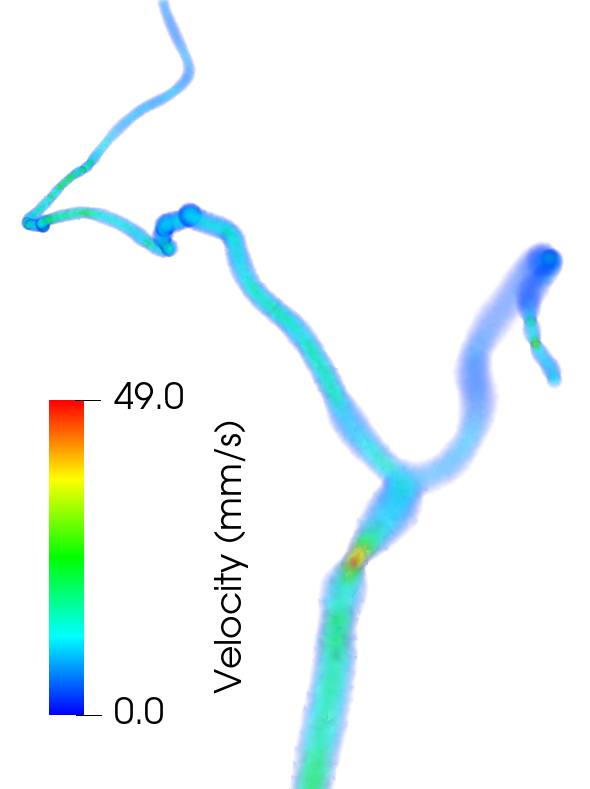}
        \caption{State velocity}
        \label{fig91:cabg_FOM}
    \end{subfigure}%
    \hfill
    \begin{subfigure}[b]{0.325\textwidth}
        \centering
        \includegraphics[width=\textwidth]{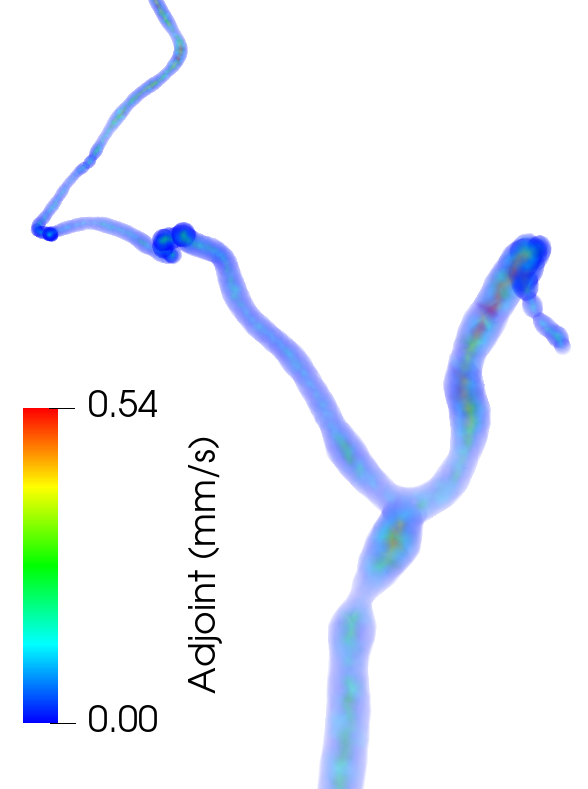}
        \caption{\textit{Adjoint velocity}}
        \label{fig91:cabg_Adj}
    \end{subfigure}
    
    \caption{Comparison of flow profiles for the unseen Reynolds numbers $\left(Re_1, Re_2\right) = \left(70, 70\right)$ and $t = 0.4$ s. (a) Target profile, (b) State velocity, and (c) Adjoint velocity (Lagrangian multiplier) profiles, respectively.}
    \label{fig91:cabg_flow}
\end{figure}
Figure~\ref{fig91:cabg_flow} illustrates a comparative visualization of the target, state velocity, and adjoint velocity profiles for the unseen Reynolds number $\left(Re_1, Re_2\right) = \left(70, 70\right)$ and $t = 0.4$ s. Figure~\ref{fig91:cabg_desired} shows the target velocity profile, representing the desired flow distribution for the optimization problem. Figure~\ref{fig91:cabg_FOM} illustrates the state velocity field obtained from the forward Navier–Stokes simulation under the influence of the optimized outflow boundary control. The close visual agreement between the target and high-fidelity fields indicates the effectiveness of the control strategy in reproducing physiologically realistic flow conditions. Figure~\ref{fig91:cabg_Adj} shows the corresponding adjoint velocity field represents the Lagrange multiplier associated with the optimization problem and encodes sensitivity information for the control objective. This plays a crucial role in driving the optimization by quantifying the influence of flow perturbations on the mismatch between the target and simulated solutions. Together, these profiles highlight the control framework's ability to accurately track desired flow dynamics while enabling gradient-based optimization via the adjoint formulation.

%%======================================================%%
\subsubsection{Flow field Characteristics}
\label{sec4.3.2:cabg_flow}
\begin{figure}[htbp]
  \centering

    % control distribution
  \begin{subfigure}[b]{0.15\textwidth}
    \centering
    \includegraphics[width=\textwidth]{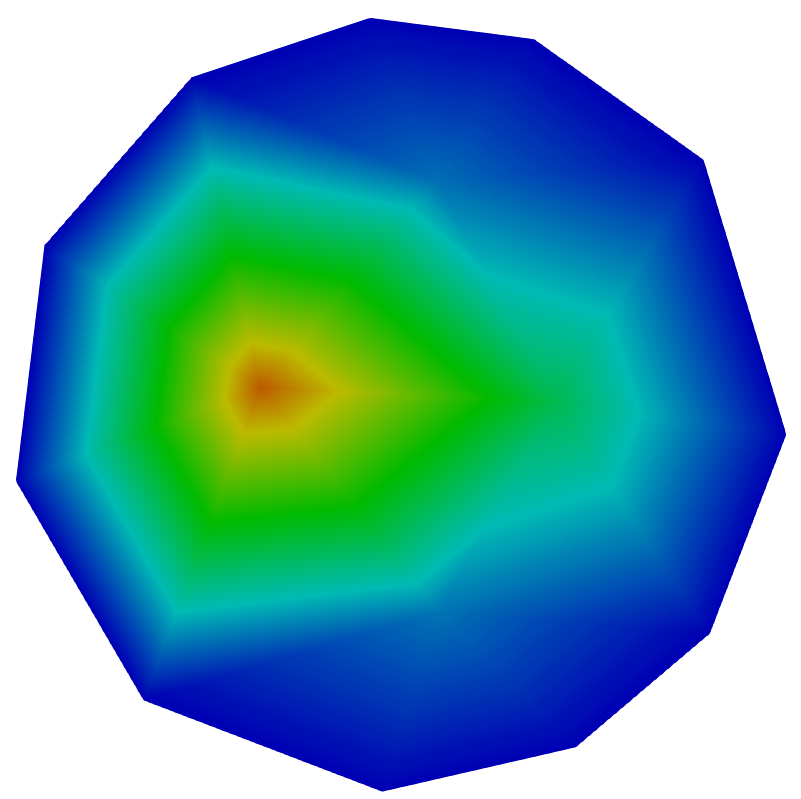}
    \caption{$t=0.04$ s}
    \label{fig11:cabg_u_t0.04}
  \end{subfigure}
  \hfill
  \begin{subfigure}[b]{0.15\textwidth}
    \centering
    \includegraphics[width=\textwidth]{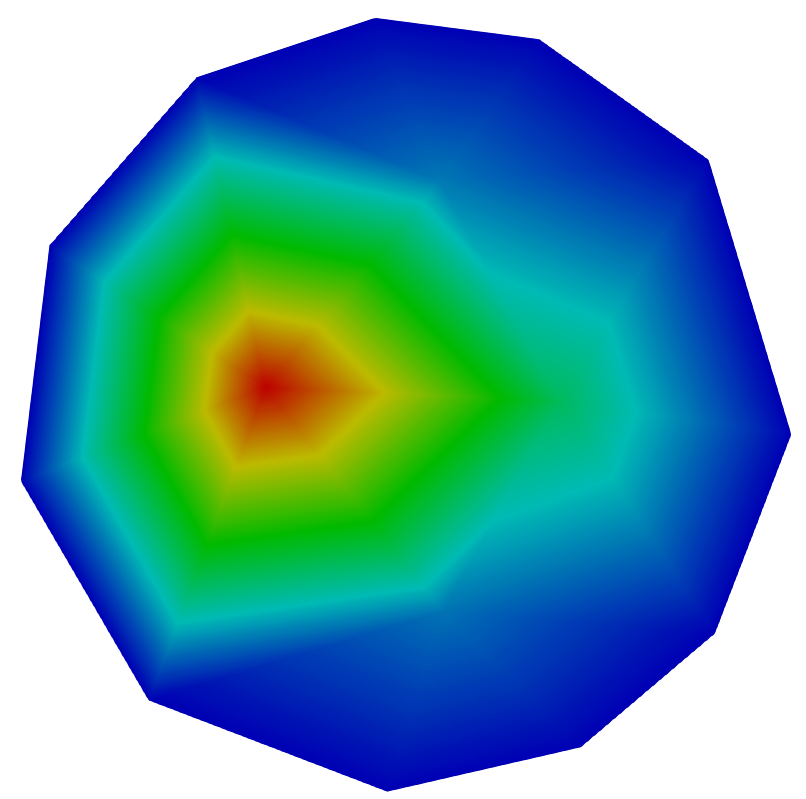}
    \caption{$t=0.2$ s}
    \label{fig11:cabg_u_t0.2}
  \end{subfigure}
  \hfill
  \begin{subfigure}[b]{0.15\textwidth}
    \centering
    \includegraphics[width=\textwidth]{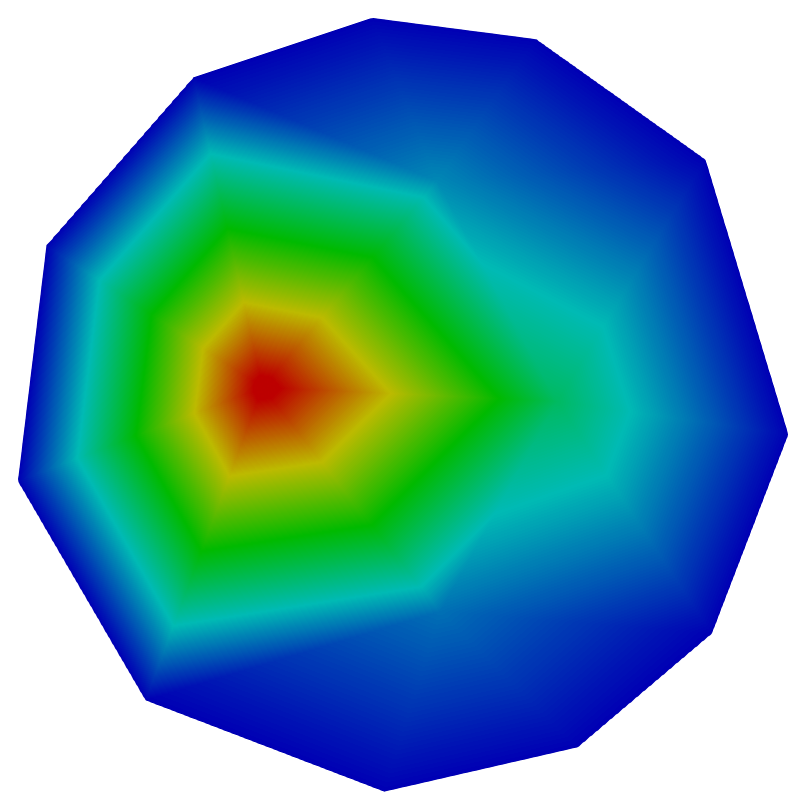}    
    \caption{$t=0.4$ s}
    \label{fig11:cabg_u_t0.4}
  \end{subfigure}
  \hfill
  \begin{subfigure}[b]{0.15\textwidth}
    \centering
    \includegraphics[width=\textwidth]{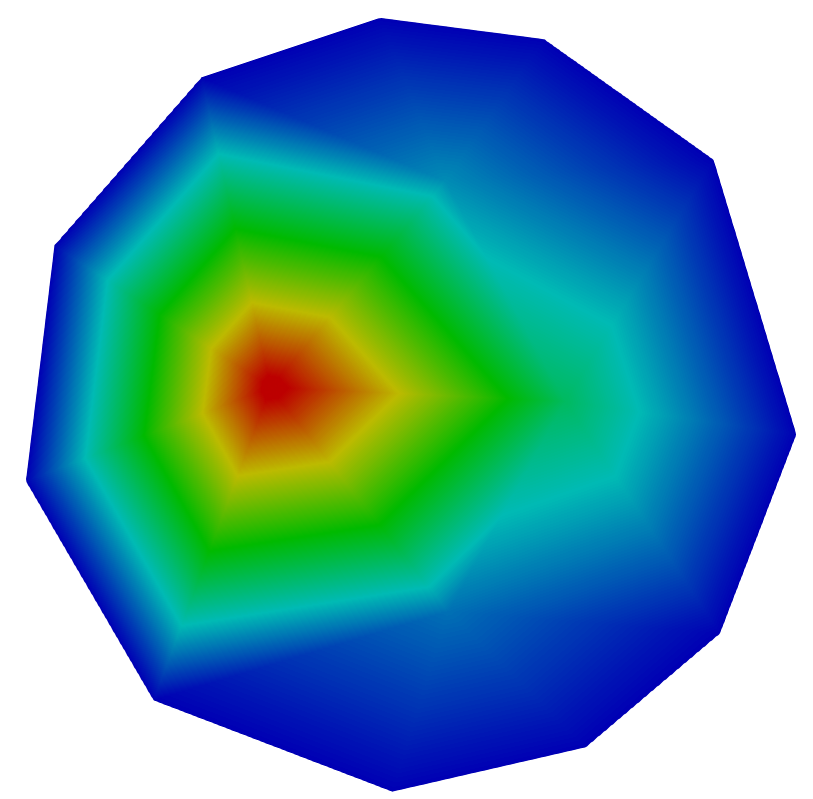}
    \caption{$t=0.6$ s}
    \label{fig11:cabg_u_t0.6}
  \end{subfigure}
   \hfill
  \begin{subfigure}[b]{0.15\textwidth}
    \centering
    \includegraphics[width=\textwidth]{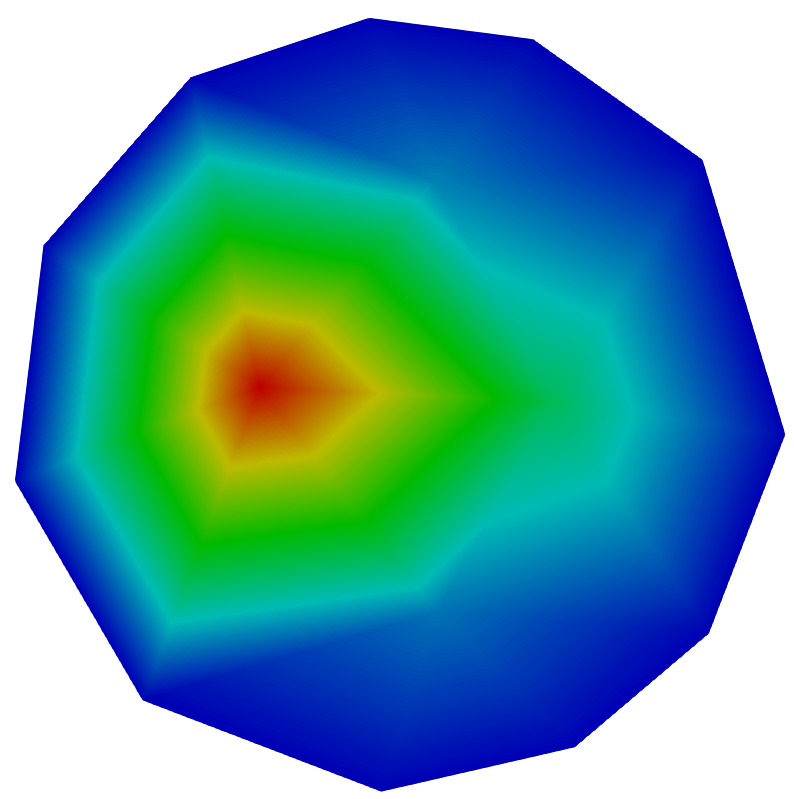}
    \caption{$t=0.8$ s}
    \label{fig11:cabg_u_t0.8}
  \end{subfigure}
   %\vskip\baselineskip  % Colorbar for WSS
  \begin{subfigure}[b]{0.3\textwidth}
    \centering
    \includegraphics[width=\textwidth]{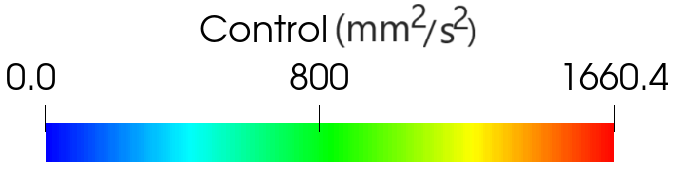}
    \label{fig11:colorbar_u_t_cabg}
  \end{subfigure}

\vskip\baselineskip
  % WSS distribution
  \begin{subfigure}[b]{0.19\textwidth}
    \centering
    \includegraphics[width=\textwidth]{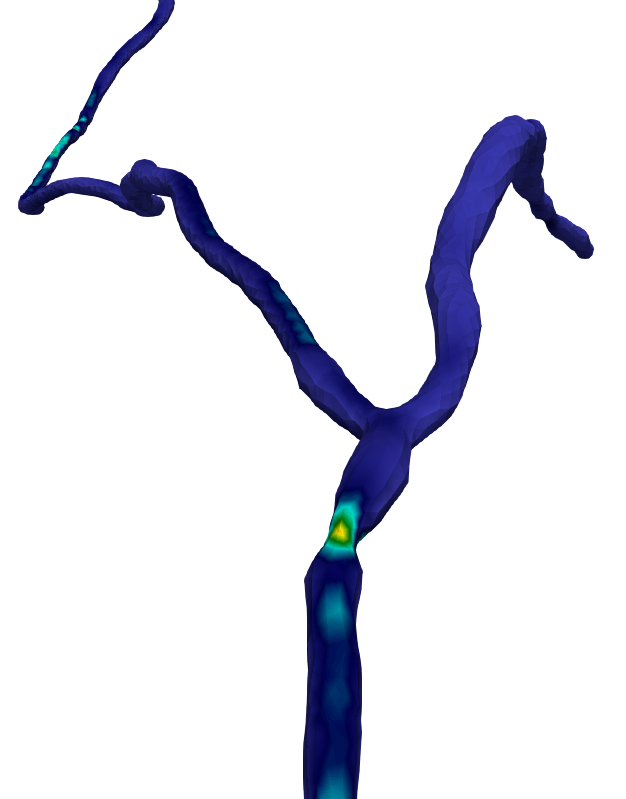}
    \caption{$t=0.04$ s}
    \label{fig11:cabg_wss_t0.04}
  \end{subfigure}
  \hfill
  \begin{subfigure}[b]{0.19\textwidth}
    \centering
    \includegraphics[width=\textwidth]{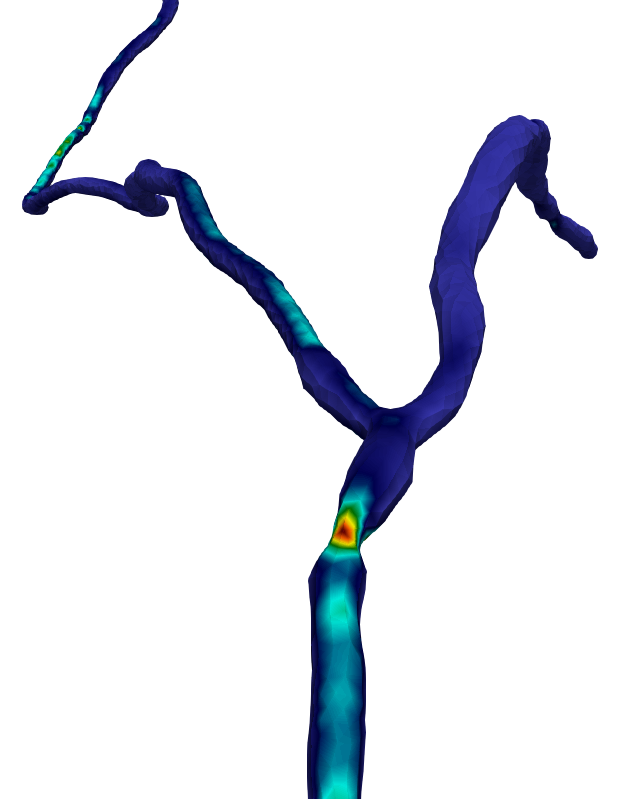}
    \caption{$t=0.2$ s}
    \label{fig11:cabg_wss_t0.2}
  \end{subfigure}
  \hfill
  \begin{subfigure}[b]{0.19\textwidth}
    \centering
    \includegraphics[width=\textwidth]{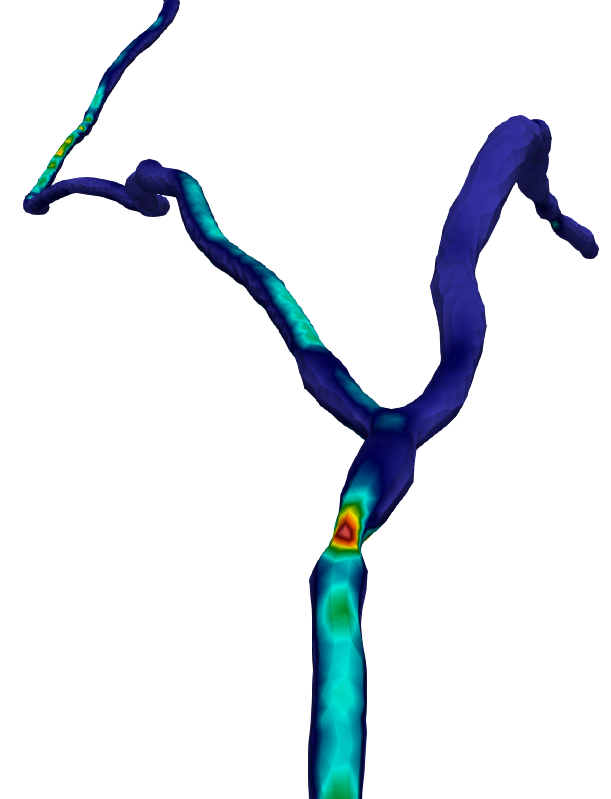}    
    \caption{$t=0.4$ s}
    \label{fig11:cabg_wss_t0.4}
  \end{subfigure}
  \hfill
  \begin{subfigure}[b]{0.19\textwidth}
    \centering
    \includegraphics[width=\textwidth]{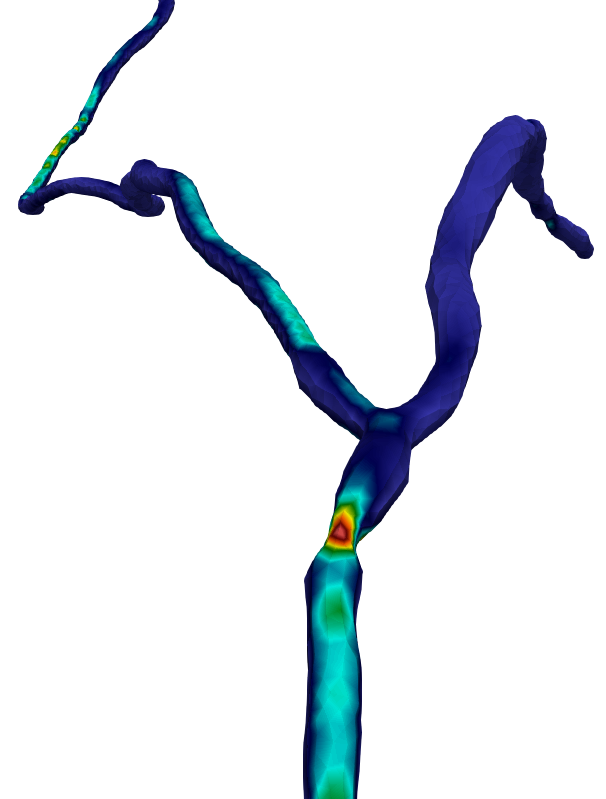}
    \caption{$t=0.6$ s}
    \label{fig11:cabg_wss_t0.6}
  \end{subfigure}
   \hfill
  \begin{subfigure}[b]{0.19\textwidth}
    \centering
    \includegraphics[width=\textwidth]{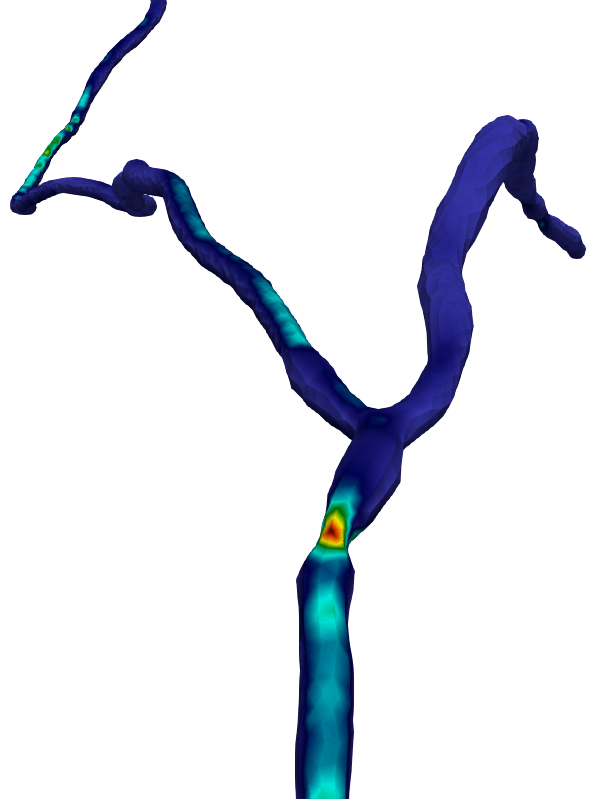}
    \caption{$t=0.8$ s}
    \label{fig11:cabg_wss_t0.8}
  \end{subfigure}
   
   \vskip\baselineskip  % Colorbar for WSS
  \begin{subfigure}[b]{0.3\textwidth}
    \centering
    \includegraphics[width=\textwidth]{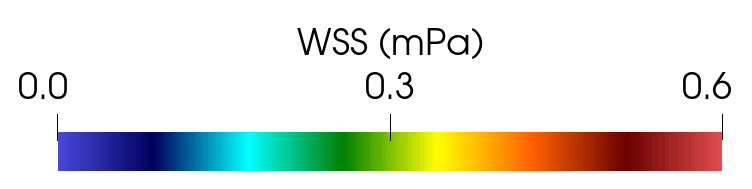}
    \label{fig11:colorbar_wss_t_cabg}
  \end{subfigure}
  
  \caption{Temporal variation with $\left(Re_1, Re_2\right) = \left(80, 50\right)$: Control distributions at $\Gamma_\mathrm{out}$ (top row) and  WSS distribution (bottom row).}
  \label{fig11:timevar_cabg}
\end{figure}

Figure~\ref{fig11:timevar_cabg} depicts the control distribution at the outflow boundary $\Gamma_\mathrm{out}$ and wall shear stress (WSS) distributions at different time instances with $\left(Re_1, Re_2 \right) = \left(80, 50 \right)$.  WSS is a crucial parameter in CV flows because it directly influences endothelial cell functions~\cite{ku1985pulsatile,gijsen2019expert,rathore2021}. It is defined as the tangential force exerted by blood flow on the endothelial surface of the vessels, expressed by:
\begin{equation}\label{eq:4.3}
\displaystyle  \bm{\tau}_w =  \nu\,\Big( \nabla \left(\bm{v}\right) + \nabla\left(\bm{v}\right)^\top \Big) \cdot \bm{n}.
\end{equation}
Figure~\ref{fig11:timevar_cabg} shows the flow dynamics as two inlet flows merge into a single outlet, significantly affecting the flow characteristics and control distributions. At time instance $t=0.04$ s, the control is distributed within the interior region of the outlet $\Gamma_\mathrm{out}$, avoiding the boundary edges, exerting a significant influence to regulate the early, developing flow dynamics. This control distribution influences the developing flow to conform to the desired velocity profile early in the cycle. Correspondingly, the WSS is relatively low throughout most of the vessel, and concentrated in regions of stenosis, where the flow merges and accelerates, leading to higher stress on the vessel walls. As time progresses, particularly by $t=0.2$ s and $t=0.4$ s, the control magnitude increases, which reflects a greater need for regulation as the inlet flows begin interacting more strongly and the system works to steer the flow toward the target profile. Despite this increase in control, the WSS distribution also rises, particularly near the bifurcations and graft anastomosis regions where the velocity gradients intensify.  By $t=0.8$ s, the control slightly reduces intensity, indicating a transition toward a more developed and regulated flow state; however, WSS remains concentrated within the same critical regions. WSS is directly proportional to velocity gradients at the wall; a steeper gradient leads to higher WSS due to the greater difference in velocity between the fluid near the wall and that further away. The region of high WSS corresponds to the high-velocity regions as depicted in Figure~\ref{fig9:cabg_fom_vel}. The observed temporal variations in control and WSS are due to the inter-dependency between the state, adjoint, and control solutions, Eqs.~(\ref{eq:2.3}), (\ref{eq:2.12}) and (\ref{eq:2.13}), which ensure that control effectively regulates flow dynamics with the objective and optimality conditions. 
\begin{figure}[htbp]
  \centering
  % First row of subfigures
  \begin{subfigure}[b]{0.312\textwidth}
    \centering
    \includegraphics[width=\textwidth]{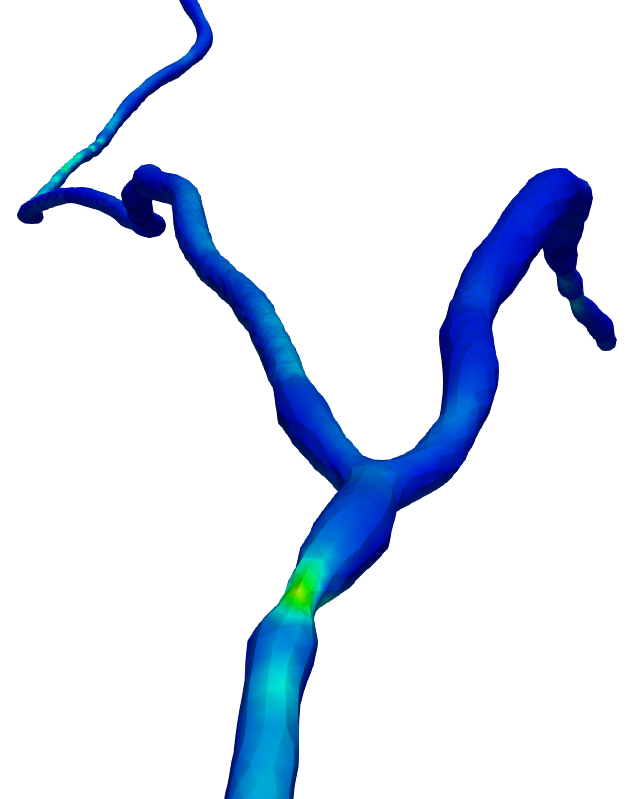}
    \caption{$\left( Re_1, Re_2 \right) = \left(50, 50 \right)$}
    \label{fig12:cabg_wss_Re_50_50}
  \end{subfigure}
  \hfill
  \begin{subfigure}[b]{0.312\textwidth}
    \centering
    \includegraphics[width=\textwidth]{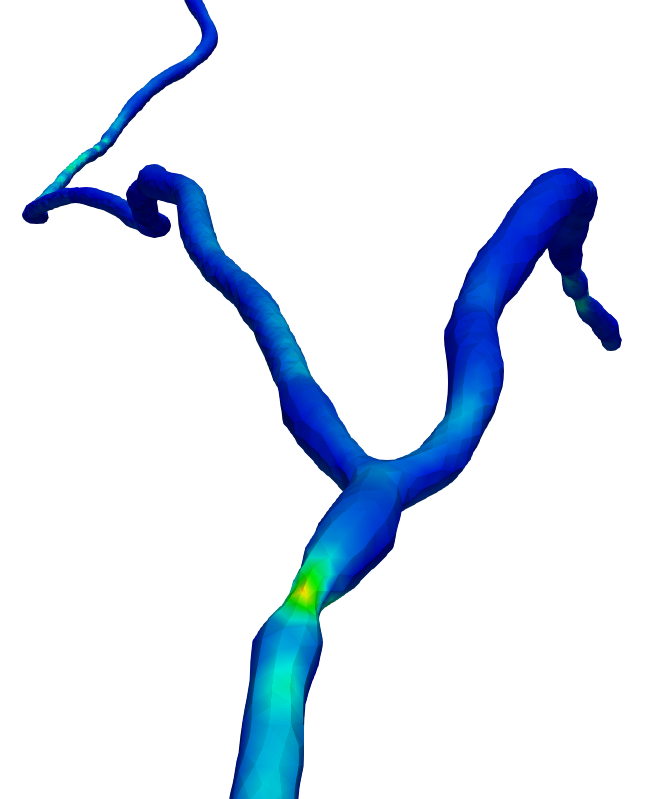}
    \caption{$  \left(Re_1, Re_2 \right) = \left(50, 70 \right)$}
    \label{fig12:cabg_wss_Re_50_70}
  \end{subfigure}
  \hfill
  \begin{subfigure}[b]{0.312\textwidth}
    \centering
    \includegraphics[width=\textwidth]{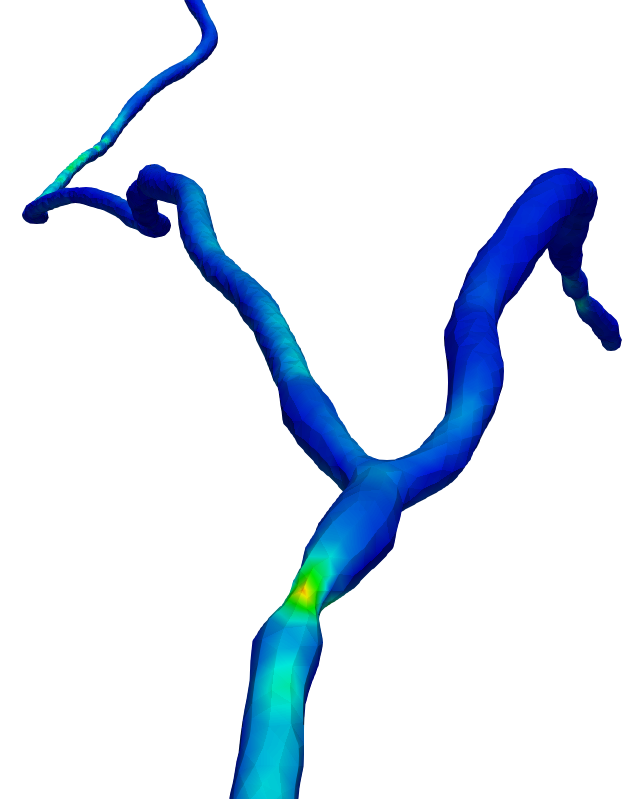}
    \caption{$\left( Re_1, Re_2 \right) = \left(60, 60 \right)$}
    \label{fig12:cabg_wss_Re_60_60}
  \end{subfigure}

  \vskip\baselineskip % Add vertical space between rows

  % Second row of subfigures
  \begin{subfigure}[b]{0.3\textwidth}
    \centering
    \includegraphics[width=\textwidth]{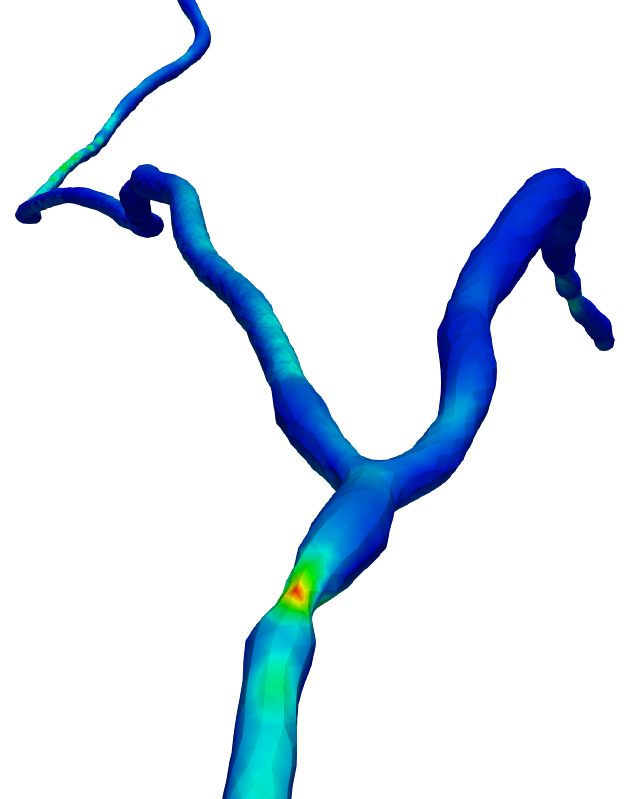}
    \caption{$ \left( Re_1, Re_2 \right) = \left(70, 70 \right)$}
    \label{fig12:cabg_wss_Re_70_70}
  \end{subfigure}
  \hfill
  \begin{subfigure}[b]{0.3\textwidth}
    \centering
    \includegraphics[width=\textwidth]{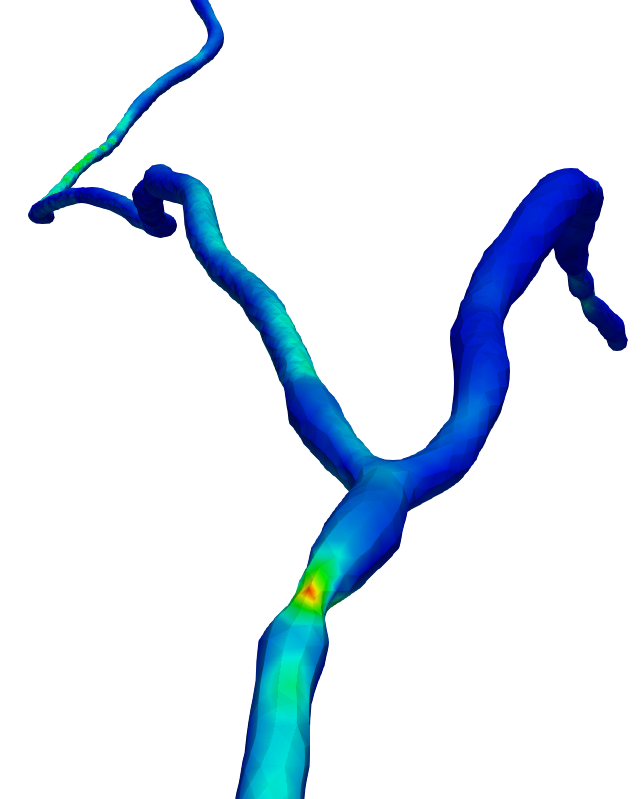}
    \caption{$\left( Re_1, Re_2 \right) = \left(80, 50 \right)$}
    \label{fig12:cabg_wss_Re_80_50}
  \end{subfigure}
  \hfill
  \begin{subfigure}[b]{0.3\textwidth}
    \centering
    \includegraphics[width=\textwidth]{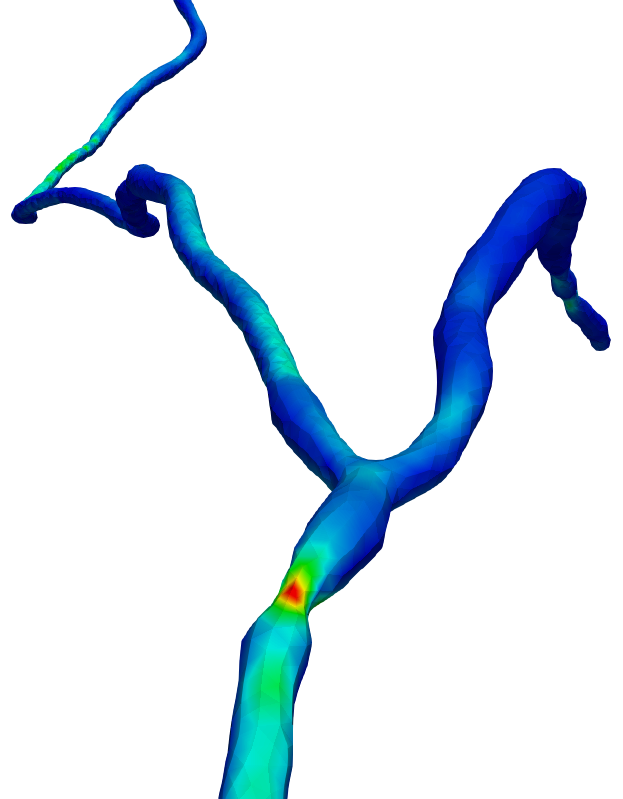}
    \caption{$\left( Re_1, Re_2 \right) = \left(80, 80 \right)$}
    \label{fig12:cabg_wss_Re_80_80}
  \end{subfigure}

 \vskip\baselineskip % Add vertical space between rows
 \begin{subfigure}[b]{0.5\textwidth}
    \centering
    \includegraphics[width=0.75\textwidth]{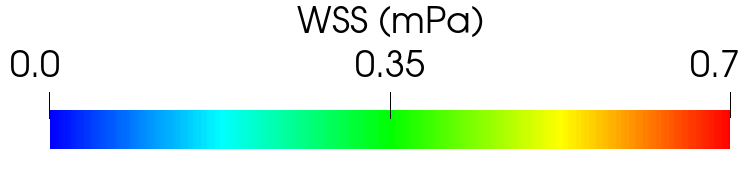}
    \label{fig12:colorbar_wss_Re_cabg}
  \end{subfigure}
  \caption{Parametric dependency of $\left( Re_1, Re_2 \right)$ values on reduced order WSS distributions at $t = 0.4$ s.}
  \label{fig12:Re_wss_cabg}
\end{figure}
In Figure~\ref{fig12:Re_wss_cabg}, the WSS distribution shows how the shear stress is affected by different combinations of Reynolds numbers. When the  $Re$ value is same at both inlets $\Gamma_\mathrm{in},$ as shown in Figures~\ref{fig12:cabg_wss_Re_50_50}, \ref{fig12:cabg_wss_Re_60_60} \ref{fig12:cabg_wss_Re_70_70} and \ref{fig12:cabg_wss_Re_80_80}, the WSS distribution is symmetric, and as it $Re$ increases, the WSS also increases, especially in critical regions such as bifurcations and stenosis. This is because higher $Re$ values lead to higher flow velocity, which increases the shear stress on the vessel walls. However, a notable difference is observed when comparing the WSS values from Figures~\ref{fig12:cabg_wss_Re_50_70} and \ref{fig12:cabg_wss_Re_80_50} and Figures~\ref{fig12:cabg_wss_Re_70_70} and \ref{fig12:cabg_wss_Re_80_50}, respectively. In Figure~\ref{fig12:cabg_wss_Re_50_70}, we have $\left(Re_1, Re_2\right) = \left(50, 70\right)$, and the WSS distribution is asymmetric due to the difference in $Re$ between the two inlets. The inlet with $Re=70$ experiences higher WSS than the inlet with $Re=50$. A similar pattern occurs for the Figure~\ref{fig12:cabg_wss_Re_80_50}; even though it has a higher $Re$ of $80$ in one inlet, the overall WSS is lower than in  Figure~\ref{fig12:cabg_wss_Re_70_70}, where both inlets have the same $Re$ values of $70.$  This occurs because the inlet flow with $Re=50$ in Figure~\ref{fig12:cabg_wss_Re_80_50} has a slower flow rate, which lowers the overall WSS. Such different flow conditions, where one inlet has a much lower $Re$, reduce the overall shear stress, despite the other inlet experiencing high stress due to a high Reynolds number. 

%%-------------------------------------------------------%%
\section{Discussion}%{Conclusion}
\label{sec5:Con}

This study successfully introduced a projection-based reduced approach for solving parameterized optimization problems related to the unsteady N-S equations, specifically focusing on computing CV flows. 
{The findings illustrated in Figures~\ref{fig1311:discrepencies} and~\ref{fig13:discrepencies} highlight that the effectiveness of outlet-based control is highly dependent on vascular geometry, showing considerably less influence within the idealised, symmetric models but demonstrating significant improvements in physiologically realistic, asymmetric configurations, \textit{i.e.} patient-specific CABG. These results underscore the importance of incorporating anatomical detail when formulating the control strategies for CV flow simulations.}  Additionally, they show minimal discrepancies between the targeted and computed velocity profiles derived from \eqref{eq:4.2}, demonstrating that we have achieved our objective and highlighting the robustness of our optimisation-based framework for CV flows. These simulations can be computationally intensive, especially for complex 3D geometries and patient-specific data. 
Cutting-edge techniques, including reduced modelling, machine learning, and deep learning, enable real-time hemodynamics by reducing computational costs and preserving accuracy for real-time CV applications \cite{rozza2005optimization,zainib2021,siena2023data,balzotti2024}. We observed that our proposed projection-based ROM effectively captures the essential dynamics of the flow fields while substantially reducing computational costs compared to high-fidelity solutions. The speed-up achieved by the ROM, which is approximately 8 times faster than a high-fidelity solution as shown in Table \ref{tab1:parameters} for both the test cases, underscores its advantage in accelerating simulations over traditional high-fidelity computations. Figures~\ref{fig2:eigen_bif} and \ref{fig8:eigen_cabg} illustrate that a $90\%$ reduction in snapshots using the nested-POD approach for temporal compression significantly enhances data management for OCP$_{(\bm{\mu})}$ problems while preserving accuracy; the same has been understood for the flow problems \cite{ballarin2017numerical,ballarin2016fast}. It also highlights that only a few POD modes are necessary to capture the crucial dynamic features for addressing CV flows. Figures~\ref{fig3:Vel_comp_bif}, ~\ref{fig4:press_comp_bif}, ~\ref{fig9:vel_comp_cabg} and ~\ref{fig10:press_comp_cabg}, demonstrate that this reduced methodology effectively reproduces the detailed features observed in the Galerkin FE formulation for both velocity and pressure distributions in vascular models at specified Reynolds number and time.

The configured outflow boundary condition ensures that the control strategy can dynamically adjust to optimize blood flow essential for maintaining stable and efficient flow within the computational domain, as shown in Figures~\ref{fig5:timevar_bif} and \ref{fig11:timevar_cabg}. These figures illustrate how variations in control distributions significantly impact flow dynamics within a system where two inlet streams merge into a single outlet, altering both velocity profiles and wall shear stress distributions.  Furthermore, it underscores the interdependence among the state, adjoint, and control solutions and elucidates how control regulates the flow in both test cases. This analysis further confirms that WSS and control distributions maintain consistency with the desired flow dynamics while the outlet pressure decreases from Figure~\ref{fig7:Re_press_bif}. These findings align with research highlighting the role of OCPs in enhancing hemodynamics, patient-specific simulations, and clinical decisions \cite{fevola2021, imperiale2023flow, zainib2021, balzotti2022data}. Given the critical role of the outflow boundary in CV flows \cite{vignon2010outflow,rathore2021,africa2024lifex}, our control strategy is crucial for capturing detailed flow dynamics and ensuring effective flow regulation. These findings highlight the importance of precise control in maintaining efficient hemodynamics and realistic physiological behaviour. Figure \ref{fig12:Re_wss_cabg} shows that WSS depends on whether the $Re$ values are the same or different between the inlets, with equal $Re$ leading to higher and more consistent WSS, while different $Re$ values can lower WSS in one inlet despite high stress in the other. In clinical practice, high WSS is associated with increased risks of vascular damage, atherosclerosis, or remodelling, particularly in bifurcation and stenotic regions with high Reynolds number flows \cite{ku1985pulsatile,gijsen2019expert}. This ROM approach demonstrates a robust and efficient methodology for simulating complex CV flows, providing significant computational savings without compromising accuracy. Its application to idealized and patient-specific CABG models highlights its potential for real-time simulations and many-query cardiovascular studies, making it a valuable tool for personalized modelling and improving clinical outcomes.

\section*{Declaration of generative AI and AI-assisted technologies in the writing process}
 During the preparation of this work, the authors utilised ChatGPT 4.0 and Grammarly to ensure consistency and check grammar. After using these tools, the authors reviewed and edited the content as needed and accepted full responsibility for the content of the publication.  
 
\section*{Competing interests}
We hereby declare that all the authors have no potential conflicts of interest. 

\section*{Ethics Statement}
The authors declare that no human or animal subjects were involved in this research, and as such, no ethical approval or informed consent are required.

\section*{Author contributions}
\textbf{SR:} Conceptualization, Methodology, Software, Formal Analysis, Investigation, Visualization, Writing-Original Draft, Review \& Editing;
\textbf{PCA:} Conceptualization, Methodology, Formal Analysis, Writing-Review \& Editing, Supervision;
\textbf{FB:} Conceptualization, Methodology, Software, Formal Analysis, Writing-Review \& Editing;
\textbf{FP:} Conceptualization, Methodology, Formal Analysis, Writing-Review \& Editing;
\textbf{MG:} Formal Analysis, Writing-Review \& Editing; 
\textbf{GR:} Funding Acquisition, Project Administration, Writing-Review \& Editing. 

\section*{Funding}

\textbf{SR} and \textbf{GR} acknowledge the financial support provided by the project ``Advanced developments for scientific computing in complex systems with applications in engineering, medicine and environmental sciences,'' funded by the Ente finanziatore: MUR; Canale di finanziamento: PRO3 (CUP  G95F21001980006). 
\textbf{PCA}, \textbf{FP} and \textbf{GR} acknowledge the support provided by the European Union - NextGenerationEU, in the framework of the iNEST - Interconnected Nord-Est Innovation Ecosystem (iNEST ECS00000043 – CUP G93C22000610007) consortium and its CC5 Young Researchers initiative.
\textbf{PCA} also acknowledges the INdAM-GNCS Project ``Metodi numerici efficienti per problemi accoppiati in sistemi complessi'' (CUP E53C24001950001).
\textbf{FB} acknowledges the European Union's Horizon 2020 research and innovation program under the Marie Sk\l{}odowska-Curie Actions, grant agreement 872442 (ARIA).
\textbf{FB} also acknowledges the PRIN 2022 Project 202249PF73 ``Mathematical models for viscoelastic biological matter'' funded by the European Union -- NextGenerationEU under the National Recovery and Resilience Plan (NRRP), Mission 4 Component 2 Investment 1.1 - Call PRIN 2022 No. 104 of February 2, 2022 of Italian Ministry of University and Research (CUP J53D23003590008). 
\textbf{MG} and \textbf{GR} acknowledge the PRIN 2022 Project ``Machine learning for fluid-structure interaction in
cardiovascular problems: efficient solutions, model reduction, inverse problem (CUP G53D23006830001)''. The authors would also like to acknowledge INdAM-GNSC for its support.

%%=================================================================%%
 \bibliographystyle{elsarticle-num} 
 \bibliography{reference}
\end{document}